\documentclass[11pt,english]{article}
\usepackage[english]{babel}
\usepackage{amsfonts}
\usepackage{epsfig}
\usepackage{graphicx}
\usepackage{fancyhdr}
\usepackage{amssymb}
\usepackage{amsmath}
 \usepackage{verbatim} 
 \usepackage{appendix} 
\usepackage{xcolor}

\newtheorem{theorem}{Theorem}

\newtheorem{corollary}{Corollary}[section]

\newtheorem{prop}{Proposition}[section]
\newtheorem{lemma}[prop]{Lemma}
\newtheorem{definition}[prop]{Definition}
\newtheorem{cl}[prop]{Claim}
\newtheorem{rem}[prop]{Remark}

\def\div{{\rm div}\,}

\def\Ag {{\cal A}} 
\def\Bg {{\cal B}} 
\def\Cg {{\cal C}} 
\def\Dg {{\cal D}} 
\def\Fg {{\cal F}} 
\def\Lg {{\cal L}} 
\def\Qg {{\cal Q}} 
\def\Rg {{\cal R}} 
\def\Vg {{\cal V}} 

\def\and {{\rm \; and \;}}

\def\exp {{\rm exp}}

\def\Ag {{\cal A}} 
\def\Bg {{\cal B}} 
\def\Cg {{\cal C}} 
\def\Dg {{\cal D}} 
\def\Fg {{\cal F}} 
\def\Lg {{\cal L}} 
\def\Qg {{\cal Q}} 

\def\Vg {{\cal V}} 

\def\and {{\rm \; and \;}}

\def\exp {{\rm exp}}

\def\T0{T_{0,1}}

\newcommand {\R}{ \mathbb{R}}
\newcommand {\C}{ \mathbb{C}}
\newcommand {\N}{ \mathbb{N}}

\newcommand {\pa}{\partial}

\newcommand {\beqna} {\begin{eqnarray}}
\newcommand {\eeqna} {\end{eqnarray}}
\newcommand {\beqtn} {\begin{equation}}
\newcommand {\eeqtn} {\end{equation}}

\newcommand {\dsp}{\displaystyle}

\begin{document}

\title{\bf Construction of  blowup solutions  for the Complex Ginzburg-Landau equation with critical parameters}
\author{Giao Ky Duong, \\ {\it \small gkd1@nyu.edu}\\ {\it \small Department of Mathematics, New York University in Abu Dhabi.}\\
Nejla Nouaili,\\ {\it \small Nejla.Nouaili@dauphine.fr}\\ {\it \small CEREMADE, Universit\'e Paris Dauphine, Paris Sciences et Lettres.} \\
Hatem Zaag,\\
{\it \small Hatem.Zaag@univ-paris13.fr}\\{\it \small Universit\'e Paris 13,}\\ {\it \small LAGA, CNRS (UMR 7539), F-93430, Villetaneuse, France.}}

\maketitle
\begin{abstract}

We construct a solution for the Complex Ginzburg-Landau (CGL) equation in a general  critical case, which blows up in finite time $T$ only at one blow-up point. We also give a sharp description of its profile. In a first part, we construct formally a blow-up solution. In a second part we give the rigorous proof. The proof relies on the reduction of the problem to a finite dimensional one, and the use of index theory to conclude. The interpretation of the parameters of the finite dimension problem in terms of the blow-up point and time allows to prove the stability of the constructed solution. We would like to mention that the asymptotic profile of our solution is different from previously known profiles for CGL or for the semilinear heat equation.

\end{abstract}

\section{Introduction}

We consider the following Complex Ginzburg-Landau (CGL) equation
\beqtn
\begin{array}{l}
u_t=(1+i\beta)\Delta u+(1+i\delta) |u|^{p-1}u+\alpha u,\\
u(.,0)=u_0\in L^\infty (\R^N,\C)
\end{array}
\mbox{     (CGL)}
\label{GL}
\eeqtn
where $\delta,\beta,\alpha\in\R$.

\medskip

This equation, most often considered with a cubic nonlinearity ($p=3$), has a long history in physics (see Aranson and Kramer \cite{AKRMP02}). The Complex Ginzburg-Landau (CGL) equation  is one of the most studied equations in physics. It describes a lot of phenomena including nonlinear waves, second-order phase transitions, and superconductivity. We note that the Ginzburg-Landau equation can be used to describe the evolution of amplitudes of unstable modes for any process exhibiting a Hopf bifurcation (see for example Section VI-C, page 37 and Section VII, page 40 from \cite{AKRMP02} and the references cited therein). The equation can be considered as a general normal form for a large class of bifurcations and nonlinear wave phenomena in continuous media systems. More generally, the Complex Ginzburg-Landau (CGL) equation is used to describe synchronization and collective oscillation in complex media.

\medskip

The study of blow-up, collapse or chaotic solutions of equation \eqref{GL} appears in many works; in the description of an unstable plane Poiseuille flow, see Stewartson and Stuart \cite{SSJFM71}, Hocking, Stewartson, Stuart and Brown \cite{HSSBJFM72} or in the context of binary mixtures in Kolodner and \textit{al}, \cite{KBSPRL88}, \cite{KSALPNP95}, where the authors describe an extensive series of experiments on traveling-wave convection in an ethanol/water mixture, and they observe collapse solution that appear experimentally.\\
For our purpose, we consider CGL independently from any particular physical context and investigate it as a mathematical model in partial differential equations with $p>1$.
 
\medskip

We note also that the interest on the study of singular solutions in CGL comes also from the analogies with the three-dimensional Navier-Stokes. The two equations have the same scaling properties and the same energy identity (for more details see the work of Plech{\'a}{\v{c}} and {\v{S}}ver{\'a}k \cite{PSCPAM01}; the authors in this work give some evidence for the existence of a radial solution which blow up in a self-similar way). Their argument is based on matching a numerical solution in an inner region with an analytical solution in an outer region.  In the same direction we can also  cite the work of  Rottsch\"afer \cite{RPD08} and \cite{REJAM13}.

\bigskip

\medskip

The Cauchy problem for equation \eqref{GL} can be solved in a variety of spaces using the semigroup theory as in the case of the heat equation (see \cite{CNYUCIM03,GV96,GVCMP97}).\\
We say that $u(t)$ blows up or collapse in finite time $T<\infty$, if $u(t)$ exists for all $t\in [0,T)$ and 
$\lim_{t\to T}\|u(t)\|_{L^\infty}=+\infty.$
In that case, $T$ is called the blow-up time of the solution. A point $x_0\in\R^N$ is said to be a blow-up point if there is a sequence $\{(x_j,t_j)\}$, such that $x_j\to x_0$, $t_j\to T$ and $|u(x_j,t_j)|\to \infty$ as $j\to\infty$. The set of all blow-up points is called the blow-up set.\\
Let us now introduce the following definition;
\begin{definition}\label{defini-cris-condition} The parameters $(\beta,\delta)$ are said to be critical  (resp. subcritical, resp. supercritical) if 
$p-\delta^2-\beta\delta (p+1)=0$ (resp. $>0$, resp. $<0$). In addition to that, we also define some critical constants  as follows:
\begin{eqnarray}
p_{cri} & =& \left\{\begin{array}{rcl}
 \sqrt{\frac{p (2p - 1)}{p-2}}  && \text{ if }  p> 2\\
 +\infty  && \text{ if } p \in (1,2]
\end{array} \right. , \label{defi-p-cri}
\end{eqnarray}
and 
\begin{equation}\label{defi-b-cri}
b_{cri}^2  =  \frac{(p-1)^4 (p+1)^2 \delta^2}{  16 (1 + \delta^2) (p(2p-1) - (p-2)\delta^2)( (p+3)\delta^2 + p(3p+1)) } > 0,
\end{equation}
for all $ \delta \in (-p_{cri}, p_{cri} )$.

\label{subcond}
\end{definition}
\begin{rem}
We choose $p_{cri}$ in such way the denominator in expression \eqref{defi-b-cri} is strictly positif. 
\end{rem}
An extensive literature is devoted to the blow-up profiles for CGL when $\beta=\delta=0$ (which is the nonlinear heat equation),  see Vel{\'a}zquez \cite{VCPDE92, VTAMS93,VINDIANA93} and Zaag \cite{ZIHP02, ZCMP02, ZMME02} for partial results). In one space dimension, given $a$ a blow-up point, this is the situation:

\medskip

\begin{itemize}

\item either
\beqtn
\sup_{|x-a|\leq K\sqrt{(T-t)\log(T-t)}}\left|(T-t)^{\frac{1}{p-1}}u(x,t)-f\left(\frac{x-a}{\sqrt{(T-t)|\log(T-t)}}\right)\right|\to 0,
\label{u2}
\eeqtn
\item or for some $m\in \N$, $m\ge 2$, and $C_m>0$
 \beqtn
\sup_{|x-a|<K(T-t)^{1/2m}}\left|(T-t)^{\frac{1}{p-1}}u(x,t)-f_m\left(\frac{C_m (x-a)}{(T-t)^{1/2m}}\right)\right|\to 0,
\label{um}
\eeqtn
as $t\to T$, for any $K>0$, where 
\begin{equation}
\label{definitionf}
\begin{array}{l}
f(z)=\left(p-1+b_0z^2\right)^{-\frac{1}{p-1}}
\mbox{ with }b_0=\frac{(p-1)^2}{4p},\\
f_m(z) =\left(p-1+|z|^{2m}\right)^{-\frac{1}{p-1}}.
\end{array}
\end{equation}
\end{itemize}

\medskip

If $(\beta,\delta)\not =( 0,0)$, some results are available in the  \textit{subcritical} case by Zaag \cite{ZAIHPANL98} ($\beta=0$) and Masmoudi and Zaag \cite{MZ07} ($\beta\not = 0$). More precisely, if
\[p-\delta^2-\beta\delta (p+1)>0,\] 
then, the authors construct a solution of equation \eqref{GL}, which blows up in finite time $T>0$ only at the origin such that for all $t\in [0,T)$,
\begin{eqnarray}
\dsp\left\| (T-t)^{\frac{1+i\delta}{p-1}}\left |\log(T-t)\right |^{-i\mu}u(x,t)-\left(p-1+\frac{b_{sub}|x|^2}{(T-t)|\log(T-t)|}\right)^{-\frac{1+i\delta}{p-1}}
\right\|_{L^\infty} & &\label{profilesubc} \\
& & \hspace{-4.5cm}\leq \frac{C_0}{1+\sqrt{|\log(T-t)|}}, \nonumber
\end{eqnarray}
where
\beqtn
\dsp b_{sub}=\frac{(p-1)^2}{4(p-\delta^2-\beta\delta(1+ p))}>0\mbox{ and }\mu=-\frac{2b_{sub}\beta}{(p-1)^2}(1+\delta^2).
\label{bs}
\eeqtn
Note that this result was previously obtained formally by Hocking and Stewartson \cite{HSPRSLS72} ($p=3$) and mentioned later in Popp et \textit{al} \cite{PSKKPD98} (see those references for more blow-up results often aproved numerically, in various regimes of the parameters).

\medskip


\medskip

In the  \textit{critical} case, few results are known about blow-up solutions for the equation. Up to our knowledge, there are two blow-up results in the critical case: one formal by Popp and \textit{al} \cite{PSKKPD98}, when $p=3$ and the result by the second and thirs author in \cite{NZARMA18} when $\beta=0$.

\medskip

Let us now give the formal result, when $p=3$, given by Popp and \textit{al} \cite{PSKKPD98} (equations (44) and (64)) in the  \textit{critical} case ($3-4\delta\beta-\delta^2=0$): the authors obtained 
\beqtn
(T-t)^{\frac{1+i\delta}{p-1}}u(x,t)\sim e^{i\psi(t)}  \left( 2+ \frac{b_{PSKK}|x|^2}{(T-t)|\log(T-t)|^{\frac 12}}\right)^{-\frac{1+i\delta}{2}}.
\label{propal}
\eeqtn
where 
\beqtn
b_{PSKK}=2\left( \sqrt{\frac {3}{2 \delta^2}(\delta^2+5)(\delta^2+1)(15-\delta^2)}\right)^{-1}.
\label{defbp}
\eeqtn
and $\psi(t)$ is given by equation (40) in \cite{PSKKPD98}. We can clearly see that this profile exist only for $\delta^2<15$, when $p=3$.
\begin{rem}
We will see later, in Section \ref{sectFormapp}, that we obtain the formally the same $b_{PSKK}$, for any $p>1$. Moreover, the constant in   \eqref{defbp}    will be proven to be  true    in the rigorous proof. 
\end{rem}
As for the result of \cite{NZARMA18}, we simply say that as far as statements are concerned, it is a particular case of our new result given in Theorem \ref{Theorem-existence} below (just take $\beta=0$). However, we stress the fact that the profs of $\beta=0$ and $\beta\not =0$ are very different technically.

\medskip
\subsection{Statement of our result}
Our main claim is to construct a solution $u(x,t)$ of \eqref{GL} in the critical case ($\beta\neq0$ and $p-\beta\delta(p+1)-\delta^2=0$) that blows up in some finite time $T$, in the sense that
\[\dsp \lim_{t\to T}\|u(.,t)\|_{L^\infty}=+\infty.\]
We also prove the stability of the constructed solution. This is our first statement:
 
\begin{theorem}[Blow-up profiles for equation \eqref{GL}] \label{Theorem-existence}Let us consider the critical case where  $p-\delta^2-\beta\delta(p+1)=0$, $\beta\neq 0$ and $\delta \in (0,p_{cri})$ with $p_{cri}$ defined as in \eqref{defi-p-cri}. Then, there exists a unique constant $\mu$ depending on $p$ and $\delta$ such that equation \eqref{GL} has a solution $u(x,t)$, which blows up in finite time $T$, only at the origin. Moreover:

(i) For all $t\in [0,T)$,

\begin{eqnarray}
\left\|(T-t)^{\frac{1+i\delta}{p-1}}e^{-i\nu\sqrt{|\log(T-t)|}}|\log (T-t)|^{-i\mu} u(.,t)-\varphi_0\left(\frac{.}{\sqrt{(T-t)}|\log (T-t)|^{1/4}}\right) \right\|_{L^\infty(\R^N)} & &\label{profile-u} \\
& & \hspace{-4.5cm}  \leq \frac{C_0}{1+|\log(T-t)|^{\frac 14}},\nonumber
\end{eqnarray}
where
\beqtn
\varphi_0(z)=\left ( p-1+b_{cri} z^2\right )^{-\frac{1+i\delta}{p-1}},
\label{deffi0}
\eeqtn
 with $b_{cri}$ is defined as in \eqref{defi-b-cri},\;\;

\beqtn
  \nu =   -\frac{ 4 b \beta (1 + \delta^2)}{(p-1)^2} \text{ and  }   a=2\kappa (1-\beta\delta)\frac{b}{(p-1)^2}.
\label{definition-nu-a}
\eeqtn

(ii) For all $x\not =0$, $u(x,t) \to u^*(x)\in C^2(\R^N\backslash\{0\})$ and
\beqtn
\dsp u^{*}(x) \sim e^{i\nu \sqrt{2|\log |x||}}|2\log |x||^{i\mu_{cri}}\left [ \frac{b_{cri}|x|^2}{\sqrt{2|\log |x||}}\right]^{-\frac{1+i\delta}{p-1}}\mbox { as }x\to 0.
\label{finalprofile}
\eeqtn 
\label{thm1}
\end{theorem}

\medskip

\begin{theorem}[First order terms] \label{Theorem-existence2}
Following Theorem \ref{Theorem-existence}, we claim that the solution decomposes in self similar variables
\[W(y,t)=(T-t)^{\frac{1+i\delta}{p-1}}u(x,t),\;\;y=\frac{x}{\sqrt{T-t}},\]
 as follows: For $M>0$
  \beqtn
 \begin{array}{ll}
\dsp \sup_{|y|<M|\log(T-t)|^{\frac 14}}\left |W(y,t)e^{-i\nu\sqrt{|\log(T-t)|}}|\log(T-t)|^{-i\mu}e^{i\theta(t)}\right .&\\
-\left .\left \{ \dsp \varphi_0\left (\frac{y}{|\log (T-t)|^{1/4}}\right )
\dsp+\frac{a(1+i\delta)}{|\log(T-t)|^{\frac 12}}
 \dsp + \frac{1}{|\log(T-t)|}\Fg(y)\right \}\right |&\\
 & \hspace{-0.8cm} \leq \dsp \frac{C}{|\log(T-t)|^{\frac 32}}(1+|y|^5),
 \end{array}
 \eeqtn
 and $\theta(t)\to \theta_0$ as $t\to T$, such that
 \[|\theta(t)-\theta_0|\leq \frac{C}{|\log(T-t)|^{\frac 14}}\]
 \noindent
 with
 \beqtn
\varphi_0(z)=\left ( p-1+b_{cri} z^2\right )^{-\frac{1+i\delta}{p-1}},
\eeqtn
where $b_{cri}$ is defined as in \eqref{defi-b-cri}, $\nu$ and $a$ are given by \eqref{definition-nu-a} and $\Fg(y)$ is a function defined as follows
 
 \beqtn
\Fg(y)= \Ag_0 h_0(y)+ \Ag_2 h_2(y)+\tilde \Ag_2\tilde h_2(y), 
 \label{defFg}
 \eeqtn
where $ \Ag_0$, $\Ag_2$ and $\tilde \Ag_2$ depend only on $\beta$ and $\delta$ and are given by \eqref{defAg} in Definition \ref{defthess} and $ h_0(y)$, $h_2(y)$ and $\tilde h_2(y)$ will be given in Lemma \ref{lemma-formula-h-j-tilde-h-j}.
 \label{Newthm}
 \end{theorem}
 \begin{rem} Theorem \ref{Theorem-existence2} is true only on the case of $\beta\neq 0$. 
 \end{rem}
 \begin{rem} 
 We will  give the proof only  in one dimension.   Our proof remains valid in higher dimensions, with exactly the same ideas, and purely technical differences, that we omit to keep this (already long) paper in reasonable length.    Indeed, the computation of the eigenfunctions of  $\Lg_{\beta,\delta}$ (see \eqref{eqqd}) and the projection of  equation \eqref{eqq} on the eigenspaces become much more complicated when $N \geq 2$. 
 \end{rem}

 \begin{rem}
In this paper we will treat  the case $\beta\neq 0$, which is far from being a simple adaptation of the case $\beta=0$, treated in \cite{NZARMA18}. Indeed, the techniques are different. There is an additional difficulty arising from the linearized operator $\Lg_{\beta,\delta}$ (see \eqref{eqqd}), which is not diagonalisable as in the case $\beta=0$. To avoid this problem, we will use ideas from the work of Masmoudi and Zaag   \cite{MZ07}. There is another difficulty, coming from the criticality of the problem, which makes the projection of the linearized equation quite complicated especially for the neutral mode.
\end{rem}

\begin{rem} We will consider CGL, given by \eqref{GL}, only when $\alpha=0$. The case $\alpha\not =0$ can be done as in \cite{EZSMJ11}. In fact, when $\alpha\not =0$, exponentially small terms will be added to our estimates in self-similar variable (see \eqref{chauto} below), and that will be absorbed in our error terms, since our trap $\Vg_A(s)$ defined in Definition \ref{defthess} is given in polynomial scales.
\label{nu0}
\end{rem}

\begin{rem} 

The derivation of the blow-up profile \eqref{deffi0} can be understood through a formal analysis, using  the matching asymptotic expansions. 
This method was used by  Galaktionov, Herrero and Vel\'azquez \cite{VGH91} to derive all the possible behaviors of the blow-up solution given by (\ref{u2},\ref{um}) in the heat equation ($\beta=\delta=0$). This formal method fails in the determination of $b_{cri}$ when $\delta\neq 0$, because of the complexity of the system of ODE in that case. 

\medskip

We will use in this case the formal method used by Popp and al in \cite{PSKKPD98}. We note that this method was used by Hocking and Stewartson \cite{HSPRSLS72} as well as Masmoudi and Zaag \cite{MZ07} ($\delta\not = 0$) to obtain the profile in the subcritical case of CGL, and also by Berger and Kohn for the nonlinear heat equation ($\beta=\delta=0$).

\end{rem}
\begin{rem}
The exhibited profile  \eqref{deffi0} is new in two respects: 
\begin{itemize}
\item The scaling law in the critical case is $\sqrt{(T-t) |\log (T-t)|^{\frac 12}}$ instead of the laws of subcritical case, $\sqrt{(T-t)| \log (T-t)|}$, (see \eqref{profilesubc}). 
\item The profile function: $\varphi_0(z)=(p-1+b_{cri} |z|^2)^{-\frac{1+i\delta}{p-1}}$ is different from the profile of the subcritical case , namely $f(z)=(p-1+b_{sub} |z|^2)^{-\frac{1+i\delta}{p-1}}$, in the sense that  $b_{cri}\not = b_{sub}$ (see \eqref{bs} and  \eqref{defi-b-cri}).
\end{itemize}
\end{rem}
\begin{rem}
In the subcritical case $p-\delta^2-\beta\delta (p+1)>0 $ ($p-\delta^2>0$, when $\beta=0$), the final profile of the CGL is given by 

\[u^*(x)\sim |2\log |x||^{i\mu}\left [ \frac b 2 \frac{|x|^2}{|\log|x||}\right]^{-\frac{1+i\delta}{p-1}}\mbox{ as $x\to 0$  }\]
with
\[b_{sub}=\frac{(p-1)^2}{4(p-\delta^2-\beta\delta(p+1))}\mbox{ and }\mu=-\frac{2b\beta}{(p-1)^2}(1+\delta^2).\]
In the critical case $p-\delta^2-(p+1)\beta\delta=0$, the final profile is given by \eqref{finalprofile}.
\end{rem}

As a consequence of our techniques, we show the stability of the constructed solution with respect to perturbations in initial data. More precisely, we have the following result.

\begin{theorem}[Stability of the solution constructed in Theorem \ref{thm1}] Let us denote by $\hat u(x,t)$ the solution constructed in Theorem 1 and by $\hat T$ its blow-up time. Then, there exists a neighborhood $\Vg_0$ of $\hat u(x,0)$ in $L^\infty$ such that for any $u_0\in\Vg_0$, equation \eqref{GL} has a unique solution $u(x,t)$ with initial data $u_0$, and $u(x,t)$ blows up in finite time $T(u_0)$ at one single blow-up point $a(u_0)$. Moreover estimate \eqref{profile-u} is satisfied by $u(x-a,t)$ and
\[T(u_0)\to \hat T,\;\; a(u_0)\to 0\mbox{ as }u_0\to \hat u_0\mbox{ in }L^\infty(\R^N,\C).\]
\label{thstablity}
\end{theorem}

\begin{rem}
We will not give the proof of Theorem \ref{thstablity}
because the stability result follows from the reduction to a finite dimensional case as in \cite{MZDuke97}  (see Theorem 2 and its proof in Section 4) 
and \cite{MZ07} (see Theorem 2 and its proof in Section 6) with the same argument. Hence, we only prove the existence result (Theorem \ref{thm1}) and kindly refer the reader to \cite{MZDuke97} and \cite{MZ07} for the proof of the stability.
\end{rem}

\medskip

Let us give an idea of the method used to prove the results. We construct the blow-up solution with the profile in Theorem \ref{thm1}, by following the method of \cite{MZDuke97}, \cite{BKN94}, though we are far from a simple adaptation, since we are studying the critical problem, which make the technical details harder to elaborate. This kind of methods has been applied for various nonlinear evolution equations. For hyperbolic equations, it has been successfully used for the construction of multi-solitons for semilinear wave equation in one space dimension (see \cite{CZCPAM13}). For parabolic equations, it has been used in \cite{MZ07} and \cite{ZCPAM01} for the Complex Ginzburg Landau (CGL) equation with no gradient structure, the critical harmonic heat flow in  \cite{RSCPAM13}, the two dimensional Keller-Segel equation in \cite{RSMA14} and the nonlinear heat equation involving nonlinear gradient term in \cite{EZSMJ11}, \cite{TZ17}. Recently, this method has been applied for various non variational parabolic system in  \cite{NZCPDE15} and \cite{GNZJDE17,GNZIHP18,GNZJDE18,GNZAPAM19}, for a logarithmically perturbed nonlinear equation in  \cite{NVTZ15,DJFA19,DJDE19,DNZTJM19}.
We also mention a result for a higher order parabolic equation \cite{GNZANA20}, two more results for equation involving non local terms in  \cite{DZMMMAS19,AZ19}.
\medskip

Unlike in the subcritical case of \cite{MZ07} and \cite{ZCPAM01}, the criticality of the problem induces substantial changes in the blow-up profile as pointed-out in the comments following Theorem \ref{Theorem-existence}. Accordingly, its control requires special arguments. So, working in the framework of \cite{MZDuke97}, \cite{MZ07}, \cite{NZARMA18}, some crucial modifications are needed. In particular, we have overcome the following challenges:
\begin{itemize}
\item The prescribed profile was not known before and is not obvious to obtain. See Section \ref{sectFormapp} for a formal approach to justify such a profile.
\item The profile is different from the profile in \cite{MZDuke97} and \cite{MZ07}, therefore new estimates are needed.
\item In order to handle the new scaling, we introduce a new shrinking set to trap the solution,  see Definition \ref{defthess}. Finding such set is not trivial, in particular in the critical case, where we need much more details in the expansions of the rest term (see Appendix \ref{appendix-expansion-R}).
\end{itemize}
Then, following \cite{MZDuke97}, the proof is divided in two steps. First, we reduce the problem to a finite dimensional case. Second, we solve the finite time dimensional problem and conclude by contradiction using index theory.  More precisely, the proof is performed in the framework of the similarity variables defined below in \eqref{chauto}. 
We linearize the self-similar solution around the profile $\varphi_0$ and we obtain $q$ (see \eqref{definitionq} below). Our goal is to guarantee that $q(s)$ belongs to some set $\Vg_A(s)$ (introduced in Definition \ref{defthess}), which shrinks to $0$ as $s\to +\infty$. The proof relies on two arguments:

\begin{itemize}
\item The linearized equation gives two positives mode; $\tilde Q_0$ and $\tilde q_1$, one zero modes ($\tilde q_2$) and an infinite dimensional negative part. The negative part is easily controlled by the effect of the heat kernel. The control of the zero mode is quite delicate (see Part 2: Proof of Proposition \ref{propode}, page \pageref{delicatq}). Consequently, the control of $q$ is reduced to the control of its positive modes.

\item The control of the positive modes $\tilde Q_0$ and $\tilde q_1$ is handled thanks to a topological argument based on index theory (see the argument at page \pageref{indextheory}). 
\end{itemize}

\medskip

The organization of the rest of this paper is as follows. In Section \ref{sectFormapp}, we explain formally how we obtain the profile. In Section \ref{formpb}, we give a formulation of the problem in order to justify the formal argument. Section \ref{existence} is divided in two subsections. In Subsection \ref{pwtr} we give the proof of the existence of the profile assuming technical details. In particular, we construct a shrinking set and give an example of initial data giving rise to the blow-up profile. Subsection \ref{tecnicsect} is devoted to the proof of technical results which are needed in the proof of existence. Finally, in Section \ref{pth1}, we give the proof of Theorem \ref{thm1} and Theorem \ref{Newthm}.

\medskip

\section{Formal approach }\label{sectFormapp}
The aim of this section is to explain formally how we derive the behavior given in Theorem \ref{thm1}.  In particular, how to obtain the profile $\varphi_0$ in \eqref{deffi0}, the parameter $b_{cri}$ in \eqref{defi-b-cri}.  In fact, let us  consider CGL, given by \eqref{GL} in the case where  $\alpha=0$ (the case $\alpha \neq 0$ is the same, thanks to what we mentioned before in Remark \ref{nu0}).

\medskip
Firstly,  we consider  an arbitrary $T>0$ and introduce the following self-similar variable  transformation of equation \eqref{GL}, defined by the following
\beqtn
w(y,s)=(T-t)^{\frac{1+i\delta}{p-1}}u(x,t),\;\;y=\frac{x}{\sqrt{T-t}},\;\; s=-\log(T-t).
\label{chauto}
\eeqtn

As a matter of fact, if  $u(x,t)$ satisfies \eqref{GL} for all $(x,t)\in \R^N\times [0,T)$, then $w(y,s)$ satisfies for all $(x,t)\in \R^N\times [-\log T,+\infty)$ the following equation
\beqtn
\frac{\pa w}{\pa s}=(1+i\beta)\Delta w-\frac 1 2 y.\nabla w-\frac{1+i\delta}{p-1}w+(1+i\delta)|w|^{p-1}w,
\label{equa-w}
\eeqtn
for all $(y,s)\in \R^N\times [-\log T,+\infty)$. Thus constructing a solution $u(x,t)$ for the equation \eqref{GL} that blows up at $T<\infty$ like $(T-t)^{-\frac{1}{p-1}}$  is reduced to constructing a global solution $w(y,s)$ for equation \eqref{equa-w} such that
\beqtn
0<\varepsilon\leq \lim_{s\to\infty} \|w(s)\|_{L^\infty(\R^N)}\leq \frac 1 \varepsilon.
\eeqtn

Inspired by  the work of Popp and \textit{al} \cite{PSKKPD98} and Nouaili and Zaag in  \cite{NZARMA18}, we consider a new form  of   $w$.  More precisely, we suppose that     $ w (y,s) =w(|y|, s) $ and $w \left(r s^{\frac{1}{4}}, s \right) = D(r,s),$ where $ r = |y|$. Then, it is easy to see that $D$ satisfied the following
\begin{eqnarray}
\partial_s D(r,s) & = & (1 + i \beta) D''(r,s) \frac{1}{\sqrt s}  + \left( \frac{1}{4 s} - \frac{1}{2}\right) r D'(r,s)  - \frac{1 + i \delta}{p-1}D(r,s)\\
& +& (1 + i \delta)|D|^{p-1}D.\nonumber
\label{eqD}
\end{eqnarray}
In addition to that, we also assume that we can write $D$ under  the following form

\begin{equation}\label{form-D=R.ei-pi}
D(r,s)= R(r,s) e^{i \varphi (r,s)},
\end{equation} 
where  $R$ and $\varphi$ are real   functions.   In particular, using equation  \eqref{eqD}, we obtain the following system, satisfied by $R$ and $\varphi$

\begin{equation} \label{equa-R-and-varphi}
\hspace{-0.2cm}\left\{
\begin{array}{lll}
\partial_s R  &=& \frac{1}{\sqrt s} \left[ R'' - R(\varphi')^2 - \beta (2 R' \varphi'  + R \varphi'')\right] + R'r \left(\frac{1}{4s} - 
\frac{1}{2} \right) - \frac{R}{p-1} + |R|^{p-1}R,\\
\partial_s \varphi &=& \frac{1}{\sqrt s} \left[ \varphi'' - \beta (\varphi')^2 + \frac{1}{R} \left( 2 R' \varphi' + \beta R'' \right)   \right]+  \varphi' r \left( \frac{1}{4s} - \frac{1}{2} \right)  - \frac{\delta}{p-1} + \delta |R|^{p-1}.
\end{array}
\right .
\end{equation}
In the following, we will consider the following ansatz, inspired by the work of Popp and \textit{al} \cite{PSKKPD98}
 \begin{eqnarray*}
R(r,s) &=& R_0(r) + \frac{R_1(r)}{\sqrt s} + \frac{R_2(r)}{s}  + ...\\
\varphi (r,s)& = & \Phi (s) + \varphi_0(r) + \frac{\varphi_1(r)}{\sqrt s}  + \frac{\varphi_2(r)}{s} + ...
\end{eqnarray*}
where $\Phi(s) = \nu \sqrt s  + \mu \ln s$ and  $\nu, \mu$ unknown. 

\medskip
Beside that, we introduce the notation that $f' = \partial_r f$. By looking at the leading order, we formally derive that $R_0$ and $\varphi_0$ should satisfy
\begin{equation}\label{equa-R-0}
 -\frac{1}{2} R_0' r  - \frac{R_0}{p-1}  + |R_0|^{p-1}R_0 = 0,
\end{equation}
and
\begin{equation}\label{equa-varphi-0}
- \frac{1}{2} \varphi_0' r - \frac{\delta}{p-1} + \delta |R_0|^{p-1} = 0.
\end{equation}
Hence, we can explicitly  solve \eqref{equa-R-0}   and obtain the following
\begin{equation}\label{solution-R-0}
R_0 (r) = \left( p-1 + b r^2\right)^{- \frac{1}{p-1}},
\end{equation}
for some constant $b \in \mathbb{R}$. However,  our aim is finding a global profile, so $b$ will be chosen as a positive constant.  From  \eqref{equa-varphi-0} and \eqref{solution-R-0}, we deduce that
\begin{equation}\label{solution-varphi-0}
\varphi_0 (r) = - \frac{\delta}{p-1} \ln \left( p-1 + b r^2 \right).
\end{equation}

We now look  at the order $\frac{1}{\sqrt s}$ in the system  \eqref{equa-R-and-varphi} and  we deduce the following
\begin{equation}\label{equa-R-1}
-\frac{1}{2} R_1' r  - \frac{R_1}{p-1} + p |R_0|^{p-1} R_1 +  R_0'' - R_0 \varphi_0'^2 - \beta \left( 2 R_0' \varphi_0' + R_0 \varphi_0''\right) = 0, 
\end{equation}
and
\begin{equation}\label{eua-varphi-1}
- \frac{1}{2} \varphi_1' r +\varphi_0'' - \beta \varphi_0'^2 + R_0^{-1} (2 R_0' \varphi_0' + \beta R_0'') + \delta (p-1) |R_0|^{p-2} R_1 - \frac{\nu }{2} = 0.
\end{equation}

\medskip

We now solve  \eqref{equa-R-1}, and obtain 

\begin{eqnarray}
& & R_1 (r) = \frac{r^2}{(p-1 + br^2)^{\frac{p}{p-1}}}\label{formu-R-1} \\
& \times & \left[ - \frac{2 b
 (\delta \beta - 1)}{p-1}r^{-2}  + \frac{8 b^2 (p - (p+1)\delta
  \beta - \delta^2)}{(p-1)^3}\left( \ln |r| - \frac{ \ln (p-1+ b r^2)}{2} 
  \right) +\mathcal{C} \right].
  \label{defC}
\end{eqnarray}

$R_1$ and its derivatives should be bounded at least in some inner region around the pulse center. Thus, the contribution $\log|r|$ are suppressed.  Consequently, we obtain

\begin{equation}\label{condition-delta-beta}
p - (p +1) \delta \beta - \delta^2 = 0,
\end{equation}
which is the   critical condition in our paper.   Besides that,   the constant $\mathcal{C}$ is an unknown    constant,  depending on $\beta, \delta$ and  we particularly  know from   \cite{NZARMA18}      that $\mathcal{C}(0, \delta)=2\frac{pb^2}{(p-1)^3}$. 

\medskip
We now solve \eqref{eua-varphi-1} and obtain
\begin{eqnarray*}
\varphi_1 (r) = \left[ - \nu - \frac{4 b \beta (1 + \delta^2)}{(p-1)^2} \right] \ln |r| + \frac{2 \beta (1 + \delta^2) b}{(p-1)^2} \ln (p-1 + b r^2) \\
 - \frac{2 b}{ (p-1)^2}\left(  (p +  3) \delta + \beta(2 p + \delta^2(p-3)) + \frac{ \mathcal{C}\delta(p-1)^3}{2b^2} \right) (p-1 + b r^2)^{-1}.
\end{eqnarray*}
 By the regularity of $ \varphi_1$ at $0$, the contribution of $\log |r|$  need to be removed. This gives us the following condition
\begin{equation}\label{conditon-on-nu}
\nu = -\frac{ 4 b \beta (1 + \delta^2)}{(p-1)^2}.
\end{equation}

\medskip

At the order  $\frac{1}{s}$, we get 
\begin{equation}\label{equation-R-2}
- \frac{1}{2} R_2' r - \frac{R_2}{p-1} + p |R_0|^{p-1}R_2 + F_3 = 0,
\end{equation}

and 

\begin{equation}\label{equation-varphi-2}
- \frac{1}{2} \varphi_2' r + F_4(r) = 0,
\end{equation}

where
\begin{eqnarray}
F_3 (r) & = & R_1'' - 2 R_0 \varphi_0' \varphi_1' - R_1 \varphi_0'^2 - 2 \beta R_0' \varphi_1' - 2\beta R_1' \varphi_0' - \beta R_0 \varphi_1'' - \beta R_1 \varphi_0'' \nonumber\\
&+& \frac{1}{4}  R_0' r + \frac{p(p-1)}{2}|R_0|^{p-3} R_0 R_1^2.\label{def-of-F_3}\\
F_4(r) &= & \varphi_1'' - 2 \beta \varphi_0' \varphi_1' + 2  \left( R_0^{-1} R_0' \varphi_1' + R_0^{-1} R_1' \varphi_0' - R_0' \varphi_0' R_0^{-2} R_1\right)  - \mu \label{def-of-F_4}\\
& + & \beta \left( R_0^{-1} R_1'' - R_0'' R_0^{-2} R_1\right) - \frac{1}{4}\varphi_0' r + \delta \left( (p-1) |R_0|^{p-3}R_0 R_2 + \frac{(p-1)(p-2)}{2} |R_0|^{p-3} R_1^2 \right)\nonumber.
\end{eqnarray}

\medskip
We solve \eqref{equation-R-2}, by using variation of constant and we obtain that
\begin{eqnarray*}
R_2 = H^{-1}(r) \left(  \displaystyle \int F_3 \frac{2 H}{r}  H dr \right),
\end{eqnarray*}
where
$$H(r)  = \frac{(p-1 + b r^2)^{ \frac{p}{p-1}}}{r^2}.$$

In particular, by a careful calculation we can write

 \begin{eqnarray*}
\frac{2 H }{r} F_3= P(p, \delta, \beta) \frac{1}{r} + \text{  ``regular term''}.
\end{eqnarray*}

 After integrating,  we can see  that  $\frac{1}{r}$ will become  $\ln r$. So, in order to remove this term, we need to have
 $$P = 0.$$
The computation of $P$ is straightforward though a bit lengthy, the interested reader will find details in Appendix \ref{expand-derivation-P-}. Then $P$ is given by the following formula;  
\begin{eqnarray*}
 P &=& 2 \left\{ - \frac{b}{2 (p-1)} - \frac{b^3}{(p-1)^4} (\delta \beta -1) \left( \frac{ - 8 \delta^4 + (12p - 4p^2) \delta^2 + 8 p^3 }{(p+1)\delta^2}  \right) \right.\\
&+& \frac{b^3}{(p-1)^5} (\delta \beta -1) \left( \frac{(8p^2 + 20 p - 8) \delta^4 + 20p^3 \delta^2 - 8 p^2 (p^2 -1)}{(p+1) \delta^2} \right)\\
& + & \frac{b^3}{(p-1)^4} (\delta \beta -1) \left( \frac{4p (1 + \delta^2)}{p+1} \right),\\
&+&\left. \frac{b^3}{(p-1)^5} (\delta \beta -1) \left( \frac{ (- 4 p^2 -16p -32) \delta^4  + (-16p^3 - 68p^2) \delta^2  - 32p^4}{(p+1) \delta^2} \right) \right\}.
\end{eqnarray*}

\noindent
Then, the condition $P=0$ is equivalent to the following

\begin{eqnarray*}
& & \frac{b}{2(p-1)} \\
&= & \frac{b^3}{(p-1)^5} \frac{(\delta \beta -1)}{(p+1)\delta^2} \left( (8p^2 + 8 p - 48) \delta^4 + (8p^3 - 80 p^2  + 8p)\delta^2  - 48p^4 + 8p^3 + 8p^2 \right)\\
&=& \frac{b^3}{(p-1)^5} \frac{8(\delta \beta -1)}{(p+1)\delta^2}  \left( (p^2 +p - 6) \delta^4  + p (p^2 - 10 p + 1) \delta^2 +p^2 (- 6p^2 + p + 1)  \right).
\end{eqnarray*}
Plugging the critical condition $ p - (p+1) \delta \beta - \delta^2 = 0$ into the above equality, we obtain

\begin{eqnarray}
b^2  &=&  \frac{(p-1)^4 (p+1)^2 \delta^2}{ - 16 (1 + \delta^2) L}, \label{formula-b-formal}
\end{eqnarray}
where
\begin{eqnarray*}
L &=& (p^2 +p - 6) \delta^4  + p (p^2 - 10p + 1) \delta^2 +p^2 (- 6p^2 + p + 1).
\end{eqnarray*}
As a matter of fact, in order to justify the positivity of $b^2$, we need to have
$$L < 0.$$

\noindent
This condition is satisfied if and only if 
\begin{equation}\label{critical-number}
\delta \in (-p_{cri}, p_{cri}), \delta \neq  0.
\end{equation}
where 
$$p_{cri}= \sqrt{\frac{p (2p - 1)}{p-2}} \text{ if } p > 2, \text{and } p_{cri} = + \infty \text{ if } p \in (1,2] $$
which was already introduced in Definition \ref{defini-cris-condition}.

\medskip
In the following, we try to determine $\mu$. We use  \eqref{equation-varphi-2}, taking $r=0$, we obtain 
 $$F_4(0) = 0.$$
By using the definition of $F_4$ we directly derive (see Appendix \ref{expand-derivation-P-} for more details),  
  
  \begin{eqnarray*}
  \mu &=& \frac{b^2}{(p-1)^4} \left\{   8 (p+1) \delta + 8p \beta  + (4p +8) \delta^2 \beta + (16p-8) \delta \beta^2 + (8p -16) \delta^3 \beta^2   \right\} \\
  &+& \frac{2 \beta (1 + \delta^2) \mathcal{C}}{p-1}.
  \end{eqnarray*}

\noindent 
  Note that in $\mu$, there is the unknown constant $\mathcal{C}$ already introduced in \eqref{defC}. 
  
\bigskip
\noindent
\textbf{Summary:} From the above approach, we can formally derive  the profile of our solution as in Theorem \ref{Theorem-existence}

\begin{eqnarray*}
w(y,s) \sim  e^{i (\nu \sqrt s + \mu \ln s)} \left( p-1 + b_{cri} \frac{|y|^2}{s}\right)^{-\frac{1+ i\delta}{p-1}}.
\end{eqnarray*}

\medskip
\noindent
In other word, we have

$$ (T-t)^{ -\frac{1+i \delta}{p-1}} e^{   -i\left( \nu\sqrt{ |\ln (T-t)|}  + \mu \ln (|\ln (T-t)|) \right )  }u(x,t) \sim \left( p-1 + b_{cri} \frac{|x|^2}{(T-t)|\ln (T-t)|^\frac{1}{2}} \right)^{-\frac{1+ i\delta}{p-1}},$$
where
\begin{eqnarray*}
\mu &=& \frac{b^2}{(p-1)^4} \left\{   8 (p+1) \delta + 8p \beta  + (4p +8) \delta^2 \beta + (16p-8) \delta \beta^2 + (8p -16) \delta^3 \beta^2   \right\} \\
& +& \frac{2 \beta (1 + \delta^2) \mathcal{C}}{p-1},\\
\nu  &=& -\frac{ 4 b \beta (1 + \delta^2)}{(p-1)^2},\\
b_{cri} &=&   \frac{(p-1)^2 (p+1) \delta}{\sqrt{ 16 (1 + \delta^2) (p(2p-1) - (p-2)\delta^2)( (p+3)\delta^2 + p(3p+1))}}.
\end{eqnarray*}
\begin{rem}We note, that coming at this level in our formal approach, we are not able to determine explicitly $\mu$. But we will see in the rigorous proof that we can obtain the existence and uniqueness of $\mu$ (see equation \eqref{finding-mu-cri}, page \pageref{mu}).

\end{rem}
\section{Formulation of the problem}\label{formpb}
We recall that we consider CGL, given by \eqref{GL}, when $\alpha=0$, as we mentioned before in Remark \ref{nu0}.\\
The preceding calculation is purely formal. However, the formal expansion provides us with the profile of the function ($w(y,s)=e^{i ( \nu \sqrt s +  \mu \ln s)}\left (\varphi_0(\frac{y}{s^{1/4}})+...\right ))$. Our idea is to linearize equation \eqref{equa-w} around that profile and prove that the linearized equation as well as the nonlinear equation have a solution that goes to zero as $s\to \infty$. Let us introduce $q(y,s)$ and $\theta (s)$ such that

\beqtn
\displaystyle w(y,s)= e^{i\left ( \nu \sqrt{s} +  \mu \log s+\theta (s)\right )}\left (\varphi(y,s)+q(y,s) \right ),
\label{formulav}
\eeqtn

where 
\begin{equation}\label{defi-varphi}
\varphi\left (y,s\right )=\varphi_0\left (\frac{y}{s^{1/4}}\right )+(1+i\delta)\frac{a}{s^{1/2}}\equiv\kappa ^{-i\delta}\left ( p-1+b\frac{|y|^2}{s^{1/2}}\right )^{-\frac{1+i\delta}{p-1}}+(1+i\delta)\frac{a}{s^{1/2}},
\end{equation}

\beqtn
  \nu =   -\frac{ 4 b \beta (1 + \delta^2)}{(p-1)^2} \text{ and  }   a=2\kappa (1-\beta\delta)\frac{b}{(p-1)^2},
\label{definitionq}
\eeqtn

\noindent
and the other constants $b,\mu$ will be defined  in the rigorous proof. 

\medskip
In order to guarantee the uniqueness of the couple $(q,\theta)$ an additional constraint is needed, see \eqref{eqmod} below; we will choose $\theta(s)$ such that we kill one the neutral modes of the linearized operator.\\

Note that $\varphi_0(z)$ has been exhibited in the formal approach and satisfies the following equation 
\beqtn
-\frac 12 z\cdot \nabla \varphi_0 -\frac{1+i\delta}{p-1}\varphi_0+(1+i\delta)|\varphi_0|^{p-1}\varphi_0=0,
\label{eqfi0}
\eeqtn
 which makes $\varphi (y,s)$ an approximate solution of  \eqref{equa-w}. In addition to that, if $w$ satisfies this equation, then $q$ satisfies the following equation
\beqtn
\label{eqq1}
\frac{\pa q}{\pa s}=\Lg_{\beta} q-\frac{(1+i\delta)}{p-1} q +L(q,\theta ', y, s)+R^*(\theta ',y,s)
\eeqtn
where
\beqtn
\begin{array}{lll}
\Lg_{\beta} q&=&(1+i\beta)\Delta q -\frac 12 y\cdot\nabla q,\\[0.2cm]
L(q,\theta',y,s)&=&(1+i\delta)\left\{|\varphi+q|^{p-1}(\varphi+q)-|\varphi|^{p-1}\varphi -i\left(  \frac{\nu}{2 \sqrt s}  + \frac \mu s+\theta'(s)\right)q\right\},\\[0.2cm]
R^*(\theta',y,s)&=&R(y,s)-i\left(   \frac{\nu}{2 \sqrt s}  +  \frac \mu s+\theta '(s)\right )\varphi,\\[0.2cm]
R(y,s)&=&-\frac{\pa \varphi}{\pa s}+(1+ i \beta)\Delta \varphi-\frac 12 y\cdot \nabla \varphi-\frac{(1+i\delta)}{p-1}\varphi+(1+i\delta)|\varphi|^{p-1}\varphi.
\end{array}
\label{eqqd1}
\eeqtn

\medskip
Our aim is to find a $\theta\in C^1([-\log T,\infty),\R)$ such that equation \eqref{eqq} has a solution $q(y,s)$ defined for all $(y,s)\in \R^N\times [-\log T,\infty)$ such that
\[q(y,s)=\frac{\Fg(y)}{s}+v(y,s),\]
where $\Fg$ is defined by \eqref{defFg} in Theorem \ref {Theorem-existence2} and
\[\|v(s)\|_{L^\infty}\to 0 \mbox{ as }s\to \infty.\]


\medskip
From \eqref{eqfi0}, one sees that the variable $z=\frac{y}{s^{1/4}}$ plays a fundamental role. Thus, we will consider the dynamics for $|z|>K$, and $|z|<2K$ separately for some $K>0$ to be fixed large.
\subsection{The outer region where $|y|>Ks^{1/4}$}
Let us consider a non-increasing cut-off function $\chi_0\in C^{\infty}(\R^+,[0,1])$ such that $\chi_0(\xi)=1$ for $\xi<1$ and $\chi_0(\xi)=0$ for $\xi>2$ and introduce
\beqtn
\chi(y,s)=\dsp\chi_0\left(\frac{|y|}{Ks^{1/4}}\right),
\label{defchi14}
\eeqtn
where $K$ will be fixed large. Let us define
\beqtn
q_e(y,s)=e^{\frac{i\delta}{p-1}s}q(y,s)\left( 1-\chi(y,s)\right),
\label{defiqe}
\eeqtn
and note that $q_e$ is the part of $q(y,s)$, corresponding to the non-blowup region $|y|>Ks^{1/4}$. As we will explain in subsection \eqref{outreg}, the linear operator of the equation satisfied by $q_e$ is negative, which makes it easy to control $\|q_e(s)\|_{L^\infty}$. This is not the case for the part of $q(y,s)$ for $|y|<2Ks^{1/4}$, where the linear operator has two positive eigenvalues, a zero eigenvalue in addition to infinitely many negative ones. Therefore, we have to expand $q$ with respect to these eigenvalues in order to control $\|q(s)\|_{L^\infty(|y|<2 K s^{1/4})}$. This requires more work than for $q_e$. The following subsection is dedicated to that purpose. From now on, $K$ will be fixed constant which is chosen such that $\|\varphi(s')\|_{L^\infty(|y|> K {s'}^{1/4})}$ is small enough, namely $\|\varphi_0(z)\|_{L^\infty(|z|>K)}^{p-1}\leq \frac{1}{C(p-1)}$ (see subsection  \eqref{outreg} below, for more details). 

\subsection{The inner region where $|y|< 2K s^{1/4}$}
If we linearize the term $L(q,\theta ', y, s)$ in equation \eqref{eqq1}, then we can write \eqref{eqq1} as

\beqtn
\label{eqq}
\frac{\pa q}{\pa s}=\Lg_{\beta, \delta} q-i\left ( \frac{\nu}{ 2 \sqrt s} +  \frac \mu s+\theta'(s)\right ) q +V_1 q+ V_2 \bar q+B(q,y,s)+R^*(\theta',y,s),
\eeqtn
where
\beqtn
\begin{array}{lll}
\Lg_{\delta, \beta}  q&=&(1+ i\beta)\Delta q -\frac 12 y\cdot\nabla q+(1+i\delta )\Re q,\\
V_1(y,s)&=&(1+i\delta)\frac{p+1}{2}\left(|\varphi|^{p-1}-\frac{1}{p-1}\right),\;\;V_2(y,s)=(1+i\delta)\frac{p-1}{2}\left(|\varphi|^{p-3}\varphi^2-\frac{1}{p-1}\right),\\
B(q,y,s)&=&(1+i\delta)\left ( |\varphi+q|^{p-1}(\varphi+q)-|\varphi|^{p-1}\varphi-|\varphi|^{p-1} q-\frac{p-1}{2}|\varphi|^{p-3}\varphi(\varphi \bar q+\bar\varphi q)\right),\\ 
R^*(\theta',y,s)&=&R(y,s)-i\left ( \frac{\nu}{2 \sqrt s} +   \frac \mu s+\theta '(s)\right )\varphi,\\
R(y,s)&=&-\frac{\pa \varphi}{\pa s}+\Delta \varphi-\frac 12 y\cdot \nabla \varphi-\frac{(1+i\delta)}{p-1}\varphi+(1+i\delta)|\varphi|^{p-1}\varphi
\end{array}
\label{eqqd}
\eeqtn
Note that the term $B(q,y,s)$ is built to be quadratic in the inner region $\dsp |y|\leq K s^{1/4}$. Indeed, we have for all $K\geq 1$ and $s\geq 1$, 
\beqtn
\sup_{|y|\leq 2 K s^{1/4}}|B(q,y,s)|\leq C(K)|q|^2.
\label{estiquadinn}
\eeqtn
Note also that $R(y,s)$ measures the defect of $\varphi(y,s)$ from being an exact solution of \eqref{equa-w}. However, since $\varphi(y,s)$ is an approximate solution of \eqref{equa-w}, one easily derives the fact that
\beqtn
\|R(s)\|_{L^\infty}\leq \frac{C}{\sqrt{s}}.
\label{estR*}
\eeqtn
Therefore, if $\theta'(s)$ goes to zero as $s\to \infty$, we expect the term $R^*(\theta',y,s)$ to be small, since \eqref{eqq} and \eqref{estR*} yield
\beqtn
|R^*(\theta',y,s)|\leq \frac{C}{\sqrt{s}}+|\theta'(s)|.
\label{Rexpect}
\eeqtn
Therefore, since we would like to make $q$ go to zero as $s\to \infty$, the dynamics of equation \eqref{eqq} are influenced by the asymptotic limit of its linear term,
\[\tilde \Lg+V_1 q+V_2\bar q,\]
as $s\to\infty$. In the sense of distribution (see the definitions of $V_1$ and $V_2$ in  \eqref{eqq}  and  $\varphi$ \eqref{defi-varphi}) this limit is $\tilde \Lg$.


\subsection{Spectral properties of $\Lg_\beta$}
Here, we will restrict to $N=1$. We consider the Hilbert space $L^{2}_{|\rho_\beta|}(\R^N,\C)$ which is the set of all $f\in L^{2}_{loc}(\R^N,\C)$ such that 
\[\int_{\R^N} |f(y)|^2|\rho_\beta(y)|dy< +\infty,\]
where

\beqtn
\displaystyle \rho_\beta(y)=\frac{e^{-\frac{|y|^2}{4(1+i\beta)}}}{(4\pi(1+i\beta))^{N/2}}\mbox{   and  }|\rho_\beta(y)|=\frac{e^{-\frac{|y|^2}{4(1+\beta^2)}}}{(4\pi\sqrt{1+\beta^2})^{N/2}}.
\eeqtn
We can diagonalize $\Lg_\beta$ in $L^{2}_{|\rho_\beta|}(\R^N,\C)$. Indeed, we can write 
\[\Lg_\beta q=\frac{1}{\rho_\beta}div (\rho_\beta \nabla q).\]
We notice that $\Lg_\beta$ is formally ``self-adjoint'' with respect to the weight $\rho_\beta$. Indeed, for any $v$ and $w$ in $L^{2}_{|\rho_\beta|}(\R^N,\C)$ satisfying $\Lg_\beta v$ and $\Lg_\beta w$ in $L^{2}_{|\rho_\beta|}(\R^N,\C)$, it holds that
\beqtn
\int v\Lg_\beta w \rho_\beta dy= \int w\Lg_\beta v \rho_\beta dy.
\eeqtn
If we introduce for each $\alpha=(\alpha_1,...,\alpha_N)\in \N^N$ the polynomial
\begin{equation}\label{defi-f-n}
\displaystyle f_{\alpha}(y)=c_\alpha \Pi_{i=1}^{N} H_{\alpha_i} \left(\frac{y_i}{2\sqrt{1+i\beta}}\right),
\end{equation}

where $H_n$ is the standard one dimensional Hermite polynomial and $c_\alpha\in \C$ is chosen so that the term of highest degree in $f_{\alpha}$ is $\Pi_{i=1}^{N} y_{i}^{\alpha_i}$, then, we get a family of eigenfunction of $\Lg_\beta$, ``orthogonal`` with respect to the weight $\rho_\beta$, in the sense that for any different $\alpha$ and $\sigma\in \N^N$
\beqtn
\begin{array}{ll}
\dsp \Lg_\beta f_{\alpha}&=-\frac{\alpha}{2} f_{\alpha},\\
\dsp \int_\R f_{\alpha}(y)f_{\sigma}(y)\rho_\beta(y)dy&=0.
\end{array}
\eeqtn

Moreover, the family $f_\alpha$ is a basis for $L^{2}_{|\rho_{\beta}|}(\R^N,\C)$ considered as a $\C$ vector space. Al the facts about the operator $\Lg_\beta$ and the family $f_\alpha$ can be found in Appendix A of  \cite{MZ07}.  

\subsection{Spectral properties of $\Lg_{\beta,\delta}$}
In the sequel, we will assume $N=1$. Now, with the explicit basis diagonalizing $\Lg_\beta$, we are able to write $\Lg_{\beta,\delta}$ in a Jordan's block's. More precisely, we recall Lemma 3.1 from \cite{MZ07}
\begin{lemma}[Jordan block's decomposition of $\Lg_{\beta,\delta}$]\label{lemma-Jordan-block's} For all $n\in\N$, there exists two polynomials
\beqtn
\begin{array}{ll}
 h_n&=if_n+\sum_{j=0}^{n-1}  d_{j,n}f_j\mbox{, where }  d_{j,n}\in \C\\
\tilde{h}_n&=(1+i\delta)f_n +\sum_{j=0}^{n-1} \tilde d_{j,n}f_j\mbox{, where } \tilde d_{j,n}\in \C,
\end{array}
\eeqtn
of degree $n$ such that 
\beqtn
\begin{array}{ll}
\dsp \Lg_{\beta,\delta} h_n=-\frac n 2  h_n,\\
\dsp \Lg_{\beta,\delta} \tilde h_n=\left(1-\frac n 2\right) \tilde h_n+c_n  h_{n-2},
\end{array}
\eeqtn
with $c_n=n(n-1)\beta(1+\delta^2)$ (and we take $h_k\equiv 0$ for $k<0$). The term of highest of $ h_n$ (resp. $\tilde h_n$) is $iy^n$  (resp. $(1+i\delta)y^n$).
\end{lemma}
\textit{Proof :} See the proof of Lemma 3.1 in \cite{MZ07}. To prove that $c_n=n(n-1)$, we look at the imaginary part of order $n-1$ in the equation  $\Lg_{\beta,\delta} \tilde h_n=\left(1-\frac n 2\right) \tilde h_n+c_n  h_{n-2}$. A simple identification gives the result.

\medskip

\textbf{Remark:} The semigroup and the fundamental solution generated by $(1+i\beta)\Delta v$ have the same regularizing effect independently from $\beta$.

\begin{lemma}[The basis vectors of degree less or equal to $6$]\label{lemma-formula-h-j-tilde-h-j} we have
\[
\begin{array}{ll}
 h_0(y)=i,       &\tilde h_0=(1+i\delta),\\
 h_1(y)=i y,    & \tilde h_1=(1+i\delta)y,\\
 h_2(y)=i y^2+\beta-i(2+\delta\beta),&\tilde h_2=(1+i\delta)(y^2-2+2\beta\delta),
\end{array}
\]
\[
\begin{array}{ll}
 h_4(y)=i y^4+y^2(c_{4,2}+i d_{4,2})+ c_{4,0}+i  d_{4,0}&,\\
c_{4,2}= 6\beta,&d_{4,2}=-6(2+\beta\delta)=-18-6(\beta\delta-1),\\
c_{4,0}=-4\beta(3+\beta\delta),&d_{4,0}=12-6\beta^2+12\beta\delta+2\beta^2\delta^2,\\
\end{array}
\]

\[
\begin{array}{lll}
& \tilde h_4(y)= (1+i\delta)y^4+y^2(12(\beta\delta -1)+i\tilde d_{4,2})+\tilde c_{4,0}+i\tilde d_{4,0}.\\
& \tilde c_{4,2}= 12 (\beta\delta-1), \quad \tilde d_{4,2}=0,\\
& \tilde c_{4,0}=6\beta^2(1+\delta^2)-12(\beta \delta-1), \quad \tilde d_{4,0}=-6\beta^2\delta(3\delta^2+7)-12\delta(\beta\delta+1)

\end{array}
\]
\[
\begin{array}{ll}
& h_6(y)=i y^6+y^4(c_{6,4}+i d_{6,4})+y^2(c_{6,2}+i d_{6,2})+ c_{6,0}+i  d_{6,0},\\
& c_{6,4}=15\beta , \quad d_{6,4}=-15(2+\beta\delta),\\
& c_{6,2}=-60\beta(3+\delta\beta) ,\quad d_{6,2}=-90\beta^2+180+180\beta\delta+30\beta^2\delta^2,\\
& c_{6,0}=180\beta+120\delta\beta^2-45\beta^3+15\beta^3\delta ,\\
& d_{6,0}=-180\beta\delta+55\delta\beta^3-60\delta^2\beta^2-5\beta^3\delta^2+180\beta^2-120,
\end{array}
\]

\[
\begin{array}{l}
\tilde h_6(y)=(1+i\delta) y^6+y^4(\tilde c_{6,4}+i \tilde d_{6,4})+y^2(\tilde c_{6,2}+i \tilde d_{6,2})+ \tilde c_{6,0}+i \tilde  d_{6,0},\\
\tilde c_{6,4}= 30(\beta\delta-1),\;\;\; \tilde d_{6,4}=0,\\
\tilde c_{6,2}=90\beta^2(1+\delta^2)-180(\beta\delta-1),\\
   \tilde d_{6,2}=-90\beta(1+\delta^2)(3\beta\delta+4)+180(\beta\delta-1)(\delta-2\beta),\\

\tilde c_{6,0}=-20\beta^2(1+\delta^2)(11\beta\delta+21)+120(\beta\delta-1)(-2\beta^2+\beta\delta+1) ,\\
\tilde d_{6,0}=270\beta(1+\delta^2)(2+\beta\delta)+\beta^2(1+\delta^2)(140\beta\delta^2-180\beta\delta+390\delta)\\
+60(\beta\delta-1)(2\beta^2\delta-\beta\delta^2+9\beta-4\delta),
\end{array}
\]

Moreover, we have 
\[
\begin{array}{l}
\Lg_{\beta,\delta}\tilde  h_0=\tilde h_0,\\
\Lg_{\beta,\delta}\tilde  h_1=\frac 12 \tilde h_1,\\
\Lg_{\beta,\delta}\tilde  h_2=2\beta(1+\delta^2) h_0=2 i\beta(1+\delta^2),\\
\Lg_{\beta,\delta}\tilde  h_4=-\tilde h_4+ 12 \beta(1+\delta^2) h_2,\\
\Lg_{\beta,\delta}\tilde  h_6=-2 \tilde h_6+30\beta(1+\delta^2) h_4.

\end{array}
\]
\end{lemma}
\textit{Proof :} The proof is straightforward though a bit lengthy.

\begin{corollary}[Basis for the set of polynomials] The family $(h_n\tilde h_n)_n$ is a basis of $\C[X]$, the $\R$ vector space of complex polynomials.

\label{basepolyn}
\end{corollary}

\subsection{Decomposition of $q$}\label{sectidecompq}
For the sake of controlling $q$ in the region $|y|<2K s^{1/4}$, we will expand the unknown function $q$ (and not just $\chi q$ where $\chi$ is defined in \eqref{defchi14})
with respect to the family $f_n$ and the with respect to the $h_n$. We start by writing 


We start by writing
\beqtn
q(y,s)=\sum_{n\leq M} \Qg_n(s) f_n(y)+q_{-}(y,s),
\label{decomp1}
\eeqtn
where $f_n$ is the eigenfunction of $\mathcal L_\beta$ defined in \eqref{defi-f-n}, $Q_n(s)\in \C$, $q_{-}$ satisfies

\[\int q_{-}(y,s)h_n(y)\rho(y)dy=0\mbox{ for all }n\leq M,\]
and $M$ is a fixed even integer satisfying
\beqtn
\dsp M\geq 4\left( \sqrt{1+\delta^2}+1+2 \max_{i=1,2,y\in \R,s\geq 1}|V_i(y,s)|\right),
\label{boundM}
\eeqtn
with $V_{i=1,2}$ defined in \eqref{eqqd}. From \eqref{decomp1}, we have
\beqtn
\Qg_n(s)=\dsp \frac{\dsp \int q(y,s) f_n(y)\rho_\beta(y)dy}{\int f_{n}(y)^{2}\rho_\beta(y)}\equiv F_n(q(s)),
\label{eqQn}
\eeqtn

The function $q_{-}(y,s)$ can be seen as the projection of $q(y,s)$ onto the spectrum of $\Lg_\beta$, which is smaller than $(1-M)/2$. We will call it the infinite dimensional part of $q$ and we will denote it $q_-=P_{-,M}(q)$. We also introduce $P_{+,M}=Id -P_{-,M}$. Notice that $P_{-,M}$ and $P_{+,M}$ are projections. In the sequel, we will denote $P_-=P_{-,M}$ and $P_+=P_{+,M}$.\\
The complementary part $q_+=q-q_-$ will be called the finite dimensional part of $q$. We will expand it as follows
\beqtn
\dsp q_+(y,s)=\sum_{n\leq M}\Qg_n(s)f_n(y)=\sum_{n\leq M} q_n(s) h_n(y)+ \tilde q_n(s) \tilde h_n(y),
\label{decomp2}
\eeqtn
where $\tilde q_n, q_n\in \R$. Finally, we notice that for all $s$, we have
\[\dsp\int q_-(y,s)q_+(y,s)\rho_\beta(y)dy=0.\] 
Our purpose is to project \eqref{eqq} in order to write an equation for  $q_n$ and $\tilde q_n$. For that we need to write down the expression of $q_n$ and $\tilde q_n$ in terms of $\Qg_n$. Since the matrix $(h_n,\tilde h_n)_{n\leq M}$ in the basis of $(i f_n,f_n)$ is upper triangular (see Lemma \ref{lemma-formula-h-j-tilde-h-j}). The same holds for its inverse.  Thus, we derive from \eqref{decomp2}

 \begin{equation}\label{decomp3}
 \begin{array}{rcl}
 \hspace{-0.2cm} & q_n = \text{ Im} \Qg_n (s) -  \delta \text{ Re} \Qg_n (s) + \sum_{j=n+1}^{M} A_{j,n} \text{Im} \Qg_j(s)+ B_{j,n} \text{ Re} \Qg_j(s) \equiv P_{n,m}(q(s)),\\[0.3cm]
\hspace{-0.2cm} & \tilde q_n (s) =  \text{ Re} \Qg_n(s) + \sum_{j=n+1}^M \tilde A_{j,n} \text{ Im}\Qg_{j} (s) + \tilde B_{j,n} \text{ Re} \Qg_j(s) \equiv \tilde P_{n,M} (q(s)),
 \end{array}
 \end{equation}
where all the constants are real. Moreover,  the coefficient of $\text{Im}\Qg_n$ and $\text{ Re} \Qg_n$ in the expression of $q_n$ and $\tilde q_n$ are explicit. This comes from the fact that the same holds for the coefficient of $i f_n$ and $f_n$ in the expansion of $h_n$ and $\tilde h_n$ (see Lemma \ref{lemma-Jordan-block's}).\\
Note that the projector $P_{n,m}(q)$ and $\tilde P_{n,m}(q)$ are well-defined thanks to \eqref{eqQn}. We will project equation  \eqref{eqq} on the different modes $h_n$ and $\tilde h_n$. Note that  from \eqref{decomp1} and \eqref{decomp2}, that
\beqtn
q(y,s)=\left(\sum_{n\leq M}q_n(s) h_n(y) + \tilde q_n(s) \tilde h_n(y)\right)+q_-(y,s),
\label{decompq}
\eeqtn

\noindent
we should keep in mind that  the presentation in \eqref{decompq}  is unique.

\section{Existence}\label{existence}
In this section, we prove the existence of a solution $q(s), \theta(s)$ of problem \eqref{eqq1}-\eqref{eqmod}  such that
\beqtn
\begin{array}{l}
\dsp q(y,s)=\frac{1}{s}\left ( \tilde \Ag_0 \tilde h_0(y)+ \Ag_2 h_2(y)+\tilde \Ag_2\tilde h_2(y)\right )+v(y,s),\\
\mbox{with, for all $M>0$ }\dsp \sup_{|y|<Ms^{\frac 14}}|v(y,s)|\leq C\frac{1+|y|^5}{s^{\frac 32}},\\
\mbox{ and }\dsp |\theta '(s)|\leq \frac{CA^{10}}{s^{\frac 54}}\mbox{ for all }s\in [-\log T,+\infty),
\end{array}
\label{limitev}
\eeqtn
where $ \tilde \Ag_0,\;\Ag_2$ and $\tilde\Ag_2$ are given by \eqref{defAg} in Definition \ref{defthess} and $ h_0(y)$, $h_2(y)$ and $\tilde h_2(y)$ are given in Lemma \ref{lemma-formula-h-j-tilde-h-j}.

\medskip

Hereafter, we denote by $C$ a generic positive constant, depending only on $p, \delta, \beta$ and $K$ introduced in \eqref{defchi14}, itself depending on $p$. In particular, $C$ neither depend on $A$ nor on $s_0$, the constants that will appear shortly and throughout the paper and need to be adjusted for the proof.\\
We proceed in two subsections. In the first, we give the proof assuming the technical details. In the second subsection we give the proof of the technical details.

\subsection{Proof of the existence assuming technical results}\label{pwtr} 
We work in the set of even functions to construct a blow-up solution. However, since we need to prove the stability of the constructed soluton in the set of all functions with no evenness assumption, we have to handle general functions.
 \begin{definition}[A set shrinking to zero]
For all $K>1$, $A\geq 1$ and $s\geq 1$, we define $\Vg_A(s)$ as the set of all $ q\in L^\infty(\R)$ such that
$$
\begin{array}{ll}
\|q_e\|_{L^\infty(\R)}\leq \frac{A^{M+2}}{s^{\frac 14}}    ,                               &\|\frac{q_-(y)}{1+|y|^{M+1}}\|_{L^\infty(\R)}\leq  \frac{A^{M+1}}{s^{\frac{M+2}{4}}}, \\
| q_j|,\;|\tilde q_j|\leq \frac{A^j}{s^{\frac{j+1}{4}}}\mbox{ for all }5\leq j\leq M,  & | q_0| \leq \frac{1}{s^{\frac 32}},| \tilde q_1| \leq \frac{A}{s^{\frac 32}},\;\;| q_1| \leq \frac{A^4}{s^{\frac 32}}.                                                                      
\end{array}
$$

\noindent
In addition to the the other modes will satisfy the following condition: 
\begin{eqnarray*}
  & &  \left| Q_4   \right| \leq  \frac{A^{7}}{s^{7/4}} \text{ and }  \left|\tilde Q_4 \right|  \leq \frac{A^4}{s^{\frac{7}{4}}},\\
  & &  \left| q_3  \right| \leq  \frac{A^{3}}{s^{\frac{3}{2}}} \text{ and }  \left|\tilde  q_3 \right|  \leq \frac{A^3}{s^{\frac{3}{2}}},\\
 & &  \left| Q_2 \right| \leq \frac{A^8}{s^\frac{7}{4}}  \text{ and } \left|\tilde Q_2 \right| \leq \frac{A^{10}}{s^\frac{5}{4}},\\
\end{eqnarray*}
and
$$  \left|  \tilde Q_0  \right| \leq  \frac{A}{s^\frac{7}{4}},$$
where 
\begin{eqnarray}
Q_4 &=&  q_4 - \left(  \frac{1}{2} D_{4,2} \frac{\tilde q_2}{\sqrt{s}} +\left[  \frac{C_{4,2} R^*_{2,1}}{2} + \frac{R^*_{4,2}}{2}\right] \frac{1}{s^{\frac{3}{2}}} \right) ,\label{defi-as-Q-4} \\
&=& q_4-(\frac{\Bg_4 }{s^{\frac{3}{2}}}+ \Cg_4  \frac{\tilde q_2}{\sqrt{s}}),\\
 \tilde Q_4 & =&  \tilde q_4  -\left(  \tilde D_{4,2} \frac{\tilde q_2}{\sqrt{s}}  + \frac{1}{s^{\frac{3}{2}}} \left[  \tilde C_{4,2} R^*_{2,1} + \tilde R_{4,2}^*   \right] \right),\label{defi-as-tilde-Q-4}  \\
&=&\tilde q_4-(\frac{\tilde \Bg_4 }{s^{\frac{3}{2}}}+ \tilde\Cg_4  \frac{\tilde q_2}{\sqrt{s}}),\\
 \tilde Q_0 &=& \tilde q_0 -  \left( \frac{\tilde q_2}{\sqrt{s}} \left[  \frac{\nu \tilde L_{0,2}}{2} -\tilde D_{0,2}  - \frac{\tilde \Theta^*_{0,0} c_2}{\kappa} \right]   - \frac{\tilde R^*_{0,1}}{s}  \right)\\
 &-& \left( \frac{1}{s^{\frac{3}{2}}} \left[ - \tilde X_0 + \frac{\nu \tilde K_{0,2}R^*_{2,1} }{2} - \tilde C_{0,2}. R_{2,1}^* \right] \right),\label{defi-as-tilde-Q-0}\\
&=&\tilde q_0-(\frac{\tilde \Ag_0}{s}+\frac{\tilde \Bg_0 }{s^{\frac{3}{2}}}+ \tilde \Cg_0  \frac{\tilde q_2}{\sqrt{s}}),
\end{eqnarray}
and 
\begin{eqnarray}
  Q_2  &=&  q_2 - \left(\frac{\tilde q_2}{\sqrt s} \left[ D_{2,2} - \frac{\nu}{2} (1 + \delta^2) + c_4 \tilde D_{4,2} + \frac{ \Theta^*_{2,0} c_2}{\kappa} \right] + \frac{R^*_{2,1}}{s} \right)\nonumber \\
  &-& \left(  \frac{1}{s^{\frac{3}{2}}} \left[X_2 +c_4 [ \tilde C_{4,2} R_{2,1}^* + \tilde R^*_{4,2} ] - D_{2,0}.\tilde R_{0,1}^* \right] \right),\label{defi-as-Q-2}\\
  &=& q_2-(\frac{ \Ag_2}{s}+\frac{\Bg_2 }{s^{\frac{3}{2}}}+ \Cg_2  \frac{\tilde q_2}{\sqrt{s}}),\\
  \tilde Q_2 & =& \tilde q_2 - \frac{\tilde{\mathcal{A}}_2}{s}\label{defi-as-tilde-Q-2},
\end{eqnarray}
and
\beqtn
\tilde\Ag_0=-\tilde R^*_{0,1},\;\;
\tilde\Ag_2=-\frac{R_{0,1}^*}{c_2},\;\; 
\Ag_2=R^*_{2,1}
\label{defAg}
\eeqtn
with $  R_{0,1}^* $, $\tilde R^*_{0,1}$, $R^*_{2,1}$ and  $c_2 $ are defined as in  page \pageref{page-constant-R^*-j-k}  and    Lemma \ref{lemma-Jordan-block's}, respectively. The constants
$C_{i,j},\;\; \tilde C_{i,j},\;\; D_{i,j}$ and $ \tilde C_{i,j}$ are defined in \eqref{Projection-of-potentials}. The constants $K_{i,j},\;\; \tilde K_{i,j},\;\; L_{i,j}$ and $ \tilde L_{i,j}$ are defined by \eqref{defKL} and \eqref{deftildeKL}.
\label{defthess}
\end{definition}

Since $A\geq 1$, the  sets $\Vg_A(s)$ are increasing (for fixed $s$) with respect to $A$ in the sense of inclusions.\\
We also show the following property of elements of $\Vg_A(s)$:\\
For all $A\geq 1$, there exists $s_{01}(A)\geq 1$, such that for all $s\geq s_{01}$ and $r\in\Vg(A)$, we have
\beqtn
\|r\|_{L^\infty(\R)}\leq C(K)\frac{A^{M+2}}{s^{\frac 14}},
\label{estilinfty}
\eeqtn
where $C$ is a positive constant (see Claim \ref{propshrinset} below for the proof).\\
By \eqref{estilinfty}, if a solution $q$ stays in $\Vg(A)$ for $s\geq s_{01}$, then it converges to $0$ in $L^\infty(\R)$.

The solution of equation \eqref{eqq} will be denoted by $q_{s_0,d_0,d_1}$ or $q$ when there is no ambiguity. We will show that if $A$ is fixed large enough, then, $s_0$ is fixed large enough depending on $A$, we can fix the parameters $(d_0,d_1)\in [-2,2]^2$, so that the solution $v_{s_0,d_0,d_1}\to 0$ as $s\to \infty$ in $L^{\infty}(\R)$, that is \eqref{limitev} holds. Our construction is built on a careful choice of the initial data of $q$ at a time $s_0$. We will choose it in the following form:

\begin{definition}[Choice of initial data] Let us define, for $A\geq 1$, $s_0=-\log T>1$ and $d_0,d_1\in \R$,  the function
\beqtn\begin{array}{l}
 \psi_{s_0,d_0,d_1}(y)=
 
 \dsp\left [\left( \frac{A}{s_{0}^{7/4}}  \tilde d_0+\frac{\tilde \Ag_0}{s_0}+\frac{\Bg_0}{s^{3/2}_0}+\frac{\Cg_0\tilde \Ag_2}{s^{3/2}_0}\right)\tilde h_0 +\frac{A}{s_{0}^{3/2}}\tilde d_1\tilde h_1(y)+d_0 h_0 \right .\\
\dsp  +\frac{\tilde \Ag_2}{s_0} \tilde h_2+ \left (\frac{ \Ag_2}{s_0}+\frac{\Bg_2}{s^{3/2}_0}+\frac{\Cg_2\tilde \Ag_2}{s^{3/2}_0} \right ) h_2\\
\dsp +\left .\left ( \frac{\tilde \Bg_4}{s^{3/2}_0}+\frac{\tilde \Cg_4\tilde \Ag_2}{s^{3/2}_0}\right )\tilde h_4  +\left ( \frac{\Bg_4}{s^{3/2}_0}+\frac{ \Cg_4\tilde \Ag_2}{s^{3/2}_0}\right )h_4 \right ]\chi (2y,s_0),
 \end{array}
\label{definitdata1}
\eeqtn
where $s_0 = - \log  T$ and  $h_i$, $\tilde h_i$, $i=0,1,2,3,4$ are given in Lemma  \ref{lemma-formula-h-j-tilde-h-j}, $\chi$ is defined by \eqref{defchi14} and $d_0=d_0(\tilde d_0,\tilde d_1)$ will be fixed later in (i) of Proposition \ref{propinitialdata}. The constants $\tilde \Ag_i,\Ag_i$, $\tilde \Bg_i,\Bg_i$, $\tilde \Cg_i,\Cg_i$, for $i=0,2,4$ are given by (\ref{defi-as-Q-4}$-$\ref{defi-as-tilde-Q-2}).
\end{definition}
\begin{rem}
Let us recall that we will modulate the parameter $\theta$ to kill one of the neutral modes, see equation \eqref{eqmod} below. It is natural that this condition must be satisfied for the initial data at $s=s_0$. Thus, it is  necessary that we choose $d_0$ to satisfy condition \eqref{eqmod}, see  \eqref{d2} below. 
\end{rem}


\medskip

\noindent So far, the phase $\theta(s)$ introduced in \eqref{formulav} is arbitrary, in fact as we will show below in Proposition \ref{enopropmod}. We can use a modulation technique to choose $\theta(s)$ in such a way that we impose the condition 

\beqtn
 P_{0,M} (q(s))=0,
\label{eqmod}
\eeqtn
which allows us to kill the neutral direction of the operator $\tilde \Lg$ defined in \eqref{eqq}.
Reasonably, our aim is then reduced to the following proposition:

\begin{prop}[Existence of a solution trapped in $\Vg_A(s)$] There exists $A_2\geq 1$ such that for $A\geq A_2$ there exists $s_{02}(A)$ such that for all $s_0\geq s_{02}(A)$, there exists $(\tilde d_0,\tilde d_1)$ such that if $q$ is the solution of \eqref{eqq}-\eqref{eqmod}, with initial data given by \eqref{definitdata1} and \eqref{d2}, then $v\in \Vg_A(s)$, for all $s\geq s_0$.
\label{propsop1}
\end{prop}
This proposition gives the stronger convergence to $0$ in $L^\infty(\R)$ thanks to \eqref{estilinfty}.\\
Let us first be sure that we can choose the initial data such that it starts in $\Vg_A(s_0)$. In other words, we will define a set where where will be selected the good parameters $(\tilde d_0,\tilde d_1)$ that will give the conclusion of Proposition \ref{propsop1}. More precisely, we have the following result:

\begin{prop}[Properties of initial data] For each $A\geq 1$, there exists $s_{03}(A)>1$ such that for all $s_0\geq s_{03}$:\\
(i) $P_{0,M}\left (i\chi (2y,s_0)\right)\not = 0$ and the  parameter $d_0(s_0,\tilde d_0,\tilde d_1)$ given by 
\beqtn
\begin{array}{lll}\label{d2}
d_0(s_0,\tilde d_0,\tilde d_1)&=&-\dsp\frac{A}{s_{0}^{3/2}}\tilde d_1  \frac{P_{0,M}\left(\tilde h_1\chi (2y,s_0) \right)}{ P_{0,M}\left (i\chi (2y,s_0)\right)} \\
&&-\dsp\left (\frac{A}{s_{0}^{7/4}} \tilde d_0+\frac{\tilde \Ag_{0}}{s_0}+\frac{\tilde \Bg_0 }{s^{3/2}_0}+\frac{\tilde \Cg_0 \tilde\Ag_2}{s^{3/2}_0}\right )\frac{P_{0,M}\left(\tilde h_0\chi (2y,s_0) \right)}{ P_{0,M}\left (i\chi (2y,s_0)\right)} \\
&&-\dsp\left (\frac{ \tilde \Ag_{2}}{s_0}\right )\frac{P_{0,M}\left( \tilde h_2\chi (2y,s_0) \right)}{ P_{0,M}\left (i\chi (2y,s_0)\right)}\\
 &&-\dsp\left (\frac{ \Ag_{2}}{s_0}+\frac{ \Bg_2 }{s^{3/2}_0}+\frac{ \Cg_2 \tilde\Ag_2}{s^{3/2}_0}\right )\frac{P_{0,M}\left( h_2\chi (2y,s_0) \right)}{ P_{0,M}\left (i\chi (2y,s_0)\right)}\\
&&-\dsp\left (\frac{ \tilde \Bg_4 }{s^{3/2}_0}+\frac{\tilde \Cg_4 \tilde\Ag_2}{s^{3/2}_0}\right )\frac{P_{0,M}\left( \tilde h_4\chi (2y,s_0) \right)}{ P_{0,M}\left (i\chi (2y,s_0)\right)}\\
&&-\dsp\left (\frac{ \Bg_4 }{s^{3/2}_0}+\frac{ \Cg_4 \tilde\Ag_2}{s^{3/2}_0}\right )\frac{P_{0,M}\left( h_4\chi (2y,s_0) \right)}{ P_{0,M}\left (i\chi (2y,s_0)\right)}
\end{array}
\eeqtn
is well defined, where $\chi$ defined in \eqref{defchi14} and the constants $\tilde \Ag_i,\Ag_i$, $\tilde \Bg_i,\Bg_i$, $\tilde \Cg_i,\Cg_i$, for $i=0,2,4$ are given by (\ref{defi-as-Q-4}$-$\ref{defi-as-tilde-Q-2}).\\
(ii) If $\psi$ is given by \eqref{definitdata1} and \eqref{d2} with $d_0$ defined by  \eqref{d2}.Then, there exists a quadrilateral $\Dg_{s_0}\subset [-2,2]^2$ such that the mapping 
\[\dsp (\tilde d_0,\tilde d_1)\to \left (\tilde \Psi_0=\tilde \psi_0-\left (\frac{\tilde \Ag_{0}}{s_0}+\frac{\tilde \Bg_0 }{s^{3/2}_0}+\frac{\tilde \Cg_0 \tilde\Ag_2}{s^{3/2}_0}\right ),\tilde \psi_1\right )\]

 (where $\psi$ stands for $\psi_{s_0,\tilde d_0,\tilde d_1}$) is linear, one to one from $\Dg_{s_0}$ onto $[-\frac{A}{s_{0}^{7/4}},\frac{A}{s_{0}^{7/4}}]\times[-\frac{A}{s_{0}^{3/2}},\frac{A}{s_{0}^{3/2}}]$. Moreover it is of degree $1$ on the boundary.\\
(iii) For all $(\tilde d_0,\tilde d_1)\in \Dg_{s_0}$, $\psi_e\equiv 0$, $\psi_0=0$, $ |\tilde\psi_i|+|\psi_j|\leq C A e^{-\gamma s_0}$

 for some $\gamma>0$, for some $\gamma>0$ and for all $3\leq i\leq M,\;i\not =4$ and $1\leq j \leq M,\;j\not =4$ and 
 \[ |\tilde\Psi_i|+|\Psi_j|\leq C A e^{-\gamma s_0}\mbox{ for }i,j=\{2,4\},\]
 where $\tilde\Psi_i$ and $\Psi_i$ are defined as in (\ref{defi-as-Q-4}$-$\ref{defi-as-tilde-Q-2}).\\
Moreover , $\|\frac{\psi_-(y)}{(1+|y|)^{M+1}}\|_{L^\infty(\R)}\leq C\frac{A}{s_{0}^{\frac M4+1}}$.\\
(iv) For all $(\tilde d_0,\tilde d_1)\in \Dg_{s_0}$, $\psi_{s_0,\tilde d_0,\tilde d_1}\in \Vg_A(s_0)$ with strict inequalities except for $(\tilde \psi_0,\tilde\psi_1)$.

\label{propinitialdata}
\end{prop}
The proof of previous proposition is postponed to subsection \ref{tecnicsect}. 

\medskip

In the following, we find a local in time solution for equation \eqref{eqq} coupled with the condition \eqref{eqmod}.

\begin{prop}(\textbf{Local in time solution and modulation for problem \eqref{eqq}-\eqref{eqmod} with initial data  \eqref{definitdata1}-\eqref{d2}})\label{enopropmod}
For all $A\geq 1$, there exists $T_3(A)\in (0,1/e)$ such that for all $T\leq T_3$, the following holds:\\
For all $(\tilde d_0,\tilde d_1)\in D_T$, there exists $s_{max}>s_0=-\log T$ such that problem \eqref{eqq}-\eqref{eqmod} with initial data at $s=s_0$,
\[(q(s_0),\theta(s_0))=(\psi_{s_0,\tilde d_0,\tilde d_1},0),\]
where $\psi_{s_0,\tilde d_0,\tilde d_1}$ is given by \eqref{definitdata1} and \eqref{d2}, has a unique solution $q(s),\theta(s)$ satisfying $q(s)\in V_{A+1}(s)$ for all $s\in[s_0,s_{max})$.
\end{prop}
The proof of this proposition will be given later in page \pageref{prolocalintime}.

\medskip

Let us now give the proof of Proposition \ref{propsop1}.\\
\textit{Proof of Proposition \ref{propsop1}}: Let us consider $A\geq 1$, $s_0\geq s_{03}$, $(\tilde d_0,\tilde d_1)\in \Dg_{s_0}$, where $s_{03}$ is given by Proposition \ref{propinitialdata}. From the existence theory (which follows from the Cauchy problem for equation \eqref{GL}), starting in $\Vg_A(s_0)$ which is in $\Vg_{A+1}(s_0)$, the solution stays in $\Vg_A(s)$ until some maximal time $s_*=s_*(\tilde d_0,\tilde d_1)$. Then, either:\\
$\bullet$ $s_*(\tilde d_0,\tilde d_1)=\infty$ for some $(\tilde d_0,\tilde d_1)\in \Dg_{s_0}$, then the proof is complete.\\
$\bullet$ $s_*(\tilde d_0,\tilde d_1)<\infty$, for any  $(\tilde d_0,\tilde d_1)\in \Dg_{s_0}$, then we argue by contradiction. By continuity and the definition of $s_*$, the solution on $s_*$ is in the boundary of $\Vg_A(s_*)$. Then, by definition of $\Vg_A(s_*)$, one at least of the inequalities in that definition is an equality. Owing to the following proposition, this can happen only for the first two components $\tilde q_0,\tilde q_1$. Precisely we have the following result

\begin{prop}[Control of $q(s)$ by $(q_0(s),q_1(s))$ in $\Vg_A(s)$]. There exists $A_4\geq 1$ such that for each $A\geq A_4$, there exists $s_{04}\in \R$ such that for all $s_0\geq s_{04}$. The following holds:\\
If $q$ is a solution of \eqref{eqq} with initial data at $s=s_0$ given by \eqref{definitdata1} and \eqref{d2} with $(\tilde d_0,\tilde d_1)\in \Dg_{s_0}$, and $q(s)\in \Vg(A)(s)$ for all $s\in [s_0,s_1]$, with $q(s_1)\in \pa \Vg_A(s_1)$ for some $s_1\geq s_0$, then:\\
(i)\textbf{(Smallness of the modulation parameter $\theta$ defined in \eqref{formulav})} For all $s\in [s_0,s_1]$, 
\[|\theta '(s) |\leq \frac{CA^{10}}{s^{\frac 54}}.\] 
(ii) \textbf{(Reduction to a finite dimensional problem)} We have: 
\[\left (\tilde Q_0(s_1),\tilde q_1(s_1)\right )\in \pa\left( \left [-\frac{A}{s_{1}^{\frac 74}},\frac{A}{s_{1}^{\frac 74}} \right ]\times\left [-\frac{A}{s_{1}^{\frac 32}},\frac{A}{s_{1}^{\frac 32}} \right ]\right).\]
(iii)\textbf{(Transverse crossing)} There exists $\omega\in \{-1,1\}$ such that
\[\omega \tilde Q_0(s_1)=\frac{A}{s_{1}^{\frac 74}}\mbox{ and }\omega \frac{d\tilde Q_0(s_1)}{ds} (s_1)>0.\]
\[\omega \tilde q_1 (s_1)=\frac{A}{s_{1}^{\frac 32}}\mbox{ and }\omega \frac{d\tilde q_1}{ds} (s_1)>0.\]

\label{propcontrol}
\end{prop}
\label{indextheory}
Assume the result of the previous proposition, for which the proof is given below in page \pageref{proofiniti}, and continue the proof of Proposition \ref{propsop1}. Let $A\geq A_4$ and $s_0\geq s_{04}(A)$. It follows from Proposition \ref{propcontrol}, part (ii) that $\left ( \tilde Q_0,\tilde q_1(s_*)\right)\in\pa\left( \left [-\frac{A}{s_{1}^{\frac 74}},\frac{A}{s_{1}^{\frac 74}} \right ]\times\left [-\frac{A}{s_{1}^{\frac 32}},\frac{A}{s_{1}^{\frac 32}} \right ]\right)$, and the following function
\[
\begin{array}{ll}
\phi&:\Dg_{s_0}\to \pa ([-1,1]^2)\\
&(\tilde d_0,\tilde d_1)\to \left( \frac{s_{*}^{7/4}}{A}\tilde Q_0,\frac{s_{*}^{3/2}}{A}\tilde q_1\right)_{(\tilde d_0,\tilde d_1)}(s_*)\mbox{, with }s_*=s_*(\tilde d_0,\tilde d_1),
\end{array}
\]
is well defined. Then, it follows from Proposition \ref{propcontrol}, part (iii) that $\phi$ is continuous. On the other hand, using Proposition \ref{propinitialdata} (ii)-(iv)
together with the fact that $q(s_0)=\psi_{s_0,\tilde d_0,\tilde d_1}$, we see that when $(\tilde d_0,\tilde d_1)$ is in the boundary of the rectangle $\Dg_{s_0}$, we have strict inequalities for the other components.\\
Applying the transverse crossing property given by (iii) of Proposition \ref{propcontrol}, we see that $q(s)$ leaves $\Vg_A(s)$ at $s=s_0$, hence $s_*(\tilde d_0,\tilde d_1)=s_0$. Using Proposition \ref{propinitialdata}, part (ii), we see that the restriction of $\phi$ to the boundary is of degree 1. A contradiction, then follows from the index theory. Thus there exists a value $(\tilde d_0,\tilde d_1)\in \Dg_{s_0}$ such that for all $s\geq s_{0}$, $q_{s_0,d_0,d_1}(s)\in \Vg_A(s)$. This concludes the proof of Proposition \ref{propsop1}.\\
Using (i) of Proposition \ref{propcontrol}, we get the bound on $\theta'(s)$. This concludes the proof of \eqref{limitev}.

\subsection{Proof of the technical results}\label{tecnicsect}
This section is devoted to the proof of the existence result given by Theorem \ref{thm1}. We proceed in 4 steps, each of them making a separate 
subsection.
\begin{itemize}
\item In the first subsection, we give some properties of the shrinking set $\Vg_A(s)$ defined in Definition \ref{defthess} and translate our goal of making $q(s)$ go to $0$ in $L^\infty(\R)$ in terms of belonging to $\Vg_A(s)$. We also give the proof of Proposition \ref{propinitialdata}.

\item In second subsection, we solve the local in time Cauchy problem for equation \eqref{eqq} coupled with some orthogonality condition.

\item In the third subsection using the spectral properties of equation \eqref{eqq}, we reduce our goal from the control of $q(s)$ (an infinite dimensional variable) in $\Vg_A(s)$ to control its two first components ($\tilde Q_0$,$\tilde q_1$) a two variables in $\left [-\frac{A}{s^{\frac 74}},\frac{A}{s^{\frac 74}} \right ]\times\left[-\frac{A}{s^{\frac 32}},\frac{A}{s^{\frac32} } \right]$.

\item  In the fourth subsection, we solve the finite dimensional problem using the index theory and conclude the proof of Theorem \ref{thm1}
.

\end{itemize}

\bigskip

\subsubsection{Properties of the shrinking set $\Vg_A(s)$ and preparation of initial data}
In this subsection, we give some properties of the shrinking set defined in Definition \ref{defthess}. Let us first introduce the following claim:
\begin{cl}[Properties of the shrinking set defined in Definition \ref{defthess}]
 For all $r\in \Vg_A(s)$,

\medskip

\noindent (i) $\|r\|_{L^\infty(|y|<2K s^{\frac 14})}\leq C(K) \frac{A^{M+1}}{s^{\frac 14}}$ and $\|r\|_{L^\infty}(\R)\leq C(K)\frac{A^{M+2}}{s^{\frac 14}}$. \\
\noindent (ii) for all $y\in \R$, $|r(y)|\leq C \frac{A^{M+1}}{s}(1+|y|^{M+1})$. 
\label{propshrinset}
\end{cl}
\textit{Proof: }Take $r\in \Vg_A(s)$ and $y\in \R$.\\
(i) If $|y|\geq 2 K s^{\frac 14}$, then we have from the definition of  $r_e$  \eqref{defiqe}, $|r(y)|=|r_e(y)|\leq \frac{A^{M+2}}{s^{\frac 14}} $.\\
Now, if $|y| < 2 K s^{\frac 14}$, since we have for all $0\leq j\leq M$, $|\tilde r_j|+| r_j|\leq C\frac{A^j}{s^{\frac{j+1}{4}}}$ from Definition \ref{defthess} (use the fact that $M\geq 4$), we write from \eqref{decompq}
\beqtn
\begin{array}{ll}
|r(y)|&\leq \dsp\left(\sum_{j\leq M}    |\tilde r_j||\tilde h_j|+ | r_j|| h_j|\right )+|r_-(y)|,\\
&\leq C\dsp\sum_{j\leq M}  \frac{A^{M+1}}{s^{\frac{j+1}{4}}}(1+|y|)^j+ \frac{A^{M+1}}{s^{\frac{M+2}{4}}}(1+|y|)^{M+1},\\
& \leq C\dsp\sum_{j\leq M}  \frac{A^{M+1}}{s^{\frac{j+1}{4}}}(1+K s^{\frac 14})^j+ \frac{A^{M+1}}{s^{\frac{M+2}{4}}}(1+K s^{\frac 14})^{M+1} \leq C \frac{(KA)^{M+1}}{s^{\frac 14}},
\end{array}
\label{estipss}
\eeqtn
which gives (i).\\
(ii) Just use \eqref{estipss} together with the fact that for all $0\leq j\leq M$, $|\tilde r_j|+|r_j|\leq C\frac{A^{M+1}}{s}$ from Definition \ref{defthess}. This ends the proof of Claim \ref{propshrinset}. $\blacksquare$

\medskip

Let us now give the proof of Proposition \ref{propinitialdata}.\\
\textit{Proof of Proposition \ref{propinitialdata}} For simplicity, we write $\psi$ instead of  $\psi_{s_0,\tilde d_0,\tilde d_1}$. We note that, from Claim \ref{propshrinset}, (iv) follows from (ii) and (iii) by taking $s_0=-\log T$ large enough (that is $T$ is small enough). Thus, we only prove  (i), (ii) and (iii). Consider $K\geq 1$, $A\geq 1$ and $T\leq 1/e$. Note that $s_0=-\log T\geq 1$.\\
The proof of (i) is a direct consequence of (iii) of the following claim
 \begin{cl}
There exists $\gamma=\frac{1}{32(1+\beta^2)}>0$ and $T_2<1/e$ such that for all $K\geq 1$ and $T\leq T_2$, if $g$ is given by $(1+i\delta)\chi(2y,s_0)$, $(1+i\delta)y\chi(2y,s_0)$, $(1+i\delta)h_2(y)\chi(2y,s_0)$ or $i\chi(2y,s_0)$, then $\left\|\frac{g_-(y)}{1+|y|^{M+1}}\right\|_{L^\infty}\leq \frac{C}{s_{0}^{\frac M4}}$ and all $ g_i$, $\tilde g_i$ for $0\leq i\leq M$ are less than $Ce^{-\gamma s_0}$. expect:\\
i) $ |\tilde g_0-1|\leq C e^{-\gamma s_0}$ when $g=\tilde h_0(y)\chi(2y,s_0)$.\\
ii)  $|\tilde g_1-1|\leq C e^{-\gamma s_0}$ when $g=\tilde h_1(y)\chi(2y,s_0)$.\\
iii)  $|\tilde g_2-1|\leq C e^{-\gamma s_0}$ when $g=\tilde h_2(y)\chi(2y,s_0)$.\\
iv)  $| g_0-1|\leq C e^{-\gamma s_0}$ when $g=h_0(y)\chi(2y,s_0)$.\\
v)$| g_2-1|\leq C e^{-\gamma s_0}$ when $g=h_2(y)\chi(2y,s_0)$.\\
vi)$|\tilde g_4-1|\leq C e^{-\gamma s_0}$ when $g=\tilde h_4(y)\chi(2y,s_0)$.\\
vii)$| g_4-1|\leq C e^{-\gamma s_0}$ when $g=h_4(y)\chi(2y,s_0)$.
\label{cl}
\end{cl}
\textit{Proof: } In all cases, we write 
\beqtn
g(y)=p(y)+r(y)\mbox{ where }p(y)=\tilde h_{j|j=0,1,2,4} \mbox{ or } h_{j|j=0,2,4} \mbox{ and }r(y)=p(y)(\chi(2y,s_0)-1).
\label{pcl1}
\eeqtn
From the uniqueness of the decomposition \eqref{decompq}, we see that $p_-\equiv 0$ and al $ p_i$, $\tilde p_i$ are zero except\\
\[\tilde p_j=1\mbox{, when }p(y)=\tilde h_j\mbox{, for }j=0,1,2\mbox{ and }4\] 
and
\[ p_j=1\mbox{, when }p(y)= h_j\mbox{, for }j=0,2\mbox{ and }4\]


Concerning the cases $2|y|<K s^{\frac 14}$ and $2|y|>K s^{\frac 14}$, we have the definition of $\chi$  \eqref{defchi14},
\[1-\chi(2y,s)\leq \left(\frac{2|y|}{K s_{0}^{\frac 14}}\right)^{M-1},\]
\[|\rho_{\beta}(y)(1-\chi(2y,s))|\leq \sqrt{|\rho_{\beta}(y)|}\sqrt{\left |\rho_{\beta}\left(\frac K2s^{\frac 14}\right)\right |}\leq Ce^{-\frac{K^2 s_0}{32(1+\beta^2)}}\sqrt{|\rho_{\beta}(y)|}.\]
Therefore, from \eqref{decompq} and \eqref{pcl1}, we see that
\beqtn
\begin{array}{l}
|r(y)|\leq C(1+|y|^2)\left(\frac{2|y|}{Ks_{0}^{\frac 14}}\right)^{M-1}\leq C\frac{(1+|y|^{M+1})}{s_{0}^{\frac M4}},\\
|\Rg_j|+| r_j|+|\tilde r_j|\leq Ce^{-\frac{K^2\sqrt {s_0}}{32(1+\beta^2)}}\mbox{ for all }j\leq M.
\end{array}
\label{pcl2}
\eeqtn
Hence, using  \eqref{pcl2} and \eqref{decomp1} and the fact that $|f_j(y)|\leq C (1+|y|)^M$, for all $j\leq M$, we get also 
\[|r_-(y)|\leq C\frac{(1+|y|)^M}{s_{0}^{\frac M4}}.\]
Using \eqref{pcl1} and the estimates for $p(y)$ stated below, we conclude the proof of Claim \ref{cl} and (i) of Proposition \ref{propinitialdata}. \\
(ii) of Proposition \ref{propinitialdata}: From \eqref{definitdata1} and \eqref{d2}, we see that
\beqtn
\left (
\begin{array}{l}
\tilde \Psi_0\\
\tilde \psi_1
\end{array}\right )
=G
\left (\begin{array}{l}
\tilde d_0\\\tilde d_1
\end{array}
\right )
\mbox{ where }G=(g_{i,j})_{0\leq i,j\leq 1}.
\eeqtn
Using Claim \ref{cl}, we see from \eqref{definitdata1} and \eqref{d2} that
\beqtn
|d_0|\leq C(|\tilde d_0|+|\tilde d_1|)e^{-\gamma s_0}
\label{initiad3}
\eeqtn
for $T$ small enough. Using again Claim \ref{cl}. We see that $G\to\left( \begin{array}{ll} \frac{A}{s_{0}^{\frac 74}}&0\\0&\frac{A}{s_{0}^{\frac 32}}  \end{array}\right )$ and

as $s_0\to\infty$ (for fixed $K$ and $A$), which concludes the proof of (ii) of Proposition \ref{propinitialdata}.\\
(iii) of Proposition \ref{propinitialdata}: Since $supp (\psi)\subset B(0,K s_{0}^{\frac 14})$ by  \eqref{definitdata1} and \eqref{d2}, we see that $\psi_e\equiv 0$ and that $\psi_0$ is zero from the definition of $d_0$  \eqref{definitdata1} and \eqref{d2}. Using the fact that $|\tilde d_{i,i=0,1}|\leq 2$ and the bound on $d_0$ by \eqref{initiad3}, we see that the estimates on $ \psi_j$ and $\tilde \psi_j$ and $\psi_-$ in (iii) follows from \eqref{definitdata1} and \eqref{d2} and Claim \ref{cl}. This concludes the proof of Proposition \ref{propinitialdata}. $\blacksquare$

\medskip

In the following we give the proof of Local in time solution for problem \eqref{eqq}-\eqref{eqmod}. In fact, we impose some orthogonality condition given by \eqref{eqmod}, killing the one of the zero eigenfunction of the linearized operator of equation \eqref{eqq}.

\medskip

\textit{Proof of Proposition \ref{enopropmod}:}\label{prolocalintime} From solution of the local in time Cauchy problem for equation \eqref{GL} in $L^\infty(\R)$, there exists $s_1>s_0$ such that equation \eqref{equa-w} with initial data (at $s=s_0$) $\varphi(y,s_0)+ \psi_{s_0,\tilde d_0,\tilde d_1}(y)$, where $\varphi(y,s)$ is given by \eqref{defi-varphi} has a unique solution $w(s)\in C([s_0,s_1),L^\infty(\R))$. Now, we have to find a unique $(q(s), \theta(s))$ such that
\beqtn
w(y,s)=e^{i(\nu \sqrt s+\mu\log s+\theta(s))}\left(\varphi(y,s)+q(y,s)\right)
\label{promod1}
\eeqtn
and \eqref{eqmod} is satisfied. Since $f_0=1$ and $\int_{\R}\rho_{\beta}(y)dy=1$, we use \eqref{decomp3} to write  \eqref{eqmod} as follows
\begin{eqnarray*}
P_{0,M} (q) &=&\Im \left(\int q(y,s)\rho_{\beta}(y)dy\right )-\delta \Re \left(\int q(y,s)\rho_{\beta}(y)dy\right )\\
&=&\Im \left((1-i\delta)\int q(y,s)\rho_{\beta}(y)dy\right)=0,
\end{eqnarray*}
or using \eqref{promod1}
\[F(s,\theta)\equiv \Im\left((1-i\delta)\int \left(e^{-i(\nu \sqrt s+\mu\log s+\theta(s))}  w(y,s)-\varphi(y,s)\right)\rho_{\beta}(y)dy \right)=0.\]
Note that
\[ \frac{\pa F}{\pa\theta} (s,\theta)=-\Re\left((1-i\delta)\int e^{-i(\nu \sqrt s+\mu\log s+\theta(s))}  w(y,s) \rho_{\beta}(y)dy\right).\]
From (iii) in Proposition \ref{propinitialdata}, $F(s_0,0)=P_{0,M}(\psi_{s_0,\tilde d_0,\tilde d_1})=0$ and 
\begin{eqnarray*}
\frac{\pa F}{\pa\theta} (s_0,0)&=&-\Re \left((1-i\delta)\int (\varphi(y,s_0)+\psi_{s_0,\tilde d_0,\tilde d_1}(y)  )\rho_{\beta}(y)dy \right) \\
&=&-\kappa +O\left(\frac{1}{s_{0}^{1/4}}\right )\mbox{ as }s_0\to \infty,
\end{eqnarray*}
for fixed $K$ and $A$.\\
Therefore, if $T$ is small enough in terms of $A$, then $\frac{\pa F}{\pa \theta}(s_0,0)\not =0$, and from the implicit function Theorem, there exists $s_2\in(s_0,s_1)$ and $\theta\in C^1([s_0,s_2),\R)$ such that $F(s,\theta(s))=0$ for all $s\in [s_0,s_2)$.Defining $q(s)$ by \eqref{promod1} gives a unique solution of the problem  \eqref{eqqd}-\eqref{eqmod} for all $s\in[s_0,s_2)$.  Now, since we have from (iv) of Proposition \ref{propinitialdata}, $q(s_0)\in V_{A}(s_0)  \begin{array}{l}\subset\\ \not =\end{array} V_{A+1}(s_0)  $, there exists $s_3\in (s_0,s_2)$ such that for all $s\in [s_0,s_3)$, $q(s)\in V_{A+1}(s)$. This concludes the proof of Proposition \ref{enopropmod}. $\blacksquare$

\subsubsection{Reduction to a finite dimensional problem}
In the following we give the proof of Proposition  \ref{propcontrol}:\\
\label{proofiniti}\\
The idea of the proof is to project equation  \eqref{eqq} on the different components of the decomposition \eqref{decompq}. More precisely, we claim that Proposition \ref{propcontrol}  is a consequence of the following

\begin{prop} There exists $A_5\geq 1$ such that for all $A\geq A_5$, there exists $s_5(A)$ such that the following holds for all $s_0\geq s_5$:\\
Assuming that for all $s\in [\tau,s_1]$ for some $s_1\geq \tau\geq s_0$, $q(s)\in \Vg_A(s)$ and $ q_0(s)=0$, then  the following holds for all $s\in [\tau, s_1]$:\\

\noindent (i) (\textbf{Smallness of the modulation parameter}):
\[|\theta '(s)|\leq \frac{ CA^{10}}{s^{\frac{5}{4}}}.\]
(ii) (\textbf{ODE satisfied by the expanding mode}): For $m=0$ and $1$, we have
\begin{eqnarray*}
\left|  \tilde Q_0'(s)   - Q_0(s)    \right|  \leq \frac{C}{s^\frac{7}{4}},
\end{eqnarray*}
and
\begin{eqnarray*}
\left|   \tilde q_1'  - \frac{1}{2} \tilde q_1   \right| \leq \frac{C}{s^\frac{3}{2}}.
\end{eqnarray*}
(iii) (\textbf{ODE satisfied by the null mode}):
$$  \left|\tilde Q_2'(s) - \tilde H_1  \frac{\tilde Q_2}{s} \right|  \leq    \frac{C A^8}{s^\frac{9}{4}},  $$
where $\tilde H_1$ is a constant depending only on $p, \delta$ and less than $ -\frac{3}{2}$. 

 (iv) (\textbf{Control of negative modes}):
\[| q_1(s)|\leq e^{-\frac{(s-\tau)}{2}}| q_1(\tau)|+\frac{CA^3}{s^{\frac 32}},\]
\[   \left| Q_2 (s)    \right| \leq     e^{-(s -\tau)}  \left| Q_2 (\tau)  \right|  + \frac{CA^7}{s^\frac{7}{4}}  ,\]
\[ | q_3|  \leq  e^{-\frac{3}{2}(s-\tau)} |q_3(\tau)| + \frac{ CA^2}{s^\frac{3}{2}}, \]
\[ | \tilde q_3|  \leq  e^{-\frac{s -\tau}{2}} | \tilde q_3(\tau)| + \frac{ CA^2}{s^\frac{3}{2}}, \]
\[  \left| Q_4 ( s)  \right|  \leq  e^{- 2(s -\tau)}  \left| Q_4(\tau)  \right| + \frac{CA^6}{s^\frac{7}{4}},      \]
\[   \left| \tilde Q_4 ( s)  \right|  \leq  e^{- (s -\tau)}  \left| \tilde  Q_4(\tau)  \right| + \frac{CA^3}{s^\frac{7}{4}},     \]
\[| q_j(s)|\leq e^{-j\frac{(s-\tau)}{2}}| q_j(\tau)|+\frac{CA^{j-1}}{s^{\frac{j+1}{4}}},\mbox{ for all } 5\leq j\leq M,\]
\[|\tilde q_j(s)|\leq e^{-(j-2)\frac{(s-\tau)}{2}}|\tilde q_j(\tau)|+\frac{CA^{j-1}}{s^{\frac{j+1}{4}}},\mbox{ for all }5\leq j\leq M,\]
\[ \left\|\frac{q_-(y,s)}{1+|y|^{M+1}} \right\|_{L^\infty} \leq e^{-\frac{M+1}{4}(s-\tau)}\left\|\frac{q_-(\tau)}{1+|y|^{M+1}}  \right\|_{L^\infty} +C\frac{A^M}{s^{\frac{M+2}{4}}}, \]
\[\|q_e(y,s)\|_{L^\infty}\leq e^{-\frac{(s-\tau)}{2(p-1)}}\|q_e(\tau)\|_{L^\infty}  +\frac{C A^{M+1}}{\tau^{\frac 14}}(1+s-\tau),\]
\noindent
where $\tilde Q_0$, $Q_2$, $\tilde Q_2$, $Q_4$ and $\tilde Q_4$ are defined by (\ref{defi-as-Q-4}-\ref{defi-as-tilde-Q-2}).
\medskip
\label{propode}
\end{prop}
The idea of the proof of Proposition \ref{propode} is to project equations \eqref{eqq1} and \eqref{eqq} according to the decomposition \eqref{decompq}. However because of the number of parameters and coordinates in \eqref{decompq}, the computation become too long. That is why Subsection \ref{proofpropode} is devoted to the proof of Proposition \ref{propode}. 

\medskip

Let us now derive Proposition \ref{propcontrol} from Proposition \ref{propode}.

\medskip

\textit{Proof of Proposition \ref{propcontrol} assuming Proposition \ref{propode}:}\\
We will take $A_4\geq A_5$. Hence, we can use the conclusion of Proposition \ref{propode}.\\
(i) The proof follows from (i) of Proposition \ref{propode}. Indeed by choosing $T_4$ small enough, we can make $s_0=-\log T$ bigger than $s_5(A)$.\\

(ii) We notice that from Claim \ref{propshrinset} and the fact that $ q_0(s)=0$, it is enough to prove that for all $s\in [s_0,s_1]$,

\beqtn\label{goalnullmode}
\left| \tilde Q_2(s)  \right|=\left| \tilde q_2(s)  - \frac{\mathcal{A}_2}{s}\right|<\frac{A^{10}}{s^{\frac 54}}.
\eeqtn

\beqtn
\begin{array}{ll}

\|q_e\|_{L^\infty(\R)}\leq \frac{A^{M+2}}{2s^{\frac 14}}, &\left\|\frac{q_-(y,s)}{1+|y|^{M+1}} \right\|_{L^\infty} \leq \frac{A^{M+1}}{2 s^{\frac{M+2}{4}}},\\
| q_j|, |\tilde q_j|\leq \frac{A^j}{2s^{\frac{j+1}{4}}}\mbox{ for all }5\leq j\leq M,& | q_1|\leq \frac{A^4}{2s^{\frac 32}},\;\;\;\;\;\;\\
 | Q_2|\leq \frac{A^8}{2s^{\frac 74}}& |q_3|, |\tilde q_3| \leq \frac{A^3}{2s^\frac{3}{2}},
 \\|Q_4|, |\tilde Q_4| \leq \frac{A^7}{2s^{\frac 74}}.&

\end{array}
\label{object}
\eeqtn
Let us first prove \eqref{goalnullmode}: Indeed,  we will use a contradictory  argument, we assume that there exists $s_* \in [s_0,s_1]$ such that 
\[ \tilde Q_2(s_*) =\left( \tilde q_2(s_*) - \frac{\mathcal{A}_2}{s_*} \right)=\omega \frac{A^{10}}{s_*^{5/4}} \mbox{ and for all } s \in [s_0,s^*[,\;\;  \left|  \tilde q_2(s) - \frac{\mathcal{A}_2}{s} \right|<\frac{A^{10}}{s^{5/4}},\]
\noindent
where $\omega = \pm 1$. As a matter of fact,  we can reduce to the  positive case where $\omega=1$ (the case $\omega = -1$ also work by the same way). Note by item (iv) in Proposition \ref{propinitialdata} that 
 $$  \left| \tilde q_2(s_0) - \frac{ \mathcal{A}_2}{s_0} \right| < \frac{A^{10}}{s_0^{\frac 54}},$$
thus $s_* >s_0$, and the interval $[s_0,s_*]$ is not empty.\\
Using the continuity of $\tilde Q_2$ and the definition of $s_*$, it is clearly that   $ \tilde Q_2(s_*)$ is the maximal value of $\tilde Q_2$ in $[s_* - \epsilon, s_*]$ with $\epsilon >0$ and small enough   
 in one hand, recalling, from (iii) Proposition \ref{propode} that
 \[
\left |\tilde Q_2' -   \tilde H_1\frac{\tilde Q_2}{s}\right |\leq \frac{CA^8}{s^{\frac94}},
\]
and from \eqref{sign-tilde-H-1} that $\tilde H_1\leq-\frac 32 $, we write
\beqtn
\tilde Q_2'(s_*)\leq  \tilde H_1\frac{\tilde Q_2}{s}+\frac{A^8}{s^{9/4}}\leq\frac{-3/2 A^{10}+CA^8}{s^{9/4}} < 0,
\label{argument-Q2-2}
\eeqtn
 for $A$ large enough. Then, $\tilde Q_2 $ has to decrease in $[s_* - \epsilon_1, s_*]$  which implies  a contradiction   with the assumption  that $\tilde Q_2$ admits  maximum  at $s_*$. In other word, \eqref{goalnullmode} holds.





\medskip

Now, let us deal with \eqref{object}. Define $\sigma =\log A$ and take $s_0\geq \sigma$ (that is $T\leq e^{-\sigma}=1/A$) so that for all $\tau \geq s_0$ and $s\in [\tau, \tau +\sigma]$, we have 
\beqtn
\tau \leq s\leq \tau+\sigma\leq \tau +s_0\leq 2\tau\mbox{ hence }\frac{1}{2\tau} \leq \frac 1s\leq \frac 1\tau\leq \frac{2}{s}. 
\label{taus}
\eeqtn
We consider two cases in the proof.\\
\textbf{ Case 1:  $s\leq s_0+\sigma$.}\\
Note that \eqref{taus} holds with $\tau=s_0$. Using (iv) of Proposition \ref{propode} and estimate (iii) of Proposition \ref{propinitialdata} on the initial data $q(.,s_0)$ (where we use \eqref{taus} with $\tau =s_0$), we write
\beqtn
\begin{array}{l}
| q_1(s)|\leq C A e^{-\gamma_1 \frac s2}+ \frac{CA^3}{s^{3/2}},\\
| Q_2(s)|\leq  C A e^{-\gamma_1 \frac s2}+ \frac{CA^7}{s^{\frac 74}},\\
|q_3| \leq C A e^{-\gamma_1 \frac s2}+ \frac{CA^2}{s^{\frac 32}},\\
|\tilde q_3| \leq C A e^{-\gamma_1 \frac s2}+ \frac{CA^2}{s^{\frac 32}},\\
| Q_4(s)|\leq  C A e^{-\gamma_1 \frac s2}+ \frac{CA^6}{s^{\frac 74}},\\
| \tilde Q_4(s)|\leq  C A e^{-\gamma_1 \frac s2}+ \frac{CA^3}{s^{\frac 74}},\\
|\tilde q_j(s)|\leq CA e^{-\gamma_1 \frac s2}+ \frac{C A^{j-1}}{s^{\frac{j+1}{4}}}\mbox{ for all }3\leq j\leq M,\;\;j\not =4\\
| q_j(s)|\leq CA e^{-\gamma_1 \frac s2}+ \frac{C A^{j-1}}{s^{\frac{j+1}{4}}}\mbox{ for all }3\leq j\leq M,\;\;j\not =4\\
\left\|\frac{q_-(s)}{1+|y|^{M+1}}\right\|_{L^\infty}\leq C\frac{A}{\left(\frac s2\right)^{\frac M4+2}}+C\frac{A^M}{s^{\frac{M+2}{4}}},\\
\|q_e(s)\|_{L^\infty}\leq \frac{CA^{M+1}}{\left(\frac s2\right)^{\frac 14}}(1+\log A).
\end{array}
\eeqtn
Thus, if $A\geq A_6$ and $s_0\geq s_6(A)$ (that is $T\leq e^{-s_6(A)}$) for some positive $A_6$ and $s_6(A)$, we see that \eqref{object} holds.\\
\textbf{ Case 2:  $s> s_0+\sigma$.}\\
Let $\tau=s-\sigma>s_0$. Applying (iv) of Proposition \ref{propode} and using the fact that $q(\tau)\in \Vg_A(\tau)$, we write (we use \eqref{taus} to bound any function of $\tau$ by a function of $s$)
\beqtn
\begin{array}{l}
| q_1(s)|\leq  e^{-\frac\sigma 2}\frac{A^6}{\left(\frac s2\right)^{3/2}}+ \frac{CA^3}{s^{3/2}},\\
| Q_2(s)|\leq e^{-\sigma}\frac{A^8}{\left(\frac s2\right)^{7/4}}+ \frac{CA^7}{s^{7/4}},\\
| q_3(s)|\leq  e^{-\frac{3\sigma}{ 2}}\frac{A^3}{\left(\frac s2\right)^{3/2}}+ \frac{CA^2}{s^{3/2}},\\
| \tilde q_3(s)|\leq  e^{-\frac{\sigma}{ 2}}\frac{A^3}{\left(\frac s2\right)^{3/2}}+ \frac{CA^2}{s^{3/2}},\\
| Q_4(s)|\leq e^{-2\sigma}\frac{A^7}{\left(\frac s2\right)^{7/4}}+ \frac{CA^6}{s^{7/4}},\\
|\tilde Q_4(s)|\leq e^{-\sigma}\frac{A^4}{\left(\frac s2\right)^{7/4}}+ \frac{CA^3}{s^{7/4}},\\
|\tilde q_j(s)|\leq e^{-\frac{(j-2)\sigma}{2}}\frac{A^j}{\left(\frac s2\right )^{\frac{j+1}{4}}}+ \frac{C A^{j-1}}{s^{\frac{j+1}{4}}}\mbox{ for all }5\leq j\leq M,\\
| q_j(s)|\leq e^{-\frac{j\sigma}{2}}\frac{A^j}{\left(\frac s2\right )^{\frac{j+1}{4}}}+ \frac{C A^{j-1}}{s^{\frac{j+1}{4}}}\mbox{ for all }5\leq j\leq M,\\
\left\|\frac{q_-(s)}{1+|y|^{M+1}}\right\|_{L^\infty}\leq e^{-\frac{M+1}{4}\sigma}\frac{A^{M+1}}{\left (\frac s 2\right)^{\frac{M+2}{4}}}+C\frac{A^M}{s^{\frac{M+2}{4}}},\\
\|q_e(s)\|_{L^\infty}\leq e^{-\frac{\sigma}{2(p-1)}}\frac{A^{M+2}}{\left (\frac s2\right)^{\frac 14}}+ \frac{CA^{M+1}}{\left(\frac s2\right)^{\frac 14}}(1+\sigma).
\end{array}
\eeqtn
For all the coordinates, it is clear that if $A\geq A_7$ and $s_0\geq s_7(A)$ for some positive $A_7$ and $s_7(A)$, then \eqref{goalnullmode} and \eqref{object} is satisfied (remember that $\sigma=\log A$).

\medskip

Conclusion of (ii): If $A\geq \max (A_6,A_7,A_8)$ and $s_0\geq \max (s_6(A), s_7(A),s_8(A))$, then \eqref{object} is satisfied. Since we know that $q(s_1)\in \pa V_A(s_1)$, we see from the definition of $V_A(s)$ that $(\tilde Q_0(s_1),\tilde q_1(s_1)) \in \pa [-\frac{A}{s_{1}^{7/4}},\frac{A}{s_{1}^{7/4}}]\times[-\frac{A}{s_{1}^{3/2}},\frac{A}{s_{1}^{3/2}}]$. This concludes the proof of (ii) of Proposition \ref{propcontrol}.

\medskip

(iii) From (ii), there is $\omega=\pm 1$ such that $\tilde Q_0(s_1)=\omega \frac{A}{s_{1}^{7/4}}$ or $\tilde q_0(s_1)=\omega \frac{A}{s_{1}^{3/2}}$. \\

Using (ii) of Proposition \ref{propcontrol}, we see that 
\[\omega \tilde Q^{'}_{0}(s_1)\geq \omega \tilde Q_0(s_1)-\frac{C}{s_{1}^{7/4}},\]
or
\[\omega \tilde q^{'}_{1}(s_1)\geq \frac 12\omega \tilde q_1(s_1)-\frac{C}{s_{1}^{3/2}}.\]
Taking $A$ large enough gives $\omega \tilde Q_{0}^{'}(s_1)>0$ or $\omega \tilde q_{1}^{'}(s_1)>0$, and concludes the proof of Proposition \ref{propcontrol}. $\blacksquare$

\subsection{Proof of Proposition \ref{propode}}\label{proofpropode}
In this section, we prove Proposition \ref{propode}. We just have to project equations \eqref{eqq1} and \eqref{eqq} to get equations satisfied by the different coordinates of the decomposition \eqref{decompq}.  We proceed as Section 5 in \cite{MZ07}, taking into account the new scaling law $\frac{y}{s^{1/4}}$. We note that the projections of $V_1 q+V_2\bar q$, $B$ and $R^*$ in  \eqref{eqqd}, will need much more effort and this is due to the fact that we are dealing with the critical case when $\beta\neq 0$.

\medskip

More precisely, the proof will be carried out in 3 subsections
\begin{itemize}
\item In the first subsection, we deal with equation \eqref{eqq} to write equations satisfied by $\tilde q_j$ and $ q_j$. Then, we prove (i), (ii), (iii) and (iv) (expect the two last identities) of Proposition \ref{propode}.
\item In the second subsection, we first derive from equation \eqref{eqq} an equation satisfied by $q_-$ and prove the last but one identity in (iv) of Proposition \ref{propode}.
\item In the third subsection, we project equation \eqref{eqq1} (which is simpler than \eqref{eqq}) to write an equation satisfied by $q_e$ and prove the last identity in (iv) of Proposition \ref{propode}.
\end{itemize}
\subsubsection{The finite dimensional part $q_+$}
We proceed in 2 parts:
\begin{itemize}
\item In Part 1, we give the  details of projections of  equation \eqref{eqq} to get ODEs,  satisfied by modes  $\tilde q_j$ and $ q_j$. 
\item In Part 2, we prove (i), (ii) and (iii) of Proposition \ref{propode}, together with the estimates concerning $\tilde q_j$ and $ q_j$ in (iv).

\end{itemize}
\textbf{Part 1: The projection of equation \eqref{eqq} on the eigenfunction of the operator $ \Lg_{\beta,\delta}$ }
In the following, we will find the main contribution in the projections $\tilde P_{n,M}$ and $ P_{n,M}$ of the six terms appearing in equation \eqref{eqq}: $\pa_s q$, $\Lg_{\beta,\delta} q$, $-i \left(\frac{\nu}{2 \sqrt s}+\frac \mu s+\theta'(s) \right)q$, $V_1 q+V_2 \bar q$, $B(q,y,s)$ and $R^*(\theta',y,s)$. Most of the time, we give two estimates of error terms, depending on whether we use or not the fact that $q(s)\in  \Vg_A(s)$.\\
\textbf{First term: $\dsp\frac{\pa q}{\pa s}$.}\\ 
From \eqref{decomp3}, its projection on $\tilde h_n$ and $\hat h_n$ is $\tilde q_{n}'$ and $ q_{n}'$ respectively:
\beqtn
\tilde P_{n,M}\left(\frac{\pa q}{\pa s}\right)=\tilde q_n'\mbox{ and }P_{n,M}\left(\frac{\pa q}{\pa s}\right)= q_n'.
\eeqtn
\textbf{Second term: $\mathcal{L}_{\beta, \delta} q$,}  where $\mathcal{L}_{\beta, \delta} $ is defined as in \eqref{eqqd}.    We will use the following Lemma from \cite{MZ07}:
\begin{lemma}[Projection of $ \Lg_{\beta,\delta} $ on $\tilde h_n$ and $h_n$ for $n\leq M$]\label{lemma-projec-P-n-L-beta-delta} \hspace{-0.04cm}

$a)$ If $n\leq M-2$,then
\[\left | P_{n,M}( \Lg_{\beta,\delta} q)-\left (-\frac n2 q_n(s)+c_{n+2}\tilde q_{n+2} \right ) \right |\leq C\left \|     \frac{q_-}{1+|y|^{M+1}} \right \|_{L^{\infty}},\]
where $c_n$ is given in  Lemma \ref{lemma-Jordan-block's}. Moreover, we have the following

\medskip
If $M-1\leq n\leq M$, then 
\[\left | P_{n,M}( \Lg_{\beta,\delta}q) +\frac n 2q_n(s)\right |\leq C\left \|     \frac{q_-}{1+|y|^{M+1}} \right \|_{L^{\infty}}. \]

\medskip
$(b)$ If $n\leq M$, then the projection of $\Lg_{\beta,\delta}$ on $\tilde h_n$ satisfies

\[\left |\tilde P_{n,M}( \Lg_{\beta,\delta}q) -\left ( 1-\frac n 2\right )\tilde q_n(s)\right |\leq C\left \|     \frac{q_-}{1+|y|^{M+1}} \right \|_{L^{\infty}}.\]

\end{lemma}
\textit{Proof}: The proof is quiet the same as the proof of Lemma 5.1 in \cite{MZ07}.

\medskip

Using   Lemma \ref{lemma-projec-P-n-L-beta-delta} and   the fact that $q(s)\in V_A(s)$ (see Definition \ref{defthess})  in addition,  then the error estimates can be improved  as follows

\begin{corollary}
For all $A\geq 1$, there exists $s_9\geq 1$ such that for all $s\geq s_9(A)$, if $q(s)\in V_A(s)$, then:\\

a) For $n=0$, we have

\[\left | P_{0,M}( \Lg_{\beta,\delta}q)-c_2 \tilde q_2  \right |\leq C\frac{A^{M+1}}{s^{\frac{M+2}{4}}}.\]

b) For $1\leq n\leq M-1$, we have
\[\left | P_{n,M}( \Lg_{\beta,\delta}q) +\frac n 2q_n(s)\right |\leq C\frac{A^{n+2}}{s^{\frac{n+3}{2}}}.\]

In particular, we have  a  smaller bound for $P_{2,M}(\mathcal{L}_{\beta, \delta} q)$:
\begin{equation*}
\left|   P_{2,M}(\mathcal{L}_{\beta, \delta} q) + q_2 - c_4 \tilde q_4 \right| \leq \frac{A^{M+1}}{s^\frac{M+2}{4}}.
\end{equation*}

c) For $n=M$, we have 
\[\left | P_{M,M}( \Lg_{\beta,\delta}q) +\frac M 2q_M(s)\right |\leq C\frac{A^{M+1}}{s^{\frac{M+2}{2}}}.\]

d) For $0\leq n\leq M$, we have 
\[\left |\tilde P_{n,M}( \Lg_{\beta,\delta}q) -\left ( 1-\frac n 2\right )\tilde q_n(s)\right |\leq C\frac{A^{M+1}}{s^{\frac{M+1}{2}}}.\]

\end{corollary}

\textbf{Third term: $ -i \left(\frac{\nu}{2 \sqrt s}+\frac{\mu }{s} +\theta '(s) \right)q$.} 
It is enough to project $ iq$, from \eqref{decomp3}, we recall Lemma 5.3 from \cite{MZ07}:
\begin{lemma}[Projection of the term $-i(\frac{\nu}{2 \sqrt s}+\frac{\mu}{ s}+\theta '(s))q$ on $h_n$ and $\tilde h_n$ for $n \leq M$] Its projection on $h_n$ is given by

\begin{eqnarray*}
& & P_{n,M}\left (-i \left(\frac{\nu}{2 \sqrt s}+\frac \mu s+\theta '(s) \right)q\right ) \\
&=& -\left(\frac{\nu}{ 2 \sqrt s}+\frac \mu s+\theta ' (s) \right)\left ( \delta q_n+(1+\delta^2)\tilde q_n+\sum_{j=n+1}^{M} K_{n,j}q_j+L_{n,j}\tilde q_j  \right),
\end{eqnarray*}

where  $ K_{n,j}$  and $L_{n,j}$   defined by 
 \begin{eqnarray}\label{defKL}
 K_{n,j} =P_{n,M}(i h_j),\\
 L_{n,j} = P_{n,M}({i\tilde h}_j).
 \end{eqnarray}
Its projection on $\tilde h_n$ is given by

\[
\tilde P_{n,M}\left (-i \left(\frac{\nu}{ 2 \sqrt s}+\frac \mu s+\theta '(s) \right)q\right )= -\left(\frac{\nu}{ 2 \sqrt s}+\frac \mu s+\theta ' \right)\left ( - q_n-\delta\tilde q_n+\sum_{j=n+1}^{M}\tilde K_{n,j}q_j+\tilde L_{n,j}\tilde q_j  \right),
\]
where  $\tilde K_{n,j}$ and $\tilde L_{n,j}$ defined  as follows

 \begin{eqnarray}\label{deftildeKL}
 \tilde K_{n,j} &=&\tilde P_{n,M}(i h_j),\\
 \tilde L_{n,j}  &=& \tilde P_{n,M}({ i\tilde h}_j).
 \end{eqnarray}

\end{lemma}

Using   the fact that  $q(s)\in \Vg_A(s)$   in addition, then the error estimates can be bounded from Definition \ref{defthess} as follows:
\begin{corollary} For all $A\geq 1$, there exists $s_{10}(A)\geq 1$ such that for all $s \geq s_{10}(A) $, if $q\in \Vg_A(s)$ and $|\theta'(s)|\leq \frac{C A^{10}}{s^{\frac54}}$, then:\\

a) For all $1\leq n\leq M$, we have 
\[\left|  P_{n,M}\left(-i \left(  \frac{\nu}{2 \sqrt s} +\frac \mu s+\theta'(s) \right)q\right)  \right|\leq C\frac{A^n}{s^{\frac{n+5}{4}}}.\]
b) For $1\leq n\leq M$, we have
\[\left| \tilde P_{n,M}\left(-i \left( \frac{\nu}{2 \sqrt s} +   \frac \mu s+\theta '(s) \right)q\right)  \right|\leq C\frac{A^n}{s^{\frac{n+5}{4}}}.\]

In particular, when $n=0,2,4$,  we can get smaller bounds as follows: 

c) For $n=0$,    we have the following in particular
\begin{eqnarray*}
& & \hspace{-1.5cm} \left|P_{0,M}\left( -i(\frac{\nu}{2 \sqrt{s}}+\frac {\mu }{s}+\theta') q \right) \hspace{-0.1cm}  +  (\frac{\nu}{2\sqrt s}+\frac \mu s+\theta') \{   \delta q_0 + (1 + \delta^2) \tilde q_0 + K_{0,2} q_2 + L_{0,2} \tilde q_2\}  \right|\\
& \leq & C \frac{A^{4}}{s^{\frac{7}{4}}},\\
 & & \hspace{-1.5cm} \left| \tilde P_{0,M}\left( -i(\frac{\nu}{2\sqrt s}+\frac \mu s+\theta')q\right )  + \left(\frac{\nu}{2\sqrt s}+\frac \mu s+\theta' \right) \{   - q_0 - \delta \tilde q_0 + \tilde K_{0,2} q_2 +  \tilde L_{0,2} \tilde q_2\}  \right|\\
& \leq & C\frac{A^{4}}{s^{\frac{7}{4}}},
\end{eqnarray*}
d) For $n =2$, we have 
\begin{eqnarray*}
& & \hspace{-1.5cm} \left|  P_{2,M}\left( -i(\frac{\nu}{2\sqrt s}+\frac \mu s+\theta')q\right )     +  (\frac{\nu}{2\sqrt s}+\frac \mu s+\theta') [ \delta q_2 + (1 + \delta^2) \tilde q_2]  \right|\\
& \leq &C \frac{A^{4}}{s^{\frac{7}{4}}},\\
 & & \hspace{-1.5cm} \left|  \tilde P_{2,M} \left(-i\left (\frac{\nu}{2\sqrt s}+\frac \mu s+\theta'\right )q\right )   + \left (\frac{\nu}{2\sqrt s}+\frac \mu s+\theta'(s)\right )\left (-q_2-\delta\tilde q_2+  \tilde K_{2,4} q_4 +  \tilde L_{2,4} \tilde q_4\right ) \right|\\
& \leq & C \frac{A^{6}}{s^{\frac{9}{4}}},
\end{eqnarray*}
e) For $n =3$, we have 
\begin{eqnarray*}
 \left|   P_{3,M}\left( -i(\frac{\nu}{2\sqrt s}+\frac \mu s+\theta')q\right )   \right| 
& \leq &C \frac{A^{2}}{s^{\frac{3}{2}}},\\
  \left|  \tilde P_{3,M} \left(-i\left (\frac{\nu}{2\sqrt s}+\frac \mu s+\theta'\right )q\right )   \right|& \leq & C \frac{A^{2}}{s^{\frac{3}{2}}},
\end{eqnarray*}
f) For $n =4$, we have 
\begin{eqnarray*}
 \left|   P_{4,M}\left( -i(\frac{\nu}{2\sqrt s}+\frac \mu s+\theta')q\right )   \right| 
& \leq &C \frac{A^{5}}{s^{2}},\\
  \left|  \tilde P_{4,M} \left(-i\left (\frac{\nu}{2\sqrt s}+\frac \mu s+\theta'\right )q\right )   \right|& \leq & C\frac{A^{5}}{s^{2}},
\end{eqnarray*}
\end{corollary}
\textbf{Fourth term: $V_1 q+V_2\bar q$.}

We claim the following
\begin{lemma}[Projection of $V_1 q$ and $V_2 \bar q$ ] \label{lemma-projec-V-1-a+-V-2-bar-q}
(i) It holds that
\beqtn
|V_i(y,s)|\leq C\frac{(1+|y|^2)}{s^{1/2}}, \mbox{ for all $y\in \R$ and $s\geq 1$,}
\label{bdVi1}
\eeqtn
and for all $k\in\N^*$
\beqtn
V_i(y,s)=\dsp \sum_{j=1}^{k}\frac{1}{s^{j/2}}W_{i,j}(y)+\tilde W_{i,k}(y,s),
\label{dcVi}
\eeqtn
where $W_{i,j}$ is an even polynomial of degree $2j$ and $\tilde W_{i,k}(y,s)$ satisfies
\beqtn
\mbox{for all $s\geq 1$ and $|y|\leq s^{1/4}$, } \left| \tilde W_{i,k}(y,s)\right|\leq  C\frac{(1+|y|^{2k+2})}{s^{\frac{k+1}{2}}}.
\eeqtn
(ii) The projection of $V_1 q$ and $V_2 \bar q$ on $(1+i\delta)h_n$ and $i h_n$, and we have
\beqtn
\begin{array}{l}
\dsp|\tilde P_n(V_1 q)|+|\hat P_n(V_1 q)|\\
\leq \dsp\frac{C}{s^{1/2}} \dsp \sum_{j=n-2}^{M}(|\tilde q_j|+|\hat q_j|)
+\sum_{j=0}^{n-3}\frac{C}{s^{\frac{n-j}{4}}}(|\tilde q_j|+|\hat q_j|)
+\frac{C}{s^{1/2}}\left\| \frac{q_-}{1+|y|^{M+1}}\right\|_{L^\infty},

\end{array}
\label{bdVi}
\eeqtn
and the same holds for $V_2\bar q$

\label{lembdVi}
\end{lemma}
\begin{rem} If $n\leq 2$, the first sum in \eqref{bdVi} runs for $j=0$ to $M$ and the second sum doesn't exist.
\end{rem}
If in addition $q(s)\in \Vg_A(s)$, then the error estimates can be bounded from Definition \ref{defthess} as follows:
\begin{corollary} For all $A\geq 1$, there exists $s_{11}(A)\geq 1$ such that for all $s \geq s_{11}(A) $, if $q\in \Vg_A(s)$, then for $3\leq n\leq M$, we have 
\[\dsp  \left|\tilde P_n(V_1 q + V_2 \bar q) \right|+ \left| P_n(V_1 q + V_2 \bar q) \right|\leq \frac{ C A^{n-2}}{s^{\frac{n+1}{4}}}.\]
\end{corollary}

\textit{Proof of Lemma \ref{lembdVi}: }\\
(i)  The estimates of $V_1q$ and $V_2 \bar q$ are the same, so we only deal with $V_1q$. Let $F(u)=\frac{(p+1)}{2}(1+i\delta)\left [|u|^{p-1}-\frac{1}{p-1}\right]$, where $u\in\C$ and consider $z=\frac{y}{s^{1/4}}$. Note that from \eqref{eqq} and \eqref{eqfi0}, we have 
\[V_1(y,s)=F(\varphi(y,s))\mbox{, where } \varphi(y,s)=\varphi_0(\frac{y}{s^{1/4}})+\frac{a}{s^{1/2}}(1+i\delta).\]
Note that there exist positive constant $c_0$ and $s_0$ such that $\varphi_0(z)|$ and $|\varphi(y,s)|=|\varphi_0(\frac{y}{s^{1/4}})+\frac{a}{s^{1/2}}(1+i\delta)|$
are both larger than $\frac{1}{c_0}$ and smaller than $c_0$, uniformly in $|z|<1$ and for $s\geq s_0$. Since $F(u)$ is $C^\infty$ for $\frac{1}{c_0}\leq |u|\leq c_0$, we expand it around $u=\varphi_0(z)$ as follows: for all $s\geq s_0$ and $|z|<1$,
\[
\begin{array}{lll}
\left|F\left(\varphi_0(z)+\frac{a}{s^{1/2}}(1+i\delta)\right)-F\left(\varphi_0(z)\right)\right|&\leq &\dsp\frac{C}{s^{1/2}},\\
\left|F\left(\varphi_0(z)+\frac{a}{s^{1/2}}(1+i\delta)\right)-F\left(\varphi_0(z)\right)-\dsp \sum_{j=1}^{n}\frac{1}{s^{j/2}}F_j(\varphi_0(z))\right|&\leq &\dsp\frac{C}{s^{\frac{n+1}{2}}},
\end{array}
\]
where $F_j(u)$ are $C^\infty$. Hence, we can expand $F(u)$ and $F_j(u)$ around $u=\varphi_0(0)$ and write for all $s\geq s_0$ and $|z|<1$,

\[
\begin{array}{l}
\left|F\left(\varphi_0(z)+\frac{a}{s^{1/2}}(1+i\delta)\right)-F\left(\varphi_0(0)\right)\right|\leq C z^2+\frac{C}{s^{1/2}},\\
\left|F\left(\varphi_0(z)+\frac{a}{s^{1/2}}(1+i\delta)\right)-F\left(\varphi_0(0)\right)-\dsp \sum_{l=1}^{n} c_{0,l}z^{2l}-\sum_{j=1}^{n}\sum_{l=0}^{n-j}\frac{c_{j,l}}{s^{j/2}} z^{2l}\right|\\
\leq C |z|^{2n+2}+\dsp\sum_{j=1}^{n}\frac{C}{j^{1/2}}|z| ^{2(n-j)+2}           \frac{C}{s^{\frac{n+1}{2}}}.
\end{array}
\]
Since $F(\varphi_0(0))=F(\kappa)=0$ and $z=\frac{y}{s^{1/4}}$, this gives us estimates in (i),  when $s\geq s_0$ and $|y|<s^{1/4}$. Since $V_1$ is bounded, the inequalities still valid when $|y|\geq s^{1/4}$ and then when $s\geq 1$.

\medskip

(ii) Note first that it is enough to prove the bound \eqref{decomp2} for the projection of $V_i q$ onto $h_n$ to get the same bound for $\tilde P_{n,M}(V_i q)$ and $P_{n,M}(V_i q)$. Since in addition, the proof for $V_2 \bar q$ is the same as for $V_1 q$, we only prove \eqref{bdVi} for the projection of $V_1 q$ onto $f_n$. Using \eqref{decompq} and \eqref{eqQn}, we see that the projection is given by
\beqtn
\begin{array}{ll}
\dsp \int f_n V_1 q\rho&=\dsp\int f_n V_1 q_-\rho+\sum_{j=0}^{M}\tilde q_j \int f_n\tilde h_j V_1\rho_{\beta}+\sum_{j=0}^{M} q_j \int f_n h_j V_1\rho_{\beta}.\\
\label{bd1}
\end{array}
\eeqtn
The first term can be bounded by
\beqtn
\dsp\int f_n V_1 \left(\frac{1+|y|^2}{s^{1/2}}\right)|q_-||\rho_{\beta}|\leq  \frac{C}{s^{1/2}}\left\| \frac{q_-}{1+|y|^{M+1}}\right\|_{L^\infty}.
\label{bd2}
\eeqtn
Now we deal with the second term. We only focus on the terms involving $h_j$ .\\
If $j\geq n-2$, we use \eqref{bdVi1} to write $|\int f_n h_j V_1\rho_{\beta}|\leq \frac{C}{s^{1/2}}$.\\
If $j\leq n-3$, then we claim that
\beqtn
\dsp\left |\int f_n h_j V_1\rho_{\beta}\right|\leq \frac{C}{s^{\frac{n-j}{4}}},
\label{bd3}
\eeqtn
(this actually vanishes if $j$ and $n$ have different parities). It is clear that \eqref{bdVi} follows from \eqref{bd1}, \eqref{bd2} and \eqref{bd3}.\\
Let us prove \eqref{bd3}. Note that $k\equiv\left [ \frac{n-j-1}{2}\right ]$ (which is in $\N^*$ since $j\leq n-3$) is the largest integer such that $j+2k <n$. We use \eqref{dcVi} to write
\beqtn
\begin{array}{l}
 \dsp \int f_nh_j V_1 \rho_{\beta}=\dsp\int_{|y|< s^{1/4}}f_nh_j V_1 \rho_{\beta}+\int_{|y|>s^{1/4}}f_nh_j V_1 \rho_{\beta},\\
= \dsp\sum_{l=1}^{k}\frac{1}{s^{l/2}}\int_{|y|< s^{1/4}}f_n h_jW_{1,l}\rho_{\beta} +O\left(\frac{1}{s^{\frac{ \left [ \frac{n-j-1}{2}\right ]+1  }{2}}}\int(1+|y|^{n-j+1})|f_n||h_j|\rho_{\beta} dy\right) \\
+ \int_{|y|>s^{1/4}}h_nh_j V_1 \rho,\\
= \dsp\sum_{l=1}^{k}\frac{1}{s^{l/2}}\int_{\R^N}h_n h_jW_{1,l}\rho +O\left(\frac{1}{s^{\frac{ \left [ \frac{n-j-1}{2}\right ]+1  }{2}}}\right)
-\sum_{l=1}^{k}\frac{1}{s^l}\int_{|y|>s^{1/4}} f_n h_j W_{1,l}\rho_{\beta} \\
+ \int_{|y|>s^{1/4}}f_nh_j V_1 \rho_{\beta},
\end{array}
\eeqtn
since $deg(h_j W_{1,l})=j+2l\leq j+2k<n=deg(h_n)$, $h_n$ is orthogonal to $h_j W_{1,l}$ and
\[\int_{\R^N}f_n h_j W_{1,l}\rho_{\beta}=0.\]
Since $|\rho_{\beta}(y)|\leq \dsp C e^{-cs^{1/2}}$ when $|y|>s^{1/4}$, the integrals over the domain ${|y|>s^{1/4}}$ can be bounded by
\[\dsp Ce^{-cs^{1/2}} \int_{|y|>s^{1/4}} |f_n||h_j|(1+|y|^{2k})\sqrt{\rho_{\beta}} \leq Ce^{-cs^{1/2}}. \]
Using that $\left [ \frac{n-j-1}{2}\right ]+1\geq \frac{n-j}{2}$, we deduce that \eqref{bd3} holds. Hence, we have proved \eqref{bdVi} and this concludes the proof of Lemma \ref{lembdVi}. $\blacksquare$.\\

As a matter of fact, we need more refinements in the cases  where $n=0$ and $2$ for the terms $\tilde P_{2,M}(V_1 q)$, $\tilde P_{2,M}(V_2\bar q)$, $P_{0,M}(V_1 q)$ and $P_{0,M}(V_2 \bar q)$. More precisely

\begin{lemma}[Projection of $V_1 q$ and $V_2\bar q$ on $\hat h_0, \tilde h_0, \hat h_2, \tilde h_2, \hat  h_4$ and $ h_4$ ]\label{lemma-projec-potential}
Using the definition of $V_1, V_2$, the following hold:

(i) It holds that for $i=1,2$
\beqtn
\forall s\geq 1 \mbox{ and } |y|< s^{1/4},\;\; \left| V_i(y,s)-\frac{1}{s^{1/2}} W_{i,1}(y)-\frac{1}{s} W_{i,2}(y) \right|\leq \frac{C}{s^{3/2}}(1+|y|^6),
\label{estiVii}
\eeqtn
where
\beqtn
\begin{array}{lll}

W_{1,1} (y) &=&  -\frac{(p+1) b}{2(p-1)^2} (1 + i \delta) (y^2 -2(1 -\delta \beta) ),\\[0.2cm]

W_{1,2}  (y)&=&(1+i\delta)\frac{b^2(p+1)}{2(p-1)^4}\left\{  (p-1)y^4- [ 2(1-\beta\delta)(p-2+\delta^2)]y^2 \right.\\
& + & \left. 2(p-2+\delta^2)(1-\beta\delta)^2 \right\},\\[0.2cm]

W_{2,1} (y)&=&-(1+i\delta)\frac{b}{2(p-1)^2}\left(p-1+2i\delta\right)(y^2-2(1-\beta\delta)),\\

W_{2,2}(y) &=&(1+i\delta)\frac{b^2}{2(p-1)^4}\left\{ (p-2+2i\delta)(p-1 +i\delta) y^4\right.\\
& - &  [ 2(p-1)(p-2) +(2p -10)\delta^2   + (8p -16 ) i \delta ](1- \delta \beta) y^2\\
& +& \left.  (1 - \delta \beta) ( 2 p^2+8 i p \delta+4 - 16 i \delta-10 \delta^2-6p+2p\delta^2)  \right\}.  
\end{array}
\label{defWii}
\eeqtn

(ii) The projection of $V_1 q$ and $V_2 \bar q$ on $\tilde h_2$ satisfy

\begin{eqnarray}
& & \left| \tilde P_n( V_1 q + V_2 \bar q ) - \frac{1}{\sqrt{s}} \sum_{j\geq 0} [ \tilde C_{n,j} q_j +\tilde D_{n,j} \tilde q_j]  - \frac{1}{s} \sum_{j \geq 0} [ \tilde E_{n,j} q_j + \tilde F_{n,j} \tilde q_j   ]  \right| \nonumber \\
& \leq & \frac{C}{s^{\frac{3}{2}}} \sum_{j \geq 0} [ |\hat q_j| + |\tilde q_j| ]   + \frac{1}{\sqrt s} \left\|  \frac{q_- (.,s)}{1 + |y|^M} \right\|_{L^\infty},\label{bdVtilde2}
\end{eqnarray}
and 
\begin{eqnarray}
& & \left| P_n( V_1 q + V_2 \bar q ) - \frac{1}{\sqrt{s}} \sum_{j\geq 0} [  C_{n,j} q_j + D_{2,j} \tilde q_j]  - \frac{1}{s} \sum_{j \geq 0} [  E_{n,j} q_j +  F_{n,j} \tilde  q_j   ]  \right| \nonumber \\
& \leq & \frac{C}{s^{\frac{3}{2}}} \sum_{j \geq 0} [ |\hat q_j| + |\tilde q_j| ]   + \frac{1}{\sqrt s} \left\|  \frac{q_- (.,s)}{1 + |y|^M} \right\|_{L^\infty}.\label{projec-P-n-V-q}
\end{eqnarray}

where for all $n,  j \geq 0$, we have

 \begin{eqnarray}
 C_{n,j} =P_{n,M}(W_{1,1}h_j+W_{2,1}\bar h_j)& \tilde C_{n,j}=\tilde P_{n,M}(W_{1,1} h_j+W_{2,1}\bar h_j),\\
 D_{n,j} = P_{n,M}(W_{1,1} \tilde h_j+W_{2,1}\bar{\tilde h}_j)& \tilde D_{n,j}=\tilde P_{n,M}(W_{1,1} \tilde h_j+W_{2,1}\bar{\tilde h}_j),\\
 E_{n,j} =P_{i,M}(W_{1,2}h_j+W_{2,2}\bar h_j)& \tilde E_{n,j}=\tilde P_{i,M}(W_{1,2} h_j+W_{2,2}\bar h_j),\\
 F_{n,j} = P_{n,M}(W_{1,2} \tilde h_j+W_{2,2}\bar{\tilde h}_j)& \tilde F_{n,j}=\tilde P_{n,M}(W_{1,2} \tilde h_j+W_{2,2}\bar{\tilde h}_j).
\label{Projection-of-potentials}
 \end{eqnarray}

\end{lemma}
In addition, by using that fact that  $q(s)\in \Vg_A(s)$, then the error estimates can be bounded from Definition \ref{defthess} as follows;
\begin{corollary}\label{corprojPtls} For all $A\geq 1$, there exists $s_{12}(A)\geq 1$ such that for all $s\geq s_{12}(A)$, if $q(s)\in \Vg_A(s)$, then 
\begin{eqnarray*}
& & \left| P_{0,M} (V_1 q + V_2 \bar q)  - \left(C_{0,0} \frac{q_0}{\sqrt{s}} + D_{0,0} \frac{\tilde q_0}{\sqrt{s}} + C_{0,2} \frac{q_2}{\sqrt{s}}+ D_{0,2} \frac{\tilde q_2}{\sqrt{s}} \right) \right| \leq  C\frac{A^4}{s^{\frac{7}{4}}},
\\
& & \left| \tilde P_{0,M} (V_1 q + V_2 \bar q)  - \left(   \tilde D_{0,0} \frac{\tilde q_0}{\sqrt{s}} + \tilde C_{0,2} \frac{q_2}{\sqrt{s}} + \tilde D_{0,2} \frac{\tilde q_2}{\sqrt{s}}  \right) \right|   \leq C \frac{A^4}{s^\frac{7}{4}},\\
& & \left|   P_{2,M}(V_1 q+V_2 \bar q) - \left( \frac{D_{2,0} \tilde q_0}{\sqrt s}  + \frac{C_{2,2} q_2 }{ \sqrt s} + \frac{D_{2,2} \tilde q_2 }{ \sqrt s}  \right)  \right|  \leq C \frac{A^4}{s^\frac{7}{4}},\\
& & \left|  \tilde P_{2,M} \left(  V_1 q + V_2 \bar q  \right)  -  \frac{1}{\sqrt{s} } \left\{     \tilde q_0 \tilde D_{2,0}   + q_2 \tilde C_{2,2}+ \tilde q_2  \tilde D_{2,2}  +  q_4 \tilde C_{2,4} + \tilde q_4 \tilde D_{2,4} \right\}  \right. \\
 &-& \left.  \frac{1}{s} \left\{    \tilde q_0 \tilde F_{2,0} + q_2 \tilde E_{2,2} + \tilde q_2 \tilde F_{2,2}        \right\}  \right|   \leq  C\frac{A^6}{s^\frac{9}{4}},\\
 & & \left|  P_{4,M} ( V_1 q + V_2 \bar q) - \left( C_{4,2}  \frac{q_2}{\sqrt{s}} + D_{4,2} \frac{\tilde q_2}{\sqrt{s}} \right) \right| \leq C \frac{A^4}{s^{\frac{7}{4}}},\\
 & & \left| \tilde  P_{4,M} ( V_1 q + V_2 \bar q) - \left( \tilde C_{4,2}  \frac{q_2}{\sqrt{s}} + \tilde D_{4,2} \frac{\tilde q_2}{\sqrt{s}} \right) \right| \leq C \frac{A^4}{s^{\frac{7}{4}}}.
\end{eqnarray*}
and
\begin{eqnarray*}
\left|  P_{3,M}(V_1 q + V_2 \bar q)\right| \leq  \frac{CA^4}{s^2},\\
\left|  \tilde P_{3,M}(V_1 q + V_2 \bar q)\right| \leq  \frac{CA^4}{s^2}.
\end{eqnarray*}
\end{corollary}

\textit{Proof of Lemma \ref{lemma-projec-potential}: }(i) This is a simple, but lengthy computation that we omit. For more details see Appendix \ref{appendix-potential}. In addition to that, item $(ii)$ directly follows Lemma \ref{lemma-projec-V-1-a+-V-2-bar-q},  and item $(i)$ and Definitions of the projection $P_{n,M} $ and $\tilde P_{n,M}$, defined as in \eqref{decomp3}. Finally, in order to know the exact formulas of the constants in item $(ii)$, we kindly address the readers to Appendix  \ref{appendix-potential}
$\blacksquare$

\medskip

\textbf{Fifth term: $B(q,y,s)$}
Let us recall from \eqref{eqqd} that:
\[B(q,y,s)=(1+i\delta)\left ( |\varphi+q|^{p-1}(\varphi+q)-|\varphi|^{p-1}\varphi-|\varphi|^{p-1} q-\frac{p-1}{2}|\varphi|^{p-3}\varphi(\varphi \bar q+\bar\varphi q)\right).\]
We have the following
\begin{lemma}\label{lemfifB} The function $B=B(q,y,s)$ can be decomposed for all $s\geq 1$ and $|q|\leq 1$ as
\begin{eqnarray*}
\dsp \sup_{|y|<s^{1/4}}
\left|B-\sum_{l=0}^{M}  \sum_{\begin{array}{l}  0\leq j,k\leq M+1\\2\leq j+k\leq M+1   \end{array}  }\frac{1}{s^{l/2}} \left [   B_{j,k}^{l}(\frac{y}{s^{1/4}}) q^j\bar q^k+\tilde B_{j,k}^{l}(y,s) q^j \bar q^k \right ]\right| & & \\
& & \hspace{-1.6cm}\leq  C|q|^{M+2}+\frac{C}{s^{\frac{M+1}{2}}},
\end{eqnarray*}
where $ B_{j,k}^{l}(\frac{y}{s^{1/4}})$ is an even polynomial of degree less or equal to $M$ and the rest $\tilde  B_{j,k}^{l}(y,s)$ satisfies 
\[\forall s\geq 1 \mbox{ and } |y|< s^{1/4},\; \left|   \tilde  B_{j,k}^{l}(y,s)\right |\leq C\frac{1+|y|^{M+1}}{s^{\frac{M+1}{2}}}.\]
Moreover,
\[\forall s\geq 1  \mbox{ and } |y|< s^{1/4},\; \left|   B_{j,k}^{l}(\frac{y}{s^{1/4}})+\tilde  B_{j,k}^{l}(y,s)\right |\leq C.\]
On the other hand, in the region $|y|\geq s^{1/4}$, we have
\beqtn
|B(q,y,s)|\leq C|q|^{\bar p},
\eeqtn
for some constant $C$ where $\bar p=\min (p,2)$.
\end{lemma}
\textit{Proof:} See the proof of Lemma 5.9, page 1646 in \cite{MZ07}. $\blacksquare$

\begin{lemma}[The quadratic term $B(q,y,s)$]\label{lembin} For all $A\geq 1$, there exists $s_{13} \ge 1$ such that for all $s\geq s_{13}$, if $q(s)\in \Vg_A(s)$, then:\\
a) the projection of $B(q,y,s)$ on $ h_n$ and on $ \tilde h_n$, for  $n \geq 3$ satisfy
\beqtn
\left|  \tilde P_{n,M} (B(q,y,s)) \right| +  \left| P_{n,M} (B(q,y,s)) \right| \leq C\frac{A^n}{s^{\frac{n+2}{4}}}.
\label{lembin1}
\eeqtn
b) For $n=0,1,3,4$, we have
\beqtn
\left| \tilde P_{n,M} (B(q,y,s)) \right| + \left| P_{n,M} (B(q,y,s)) \right| \leq \frac{C}{s^2},
\label{lembin2}
\eeqtn
c) For  $ n = 2$, we have
\beqtn
\left|  P_{2,M} (B(q,y,s)) \right|\leq \frac{C}{s^2}.
\label{lembin3}
\eeqtn
\begin{equation}\label{estima-B-h-2}
  \left|  \tilde P_{2,M} (B(q)) - \left(  \tilde B_2 \tilde q_2^2  +  B_1  \frac{\tilde q_2}{s} +\frac{B_2}{s^2} \right) \right| \leq   \frac{C A^6}{s^{\frac{9}{4}}} ,    
\end{equation}
where
 \begin{eqnarray*}
\tilde B_2   & = &   \frac{1}{\kappa} \left\{    4(p - \delta^2) - \delta \beta (6+ 4 p+ 2\delta^2)    \right\}, \\
B_1 & =& \frac{1}{\kappa} \left\{    (-7\delta^2 \beta + p \beta -6 \beta) R_{2,1}^* - (  p- \delta^2)  \tilde R_{0,1}^*         \right\},\\
B_2 & = & \frac{(R_{2,1}^*)^2}{s^2} \frac{(32-64 \delta \beta)}{8 \kappa}.
\end{eqnarray*}
\end{lemma}
\textit{Proof :}  As a matter of fact, we only prove estimate \eqref{lembin1}. In addition to that, the proofs of \eqref{lembin2}, \eqref{lembin3} and \eqref{estima-B-h-2}  will be  given in  the same way that  using the Taylor expansion of $B$, established in Appendix \ref{Appendix-B-q}.

\medskip

Let us start the  proof of  estimate \eqref{lembin1} for the projection on $h_n$,  since it implies the same estimate on $\tilde P_n$ and $ P_n$ through \eqref{decomp3}. We have
\[ \dsp\int h_n B(q,y,s)\rho dy =\int_{|y|<s^{1/4}}h_n B(q,y,s)\rho dy +\int_{|y|>s^{1/4}}h_n B(q,y,s)\rho dy.\]
Using Lemma \ref{lemfifB}, we deduce that
\begin{eqnarray*}
&& \left|  \dsp\int_{|y|<s^{1/4}}h_n B(q,y,s)\rho dy \right.\\
& - & \left.   \int_{|y|<s^{1/4}}h_n \rho\sum_{l=0}^{M}  \sum_{\begin{array}{l}  0\leq j,k\leq M+1\\2\leq j+k\leq M+1   \end{array}  }\frac{1}{s^{l/2}} \left [   B_{j,k}^{l}(\frac{y}{s^{1/4}}) q^j\bar q^k+\tilde B_{j,k}^{l}(y,s) q^j \bar q^k \right ]   \right|\\
& \leq & C\dsp\int_{|y|<s^{1/4}}|h_n|| \rho |(|q|^{M+2}+\frac{1}{s^{\frac{M+1}{2}}}).
\end{eqnarray*}
Let us write
\[\dsp B_{j,k}^{l}(\frac{y}{s^{1/4}}) =\sum_{i=0}^{M/2} b_{j,k}^{l,i}(\frac{y}{s^{1/4}}) \left(\frac{y}{s^{1/4}}\right)^{2i},\]

\[\dsp q^j=\left( \sum_{m=0}^{M}\tilde q_m \tilde h_m+ \hat q_m \hat h_m+q_-\right)^j,\;\;q^k=\left( \sum_{m=0}^{M}\tilde q_m \tilde h_m+ \hat q_m \hat h_m+q_-\right)^k,\]
where $b_{j,k}^{l,i}$ are the coefficients of the polynomials $B_{j,k}^{l}$. Using the fact that $\|q(s)\|_{L^\infty}\leq 1$ (which holds for $s$ large enough, from the fact that $q(s)\in \Vg_A(s)$ and Definition \ref{defthess}). We deduce that
\[|q^j-q_{+}^{j}|\leq C(|q_-|^j+|q_-|).\]
Using that $q(s)\in \Vg_A(s)$ and the fact that $\sqrt s\geq 2 A^2$, we deduce that in the region $|y|\leq s^{1/4}$, we have $ |q_-|\leq \frac{1}{s^{1/4}}(\frac{A}{s^{1/4}})^{M+1}(1+|y|)^{M+1}$ and that
\[\dsp |q^j-(\sum_{m=0}^{M} \tilde q_m \tilde h_m+\hat q_m \hat h_m)^j| \leq C  \left(\frac{A}{s^{1/4}}\right)^{M+1} \frac{1}{s^{1/4}}(1+|y|)^{jM+j}.\]
In the same way, we have
 \[\dsp |q^k-(\sum_{m=0}^{M} \tilde q_m \tilde h_m+\hat q_m \hat h_m)^k| \leq C  \left(\frac{A}{s^{1/4}}\right)^{M+1} \frac{1}{s^{1/4}}(1+|y|)^{kM+k},\]
hence, the contribution coming from $q_-$ is controlled by the right-hand side of  \eqref{lembin1}. Moreover for all $j,k$ and $l$, we have
\beqtn
\dsp\left | \int_{|y|< s^{1/4}} h_n\rho  B_{j,k}^{l}(\frac{y}{s^{1/4}}) q_{+}^{j}\bar q^{k}_{+}-\int h_n\rho  B_{j,k}^{l}(\frac{y}{s^{1/4}}) q_{+}^{j}\bar q^{k}_{+} \right |\leq C e^{-C\sqrt s}.
\label{lembin4}
\eeqtn
To compute the second term on the left had side of \eqref{lembin4}, we notice that $B_{j,k}^{l}(\frac{y}{s^{1/4}}) q_{+}^{j}\bar q^{k}_{+}$  is a polynomial in $y$ and that the coefficient of the term of degree $n$ is controlled by the right had side of \eqref{lembin1} since $q\in \Vg_A$.\\
Moreover, using that $\sqrt s\geq 2 A^2$, we infer that $|q|\leq \frac{1}{s^{1/4}}(1+|y|)^{M+1}$ in the region $|y|\leq s^{1/4}$
 and hence for all $j,k$ and $l$, we have
 \[\dsp \left| \int_{|y|< s^{1/4}}h_n \rho \frac{1}{s^{l/2}} \tilde B_{j,k}^{l}(y,s) q^j \bar q^k \right| \leq C\frac{1}{s^{\frac l 2+\frac{M+1+j+k}{4}}} \]
and
\[\dsp \left|\int _{|y|< s^{1/4}} h_n \rho(|q|^{M+2}+\frac{1}{s^{\frac{M+1}{2}}}) \right|\leq C \frac{1}{s^{\frac{M+2}{4}}}.\]
The terms appearing in these two inequalities are controlled by the right hand side of \eqref{lembin1}.\\
Using the fact that $\|q(s)\|_{L^\infty}\leq 1$ and \eqref{eqqd}, we remark that $|B(q,y,s)|\leq C$. Since $|\rho(y)|\leq C e^{-s\sqrt s} $ for all $|y|>s^{1/4}$, it holds that
\[\left | \int_{|y|>s^{1/4}}h_n B(q,y,s)\rho dy \right|\leq Ce^{-C\sqrt s}.\]
This concludes the proof of Lemma \ref{lembin}. $\blacksquare$

\medskip

\textbf{Sixth term: $R^*(\theta',y,s)$}\\
In the following, we expand $R^*$ as a power series of $\frac 1 s$ as $s\to\infty$, uniformly for $|y|\leq s^{1/4}$.
\begin{lemma}[Power series of $R^*$ as $s\to\infty$] For all $n\in\N$,
\beqtn
R^*(\theta',y,s)=\Pi_n(\theta',y,s)+\tilde\Pi_n(\theta',y,s),
\eeqtn
where,
\beqtn
\Pi_n(\theta',y,s)=\sum_{k=0}^{n-1}\frac{1}{s^{\frac{k+1}{2}}} P_k(y)-i\theta'(s)\left(\frac{a}{s^{1/2}}(1+i\delta)+\sum_{k=0}^{n-1}e_k\frac{y^{2k}}{s^{k/2}}\right),
\eeqtn
and
\beqtn
\forall |y|<s^{1/4},\;\;\left |\tilde\Pi_n(\theta',y,s)\right|\leq C(1+s|\theta'(s)|)\frac{(1+|y|^{2n})}{s^{\frac{n+1}{2}}},
\eeqtn
where $P_k$ is a polynomial of order $2k$ for all $k\geq 1$ and $e_k\in \R$.\\
In particular,
\beqtn
\begin{array}{l}
\dsp \sup_{|y|\leq s^{1/4}}\left | R^*(\theta',y,s)-\sum_{k=0}^{1} \frac{1}{s^{\frac{k+1}{2}}}P_k(y)+i\theta'\left[\kappa+\frac{(1+i\delta)}{s^{1/2}}\left( a-\frac{b\kappa y^2}{(p-1)^2}\right) \right]\right |\\
\leq C\left (\frac{1+|y|^4}{s^{3/2}}+C|\theta '|\frac{y^4}{s^2}\right).
\end{array}
\eeqtn
\label{decompR}
\end{lemma}
\textit{Proof: }
Using the definition of $\varphi$ \eqref{defi-varphi}, the fact that $\varphi_0$ satisfies \eqref{eqfi0} and \eqref{eqq}, we see that $R^*$ is in fact a function of $\theta'$, $z=\frac{y}{s^{1/4}}$ and $s$ that can be written as 
\beqtn
\begin{array}{l}
\displaystyle R^*(\theta',y,s)=\frac{1}{4}\frac{z}{s}\nabla_z \varphi_0(z)+\frac{1}{2}\frac{a}{s^{3/2}}(1+i\delta)+ \frac{1}{s^{1/2}}\Delta_z \varphi_0(z)\\
\displaystyle -\frac{a(1+i\delta)^2}{(p-1)s^{1/2}}+\left( F\left(\varphi_0(z)+\frac{a}{s^{1/2}}(1+i\delta)\right)-F(\varphi_0)\right)\\
\displaystyle -i\frac{\mu}{s}\left( \varphi_0(z)+\frac{a}{s^{1/2}} (1+i\delta)\right)-i\theta'(s)\left( \varphi_0(z)+\frac{a}{s^{1/2}} (1+i\delta)\right),
\end{array}
\label{eqreste}
\eeqtn

\noindent
where  $F(u)=(1+i\delta)|u|^{p-1}u$.

\medskip
Since $|z|<1$, there exists positive $c_0$ and $s_0$ such that $|\varphi_0(z)|$ and $|\varphi_0(z)+\frac{a}{s^{1/2}}(1+i\delta)|$ are both larger that $\frac{1}{c_0}$ and smaller than $c_0$, uniformly in $|z|<1$ and $s>s_0$. Since $F(u)$ is $C^\infty$ for $\frac{1}{c_0}\leq |u|\leq c_0$, we expand it around $u=\varphi_0(z)$ as follows
\[\left | F\left(\varphi_0(z)+\frac{a}{s^{1/2}}(1+i\delta)\right)-F\left(\varphi_0(z)\right)-  \sum_{j=1}^{n}\frac{1}{s^{j/2}}F_j\left(\varphi_0(z)\right)\right |\leq C\frac{1}{s^{\frac{n+1}{2}}},\]
where $F_j(u)$ are $C^{\infty}$. Hence, we can expand $F_j(u)$ around $u=\varphi_0(0)$ and write 
\[\left | F\left(\varphi_0(z)+\frac{a}{s^{1/2}}(1+i\delta)\right)-F\left(\varphi_0(z)\right)-  \sum_{j=1}^{n}\sum_{l=0}^{n-j}  \frac{c_{j,l}}{s^{j/2}} z^{2l}\right |\leq\sum_{j=1}^{n}\frac{C}{s^{\frac{j}{2}}}|z|^{2(n-j)+2}+ \frac{C}{s^{\frac{n+1}{2}}},\]
Similarly, we have the following
\[\left|\frac{z}{s}\nabla_z\varphi_0(z)-\frac{|z|^2}{s}\sum_{j=0}^{n-2}d_j z^{2j} \right |\leq \frac{C}{s}|z|^{2n}, \]
\[\left |\frac{1}{s^{1/2}}\Delta_z\varphi_0(z)-\frac{1}{s^{1/2}}\sum_{j=0}^{n-1}b_j z^{2j}\right |\leq \frac{C}{s^{1/2}}|z|^{2n}\mbox{  and  }\left|\varphi_0(z)-\sum_{j=0}^{n-1}e_j z^{2j}\right|\leq C|z|^{2n}.\]
Recalling that $z=\frac{y}{s^{1/4}}$, we get the conclusion of the Lemma. $\blacksquare$

\medskip

In the following, we introduce $F_j(R^*)(\theta,s)$ as the projection of the rest term $R^*(\theta',y,s)$ on the standard Hermite polynomial, introduced in Lemma \ref{lemma-Jordan-block's}.
\begin{lemma}[Projection of $R^*$ on the eigenfunction of $\Lg$]
It holds that $F_j(R^*)(\theta',s)\equiv 0$ when $j$ is odd, and $|F_j(R^*)(\theta',s)|\leq C\frac{1+s|\theta'(s)|}{s^{\frac j 4+\frac 12}} $, when $j$ is even and $j\geq 4$.
\end{lemma}
\textit{Proof: } Since $R^*$ is even in the $y$ variables and $f_j$ is odd when $j$ is odd, $F_j(R^*)(\theta',s)\equiv 0$, when $j$ is odd.\\
Now, when $j$ is even, we apply Lemma \ref{decompR} with $n=[\frac j2]$ and write
\[R^*(\theta',y,s)=\Pi _{\frac j2}(\theta',y,s)+O\left(\frac{1+s|\theta'(s)|+|y|^j}{s{\frac j4+\frac 12}}\right),\]
where $\Pi_{\frac j2 }$ is a polynomials in $y$ of degree less than $j-1$. Using the definition of $F_j(R^*)$ (projection on the $h_j$ of $R^*$), we write
\beqtn
\begin{array}{l}
\dsp\int_{\R^N}R^* h_j \rho=\int_{|y|<s^{1/4}} R^* h_j \rho dy+\int_{|y|>s^{1/4}} R^* h_j \rho dy\\

=\displaystyle\int_{|y|< s^{1/4}} \Pi _{\frac j2} h_j \rho dy+O\left( \int_{|y|<s^{1/4}}  \frac{1+s|\theta'(s)|+|y|^j}{s{\frac j4+\frac 12}}h_j \rho dy\right)+\int_{|y|>s^{1/4}} R^* h_j \rho dy\\

\dsp=\displaystyle\int_{\R^N} \Pi _{\frac j2} h_j \rho dy+O\left(  \frac{1+s|\theta'(s)|}{s{\frac j4+\frac 12}}\right)+\int_{|y|>s^{1/4}} R^* h_j \rho dy+\int_{|y|>s^{1/4}} \Pi_{\frac j2} h_j \rho dy.
\end{array}
\eeqtn
We can see that $\int_{\R^N} \Pi _{\frac j2} h_j \rho dy=0$ because $h_j$ is orthogonal to all polynomials of degree less than $j-1$. Then, note that both integrals over the domain $\{|y|>s^{1/4}\}$ are controlled by
\[\dsp\int _{s^{1/4}}\left(|R^*(\theta',y,s)|+1+|y|^j\right)(1+|y|^j)\rho dy.\]
Using the fact that $R(y,s)$ measures the defect of $\varphi(y,s)$ from being an exact solution of \eqref{equa-w}. However, since $\varphi$ is an approximate solution of \eqref{equa-w}, one easily derive the fact that 

\beqtn
\begin{array}{l}
\|R(s)\|_{L^\infty}\leq \frac{ C}{\sqrt{s}}\mbox{, and}\\
|R^*(\theta', y,s)|\leq \frac{C}{\sqrt{s}} +|\theta'(s)|. 
\end{array}
\eeqtn 
Using the fact that $|\rho(y)|\leq C e^{\sqrt{-s}}$, for $|y|>s^{1/4}$, we can bound our integral by
\[\dsp C(1+|\theta'(s)|)\int_{\R^N}(1+|y|^j)^2 e^{c\sqrt{-s}}\sqrt{\rho}dy=C(j)(1+|\theta'(s)|)e^{c\sqrt{-s}}.\]
This inequality gives us the result for $j\geq 4$. $\blacksquare$

\begin{lemma}[Projection of $R^*$ on the eigenfunction $\tilde h$ and $h_n$]\label{lemma-prejection-R-*constant}
Let us consider $R^*$ defined as in the above, then the following hold:

\begin{itemize}
\item[(i)] For  $j \geq 4$ which is even, then $\tilde P_j(R^*)(\theta',s)$ and $ {P}_j(R^*)(\theta',s)$ are $O\left(\frac{1+s|\theta'|}{s^{\frac j4+\frac 12}}\right)$.
\item[$(ii)$]  For all $j$ odd, we have   $\tilde P_j(R^*)(\theta',s)= {P}_j(R^*)(\theta',s)=0$.
\item[$(iii)$] For $j=0$, we have 
\begin{eqnarray*}
P_{0,M} (R^*(\theta'(s),s)) & =& \frac{R_{0,0}^*}{s^\frac{1}{2}} + \frac{R^*_{0,1}}{s} + \frac{R_{0,2}^*}{s^\frac{3}{2}}  \\
&+& \theta'(s) \left(  - \kappa +  \frac{\Theta^*_{0,0}}{\sqrt s}   + O\left(\frac{1}{s} \right) \right) + O \left( \frac{1}{s^2} \right) ,\\
\tilde P_{0,M} (R^*(\theta'(s),s)) & =& \frac{\tilde R_{0,0}^*}{s^\frac{1}{2}} + \frac{\tilde R^*_{0,1}}{s} + \frac{\tilde R_{0,2}^*}{s^\frac{3}{2}}  \\
&+& \theta'(s) \left(  \frac{\tilde \Theta^*_{0,0}}{\sqrt s}  + O\left(\frac{1}{s} \right) \right)  +O \left( \frac{1}{s^2} \right).
 \end{eqnarray*}
 \item[$(iv)$] For $j =2$, we have 
\begin{eqnarray*}
P_{2,M} (R^*(\theta'(s),s)) & =&  \frac{R^*_{2,1}}{s} + \frac{R_{2,2}^*}{s^\frac{3}{2}}  +  \theta'(s) \left(  \frac{ \Theta^*_{2,0}}{\sqrt s}  + O\left(\frac{1}{s} \right) \right)   +O \left( \frac{1}{s^2} \right),\\
\tilde P_{2,M} (R^*(\theta'(s),s)) & =&  \frac{\tilde R^*_{2,1}}{s} + \frac{\tilde R_{2,2}^*}{s^\frac{3}{2}}  + \frac{\tilde R_{2,3}^*}{s^2} \\
&+ & \theta'(s) \left(  \frac{ \tilde \Theta^*_{2,0}}{\sqrt s}  + \frac{\tilde \Theta_{2,1}^*}{s} +  O\left(\frac{1}{s^\frac{3}{2}} \right) \right)   +O \left( \frac{1}{s^\frac{5}{2}} \right),
\end{eqnarray*} 
 where $R_{j,k}^*, \tilde R_{j,k}^*, \Theta^{*}_{j,k}, \tilde \Theta^*_{j,k}$ are constants, depending on  $p, \delta, \beta $ only. For more details see page \pageref{page-constant-R^*-j-k} and equation \eqref{Def-Theta-ij}.
 
 and the other constants are complicated, will we give them in the proof of this Lemma.
\end{itemize}
\end{lemma}

\textit{Proof: }
See Appendix   \ref{appendix-expansion-R}. $\square $

\medskip

\textbf{Part 2: Proof of Proposition \ref{propode}}
\label{delicatq}
\medskip

In this part, we consider  $A\geq 1$ and take $s$ large enough so that Part 1   is satisfied.\\

+  The proof of  item $(i)$:  We control $\theta '(s)$, from the projection of \eqref{eqq} on $h_0(y)=i$, we obtain

\beqtn
  q_0'=
 c_2 \tilde q_2 -P_{0,M}\left (\left(\frac{\nu}{2\sqrt s}+\frac \mu s+\theta' \right)q \right )+ P_{0,M}(V_1 q+V_2 \bar q)+P_{0,M}(B)+ P_{0,M}(R^*(\theta'(s),s)), 
 \label{Projtild0}
 \eeqtn
 where $c_2=2\beta(1+\delta^2)$, defined as in Lemma  \ref{lemma-Jordan-block's}. In addition to that,  from the fact that $q_0\equiv 0$ by the modulation, we also obtain that 
 $$q_0'\equiv 0.$$
 By the definition of the shrinking set $\Vg_A(s)$ in  Definition \ref{defthess}, Corollary \ref{corprojPtls}, Lemma \ref{lembin} and Lemma \ref{lemma-prejection-R-*constant}, we obtain the following: 
\begin{eqnarray}
& & \kappa \theta' (s) = c_2 \tilde q_2  + \frac{\tilde q_2}{\sqrt{s}} \left[ D_{0,2} - \frac{\nu}{2} L_{0,2} +\frac{\Theta^*_{0,0} c_2}{\kappa} \right] + \frac{R^*_{0,1}}{s} \label{asymptotic-theta's-in-pro-theta-0}\\
& +&  \frac{1}{s^{\frac{3}{2}}} \left[ \frac{\nu}{2} (1 + \delta^2) \tilde R^*_{0,1} -\frac{\nu}{2} K_{0,2} R^*_{2,1} - D_{0,0} \tilde R^*_{0,1} + C_{0,2} R^*_{2,1}+R^*_{0,2} + \frac{\Theta^*_{0,0} R^*_{0,1}}{\kappa}\right] + O\left( \frac{1}{s^{\frac{7}{4}}}\right)\nonumber.
\end{eqnarray}
In addition to that, we derive again from the shrinking set that
$$ \left| c_2  \tilde q_2 (s)   + \frac{R^*_{0,1}}{s}  \right| \leq  \frac{ CA^{10}}{s^{\frac{5}{4}} }.$$
Thus, we obtain
$$ \left|  \kappa \theta'(s)          \right|  \leq \frac{CA^{10}}{s^\frac{5}{4}}. $$
and item $ (i)$  of Proposition \ref{propode}. 

\bigskip


   
 

+ The proof of item $(iii)$: 
Let us  project  equation \eqref{eqq} on $\tilde h_2$, we get
\begin{eqnarray}
\tilde q_2' &=&  \tilde P_{2, M} (\mathcal{L}_{\beta, \delta} q )  + \tilde P_{2,M} \left(-i\left (\frac{\nu}{2\sqrt s}+\frac \mu s+\theta'(s)\right )q\right ) + \tilde P_{2,M} (V_1 q + V_2 \bar q) \nonumber\\
&+&  \tilde P_{2,M}  (B(q))  +  \tilde P_{2,M} (  R^*(\theta'(s),s)  ). \label{Projtild2}
\end{eqnarray}
Using   the fact that $q(s)\in \Vg_A(s)$  for all $s \in[\tau, s_1]$, Corollary \ref{corprojPtls}, Lemma \ref{lembin} and Lemma \ref{lemma-prejection-R-*constant},  we can obtain some bounds for the terms in the right hand side of \eqref{Projtild2}:

\begin{eqnarray}
& &  \tilde P_{2,M} ( \partial_s q)  =  \partial_s \tilde q_2 \label{estima-'-tilde-q-2}\\
& & \left|  \tilde P_{2, M} ( \mathcal{L}_{\beta, \delta} q)  \right| \leq \frac{A^{M+1}}{s^{\frac{M+2}{4}}},\label{estima-L-beta-q-2-tilde}
\end{eqnarray}

 and 
\begin{eqnarray}
& & \tilde P_{2,M} \left(-i\left (\frac{\nu}{2\sqrt s}+\frac \mu s+\theta'(s)\right )q\right )\nonumber\\
& =&  \tilde q_2^2 \left\{   \frac{c_2 \delta}{\kappa}    \right\}\nonumber\\
&+&  \frac{\tilde q_2}{\sqrt{s}} \left\{ \frac{\nu \delta}{2} \right\} \nonumber\\
&+&\frac{\tilde q_2}{s} \left\{ \frac{\nu}{2}\left [ c_4\tilde D_{4,2}-(1+\delta^2)\frac{\nu}{2} +D_{2,2}+ \frac{\Theta^*_{2,0} c_2}{\kappa}    \right ] - \frac{\nu}{2} \left[  \frac{\tilde K_{2,4} D_{4,2}}{2} + \tilde L_{2,4} \tilde D_{4,2}\right] \right.\label{estima-P-i-nu-mu-q-2-tilde}\\
&+& \left. \mu \delta + \frac{c_2}{\kappa} R^*_{2,1} + \frac{\delta R_{0,1}^*}{\kappa}\right\}  \nonumber\\
&+&\frac{1}{s^{\frac{3}{2}}}\left\{ \frac{\nu}{2} R^*_{2,1}     \right\}\nonumber\\
&+&  \frac{1}{s^{2}}\left\{ \frac{\nu}{2}\left[  X_2 + c_4 [\tilde C_{4,2}R^*_{2,1} +\tilde R_{4,2}^*]  -D_{2,0} \tilde R_{0,1}\right]   - \frac{\nu}{2} \left[ \tilde K_{2,4} \left( \frac{C_{4,2} R^*_{2,1}}{2} + \frac{R^*_{4,2}}{2} \right)\right]   \right.\nonumber\\
& - & \left. \frac{\nu}{2} \left[   \tilde L_{2,4} \left( \tilde C_{4,2} R^*_{2,1} + \tilde R_{4,2}^*  \right)\right] + \mu R^*_{2,1} + \frac{R^*_{0,1} R^*_{2,1}}{\kappa}  \right\} + O\left( \frac{A^8}{s^{\frac{9}{4}}} \right). \nonumber
 \end{eqnarray}

+ Let us now give the estimate for $ \tilde P_{2,M}(V_1 q + V_2 \bar q) $: 

\begin{eqnarray}
& & \tilde P_{2,M} \left(  V_1 q + V_2 \bar q  \right) \nonumber \\
& = & \frac{1}{\sqrt{s} } \left\{     \tilde q_0 \tilde D_{2,0}   + q_2 \tilde C_{2,2}+ \tilde q_2  \tilde D_{2,2}  +  q_4 \tilde C_{2,4} + \tilde q_4 \tilde D_{2,4} \right\} \nonumber\\
& + & \frac{1}{s} \left\{   \tilde q_0 \tilde F_{2,0} + q_2 \tilde E_{2,2} + \tilde q_2 \tilde F_{2,2}        \right\} + O\left( \frac{A^6}{s^{\frac{9}{4}}}\right).\label{estima-P-V-1-2-tilde-2}
\end{eqnarray}

 + In the following, we aim at   estimating to  $\tilde P_{2,M}( B(q) )$: Using Lemma  \ref{lembin}, we deduce then
 \begin{eqnarray}
\tilde P_{2,M} (B(q)) & =&   \tilde B_2 \tilde q_2^2  +  B_1  \frac{\tilde q_2}{s} +\frac{B_2}{s^2} +O\left(  \frac{A^6}{s^{\frac{9}{4}}} \right),\label{estima--P-2-tilde-B-q}  
\end{eqnarray}
where  \begin{eqnarray*}
\tilde B_2   & = &   \frac{1}{\kappa} \left\{    4(p - \delta^2) - \delta \beta (6+ 4 p+ 2\delta^2)    \right\}, \\
B_1 & =& \frac{1}{\kappa} \left\{    (-7\delta^2 \beta + p \beta -6 \beta) R_{2,1}^* - (  p- \delta^2)  \tilde R_{0,1}^*         \right\},\\
B_2 & = & \frac{(R_{2,1}^*)^2}{s^2} \frac{(32-64 \delta \beta)}{8 \kappa}.
\end{eqnarray*}
 
 + Finally, we give the estimate of  $\tilde P_{2,M} (R^*,s)$: From Lemma \ref{lemma-prejection-R-*constant}, we get 
 
 \begin{eqnarray*}
 & & \tilde P_{2,M} (R^* )   =  \frac{\tilde R^*_{2,1}}{s}+   \frac{\tilde R^*_{2,2}}{s^{\frac{3}{2}}}    +  \frac{\tilde R_{2,3}}{s^2} +   \frac{\theta' (s) \kappa}{\sqrt{s}} \frac{- \delta b}{(p-1)^2}  +    \theta'(s) \kappa  \frac{(p+1) \delta [12 - 6 \delta \beta + 6\beta^2] b^2}{2(p-1)^4}  \frac{1}{s}  \\
&+& O\left(   \frac{\theta'(s)}{s^{\frac{3}{2}}} \right)         +      O\left( \frac{1}{s^{\frac{5}{2}}}\right).\\
 \end{eqnarray*}
 However, by the definition of $\nu$ and $a$,  it follows that
 $$    \tilde R^*_{2,1} = 0.  $$
 
\noindent
In addition to that, we use  \eqref{asymptotic-theta's-in-pro-theta-0} to deduce the following 
\begin{eqnarray}
& &\tilde P_{2,M} (R^* (\theta'(s)))  = \frac{1}{s^{\frac{3}{2}}} \left\{      - \frac{\delta b}{(p-1)^2} R_{0,1}^*    +\tilde R^*_{2,2}   \right\}  \nonumber   \\
 & +&  \frac{\tilde q_2}{\sqrt{s}} \left[  - c_2\frac{\delta b}{(p-1)^2} \right] \nonumber \\
 &+& \frac{\tilde q_2}{s} \left\{  -\frac{\delta b}{(p-1)^2}\left[ D_{0,2} - \frac{\nu}{2} L_{0,2} +\frac{\Theta^*_{0,0} c_2}{\kappa}  \right]   +\frac{b^2}{2(p-1)^4}  c_2 (p+1) \delta [ 12 - 6 \delta \beta + 6\beta^2    ]  \right\} \nonumber \\
& +& \frac{1}{s^2}     \left\{       -\frac{\delta b}{(p-1)^2} \left[ \frac{\nu}{2} (1 + \delta^2) \tilde R^*_{0,1} -\frac{\nu}{2} K_{0,2} R^*_{2,1} - D_{0,0} \tilde R^*_{0,1} + C_{0,2} R^*_{2,1}+R^*_{0,2}  +\frac{\Theta_{0,0}^* R_{0,1}^*}{\kappa} \right]      \right. \nonumber \\
& +&  \tilde R_{2,3}^*  +  \left.      \frac{    (p+1) \delta [ 12 - 6 \delta \beta + 6 \beta^2 ]   R^*_{0,1}    b^2}{2(p-1)^4}  \right\}  + O\left( \frac{A^4}{s^{\frac{9}{4}}}\right).\label{estima-R-*-2-tilde}
\end{eqnarray}

\begin{center}
\textbf{ODE of $ \tilde q_2$}
\end{center}

  Adding  estimates   \eqref{estima-'-tilde-q-2}, \eqref{estima-L-beta-q-2-tilde}, \eqref{estima-P-i-nu-mu-q-2-tilde},            \eqref{estima-P-V-1-2-tilde-2},  \eqref{estima--P-2-tilde-B-q}
  and \eqref{estima-R-*-2-tilde}, we obtain   the following 
  
  \begin{eqnarray}
\tilde q_2' &=&  \left\{  \frac{\nu \delta}{2} + \tilde D_{2,2} - c_2 \frac{\delta b}{(p-1)^2}     \right\} \frac{\tilde q_2}{\sqrt{s}}  \nonumber\\[0.2cm]
& +& \tilde q_2^2 \left\{  \frac{c_2 \delta}{\kappa}  +\frac{4(p-\delta^2) - \delta \beta (6 +4p  +2\delta^2)}{\kappa}    \right\} \nonumber\\[0.2cm]
& +&  \frac{1}{s^{\frac{3}{2}}} \left\{\frac{\nu}{2} R_{2,1}^* + \tilde C_{2,2}R^*_{2,1} - \tilde D_{2,0} \tilde R_{0,1}^*  -  \frac{\delta b}{(p-1)^2} R^*_{0,1}  + \tilde R^*_{2,2}  \right\} \nonumber\\[0.2cm]
&+& \frac{\tilde q_2}{s}\tilde H_1 +  \frac{\tilde H_2}{s^2} + O \left( \frac{A^8}{s^{\frac{9}{4}}}\right),\label{ode-principal-tilde-q-2}
\end{eqnarray}

where
\begin{eqnarray*}
& & \tilde H_1  =  \frac{\nu}{2}\left [ c_4\tilde D_{4,2}-(1+\delta^2)\frac{\nu}{2} +D_{2,2}+\Theta^*_{2,0} \frac{c_2}{\kappa}\right ] - \frac{\nu}{2} \left[  \frac{\tilde K_{2,4} D_{4,2}}{2} + \tilde L_{2,4} \tilde D_{4,2}\right] \\
&+& \mu \delta + \frac{c_2}{\kappa} R^*_{2,1} + \frac{\delta R_{0,1}^*}{\kappa} \\
& +&  \tilde D_{2,0} \left [  \frac{\nu \tilde L_{0,2}}{2} - \tilde D_{0,2} -\tilde \Theta^*_{0,0} \frac{c_2}{\kappa}\right ] 
+ \tilde C_{2,2}  \left[   D_{2,2}   -\frac{\nu}{2} (1+ \delta^2) + c_4 \tilde D_{4,2}  + \frac{\Theta^*_{2,0} c_2 }{\kappa}  \right] \\
&+&  \frac{\tilde C_{2,4} D_{4,2} }{2} + \tilde D_{2,4} \tilde D_{4,2} \\
& +&  \tilde F_{2,2}  + B_1 -\frac{\delta b}{(p-1)^2}\left[ D_{0,2} - \frac{\nu}{2} L_{0,2}  +\Theta^*_{0,0} \frac{c_2}{\kappa}\right]   +\frac{b^2}{2(p-1)^4}  c_2 (p+1) \delta [ 12 - 6 \delta \beta + 6\beta^2    ] ,\\[0.4cm]
& &\tilde H_2 = \frac{\nu}{2}\left[  X_2 + c_4 [\tilde C_{4,2}R^*_{2,1} +\tilde R_{4,2}^*]  - D_{2,0} \tilde R_{0,1}\right]   - \frac{\nu}{2} \left[  \tilde K_{2,4} \left( \frac{C_{4,2} R^*_{2,1}}{2} + \frac{R^*_{4,2}}{2} \right)\right] \\
& -&  \frac{\nu}{2} \left[  \tilde L_{2,4} \left( \tilde C_{4,2} R^*_{2,1} + \tilde R_{4,2}^*  \right)\right] + \mu R^*_{2,1} + \frac{R^*_{0,1} R^*_{2,1}}{\kappa} \\ 
& + &   \tilde D_{2,0} \left[  -\tilde X_0 + \frac{ \nu \tilde K_{0,2} R_{2,1}^*}{2}       -  \tilde C_{0,2} R^*_{2,1}         \right]    + \tilde C_{2,2} \left[   X_2 + c_4 ( \tilde C_{4,2} R^*_{2,1} + \tilde R_{4,2})  - D_{2,0} \tilde R_{0,1}^*   \right]      \\
& +&  \tilde C_{2,4} \left[   \frac{C_{4,2} R^*_{2,1}}{2}  + \frac{R^*_{4,2}}{2}  \right]   + \tilde D_{2,4} (   \tilde C_{4,2} R^*_{2,1} + \tilde R^*_{4,2} ) +\tilde E_{2,2} R_{2,1}^* - \tilde F_{2,0} \tilde R_{0,1}^*  \\
& +& \frac{1}{8\kappa} (32 - 64 \delta \beta) (R_{2,1}^*)^2 \\
&-&   \frac{\delta b}{(p-1)^2} \left[ \frac{\nu}{2} (1 + \delta^2) \tilde R^*_{0,1} -\frac{\nu}{2} K_{0,2} R^*_{2,1} - D_{0,0} \tilde R^*_{0,1} + C_{0,2} R^*_{2,1}+R^*_{0,2} + \frac{\Theta^*_{0,0} R^*_{0,1}}{\kappa} \right] \\
& + &   \tilde R_{2,3}^*  +  \frac{    (p+1) \delta [ 12 - 6 \delta \beta + 6 \beta^2 ]   R^*_{0,1}    b^2}{2(p-1)^4}.
\end{eqnarray*}
  
  \medskip
 We return to    ODE \eqref{ode-principal-tilde-q-2}, we need to cancel some coefficient in order to finish the proof.

 +\textbf{ Cancellation  of the coefficient of $\frac{1}{s}$}: It is equivalent to 
 $$ \tilde R^*_{2,1}  = 0,$$
  which holds because of  the definition of $\tilde R^*_{2,1}$, $\nu , a$ as in  \eqref{definitionq} and the critical condition  as in Definition  \ref{defini-cris-condition}
  
  \medskip
  +\textbf{ Cancellation  of the coefficient of the coefficient of  $  \frac{\tilde q_2}{\sqrt s}$}: It is equivalent to the following
  
$$ \frac{\nu \delta}{2} + \tilde D_{2,2} - c_2 \frac{\delta b}{(p-1)^2} =0,$$
which is alway  satisfied because of the definitions of $\tilde D_{2,2}, \nu, c_2 $ and the  critical condition  given Definition \ref{defini-cris-condition}.
\[p-\beta\delta(p+1)-\delta^2=0.\]
\medskip
+\textbf{ Cancellation of the coefficient of $\tilde q_2^2$}: It is due the critical condition.\\ 
+\textbf{ Cancellation of the coefficient of $\frac{1}{s^{3/2}} $}:  It is equivalent to the following 

\begin{equation}\label{equality-order-s-3-2-tilde-q-2}
\frac{\nu}{2} R_{2,1}^* + \tilde C_{2,2}R^*_{2,1} - \tilde D_{2,0} \tilde R_{0,1}^*  -  \frac{\delta b}{(p-1)^2} R^*_{0,1}  + \tilde R^*_{2,2}  = 0.
\end{equation}

In fact, plugging the definition of constants in the right hand side and using the critical condition, we obtain the equation satisfied by $b$. We note that the resolution of this equation gives us the $b_{cri}$, the \textbf{ same }obtained by the formal approach and given by formula \eqref{formula-b-formal}.

After all this cancellation the equation satisfied by $\tilde q_2$ become 

  \begin{eqnarray}
\tilde q_2' &=&   \frac{\tilde q_2}{s}\tilde H_1 +  \frac{\tilde H_2}{s^2} + O \left( \frac{A^8}{s^{\frac{9}{4}}}\right),\label{Node-principal-tilde-q-2}
\end{eqnarray}

\medskip
Using the critical condition, we can write $\tilde H_1$ as follows:

\begin{eqnarray}
 & & \tilde H_1  \label{defi-tilde-H-1}\\
  & &= -  \frac{1}{4} \frac{[ \delta^8+(1-2p)\delta^6+(36-8p-5p^2)\delta^4+(-6p+61p^2-6p^3)\delta^2+36p^4-6p^3-6p^2  
 ] }{   (-p+6-p^2)\delta^4+(-p+10p^2-p^3)\delta^2-p^2+6p^4-p^3}.\nonumber
\end{eqnarray} 

In addition to that, we claim to prove the following property: for all $ p >1$ and $\delta \in (- p_{cri}, p_{cri})$ defined as in Definition \ref{defini-cris-condition}

\begin{equation}\label{sign-tilde-H-1}
 \tilde H_1 (p, \delta) \leq  -\frac{3}{2}.
\end{equation}
Indeed, let us consider
\begin{eqnarray*}
\tilde H_1 + \frac{3}{2} &=& 	-\frac{1}{4} \frac{\delta^2 (  \delta^6+\delta^4- 2 \delta^4 p+p^2 \delta^2-2 p \delta^2+p^2 )}{   (-p+6-p^2)\delta^4+(-p+10p^2-p^3)\delta^2-p^2+6p^4-p^3}\\[0.4cm]
& =&-\frac{1}{4} \frac{\delta^2(1 + \delta^2) ( \delta^2 -p )^2 }{   (-p+6-p^2)\delta^4+(-p+10p^2-p^3)\delta^2-p^2+6p^4-p^3}\\
& =& \frac{1}{4} \frac{\delta^2(1 + \delta^2) ( \delta^2 -p )^2 }{   ((p-2) \delta^2 - p(2p-1))((p+3) \delta^2 + p (3p+1))} \leq 0,
\end{eqnarray*} 
provided that $ \delta \in (- p_{cri}, p_{cri})$.

\medskip
\noindent
 We recall $\tilde Q_2$ defined as in  \eqref{defi-as-tilde-Q-2}, then, we derive from \eqref{ode-principal-tilde-q-2} that 
\[\tilde Q_2' =  \tilde H_1 \frac{\tilde Q_2}{s}+\frac{\left (\tilde\Ag_2(\tilde H_1+1)+\tilde H_2\right )}{s^2}+O\left( \frac{A^8}{s^{9/4}} \right).
\]

As a matter of fact, in order to  end  this part, we will choose $\mu=\mu_{cri}$, which cancel order $\frac{1}{s^2}$ in the above ODE. As the calculation are lengthy and to keep our work in a reasonable length, we will only prove the existence and unicity of such $\mu$. Thus why we will try to find the contribution of $\mu$ in the term
\[\tilde\Ag_2(\tilde H_1+1)+\tilde H_2.\]
In fact, we claim to to write

\begin{equation}\label{proof-find-mu}
 \tilde\Ag_2(\tilde H_1+1)+\tilde H_2  =  a_0 (p,\delta) \mu + a_1 (p, \delta, \beta),
\end{equation}
where  $ a_0 \neq 0$.   We kindly direct the reader to see the proof of \eqref{proof-find-mu}  in   item $(b)$,  Appendix  \ref{cancellation-ODE-tilde-q-2}.  We now directly derive  from \eqref{proof-find-mu} that 
\label{mu}
\begin{equation}\label{finding-mu-cri}
\mu =   - \frac{a_1(p,\delta, \beta)}{a_0 (p,\delta)}  \equiv \mu_{cri},
\end{equation}
this gives us the existence and unicity of $\mu$  and we finally obtain
$$  \tilde Q_2' =   \frac{\tilde H_1}{s} \tilde Q_2+ O\left( \frac{A^8}{s^{9/4}} \right),   $$
which implies item $(iii)$ of Proposition \ref{propode} .

\bigskip

$-$ The proof of item $(ii)$ of Proposition \ref{propode}: We start to the ODE   of $\tilde q_0$. In fact, using again  the shrinking set $\Vg_A$,  Corollary \ref{corprojPtls}, Lemma \ref{lembin} and Corollary \ref{lemma-prejection-R-*constant}, since $q(s)\in \Vg_A(s)$  for all $s \in[\tau, s_1]$, we get

 \begin{eqnarray}
 \tilde q_0' &=&  \tilde q_0  + \frac{\delta \nu }{2 \sqrt{s}} \tilde{q}_0 -\frac{\nu \tilde K_{0,2} }{2 \sqrt{s}} q_2  - \frac{\nu \tilde L_{0,2} }{2 \sqrt{s}} \tilde q_2 + \tilde D_{0.0} \frac{\tilde q_0}{\sqrt{s}} + \tilde C_{0.2} \frac{q_2}{\sqrt{s}} +\tilde D_{0.2} \frac{\tilde q_2}{\sqrt{s}}\nonumber\\
 & +&  \frac{\tilde \Theta^{*}_{0,0}}{\kappa} \left( c_2 \frac{\tilde q_2}{\sqrt{s}} +\frac{R^*_{0,1}}{s^\frac{3}{2}}\right) + \frac{ \tilde R_{0,1}^*}{s} + \frac{ \tilde R_{0,2}^*}{s^{\frac{3}{2}}} + O\left( \frac{1}{s^{\frac{7}{4}}}\right). 
\label{odetildeq0}
 \end{eqnarray}
Let us recall $\tilde Q_0$ defined as in Definition \ref{defthess}  and using the fact that $ q \in V_A(s)$ and the ode \eqref{Node-principal-tilde-q-2}, we can   derive the following
$$  \left|    \tilde Q_0'(s)  -   \tilde Q_0 (s)      \right| \leq \frac{C}{s^\frac{3}{2}}  . $$
  Next, considering $ \tilde q_1$ and proceeding similarly as for $\tilde q_0$, we can deduce that
$$ \left|  \tilde q_1'(s)   - \frac{1}{2} \tilde q_1(s)  \right|  \leq \frac{C}{s^\frac{3}{2}}.  $$
This finishes the proof of  item $(ii)$.

\bigskip

+ The proof of item (iv):
 \textit{Estimates of $ q_1$, $ q_2$, $ q_j$ and $\tilde q_j$ for $3\leq j\leq M$:} 
We  first project  equation \eqref{eqq} on $\tilde h_j$ and $h_j$
 \begin{eqnarray*}
 \tilde q_j'  &=& \tilde P_{j,M}(\mathcal{L}_{\beta, \delta}) + \tilde P_{j,M}\left( -i\left( \frac{\nu}{2 \sqrt s} + \frac{\mu}{s }  + \theta'(s)\right) q \right) \\
 &+& \tilde P_{j,M}(V_1 q + V_2 \bar q) + \tilde P_{j,M}(B)+\tilde P_{j,M}( R^*),\\
  q_j'  &=&  P_{j,M}(\mathcal{L}_{\beta, \delta}) + P_{j,M}\left( -i\left( \frac{\nu}{2 \sqrt s} + \frac{\mu}{s }  + \theta'(s)\right) q \right)\\
  & +&  P_{j,M}(V_1 q + V_2 \bar q) +  P_{j,M}(B)+\tilde P_{j,M}( R^*).
 \end{eqnarray*}

 \noindent
 Using  the definition of the shrinking set $\Vg_A$ Definition \ref{defthess} and Corollary \ref{corprojPtls}, Lemma \ref{lembin}, Corollary \ref{lemma-prejection-R-*constant} from Part 1 and the fact that $q(s)\in \Vg_A(s)$, we see that for all $s\in [\tau,s_1]$, we have
\beqtn
\begin{array}{lll}
\dsp \left| q_1^{'}+\frac 12  q_1\right|\leq C\frac{A^3}{s^{\frac 32}},&&  \\
\dsp\left| q_{j}^{'}+\frac j2 q_j\right|\leq C\frac{A^{j-1}}{s^{\frac{j+1}{4}}},&\mbox{if}&5\leq j\leq M,\\
 \dsp\left|\tilde q_{j}^{'} +\frac{j-2}{2}\tilde q_j\right|\leq C\frac{A^{j-1}}{s^{\frac{j+1}{4}}},&\mbox{if}&5\leq j\leq M.
\end{array}
\eeqtn

In particular,  for $q_2, q_3, \tilde q_3, q_4$ and $\tilde q_4$, we have  the following
\begin{eqnarray}
 q_2' &=& -q_2 + \frac{R_{2,1}^*}{s} +   \frac{R_{2,2}^*}{s^{\frac{3}{2}}} + c_4 \tilde q_4  - \frac{\nu}{2\sqrt{s}} (1 + \delta^2)\tilde q_2 + \frac{D_{2,0} \tilde q_0}{\sqrt{s}} +  \frac{D_{2,2} \tilde q_2}{\sqrt{s}} \\
 & + & \left( C_{2,2} - \frac{\delta \nu}{2} \right) \frac{q_2}{\sqrt{s}} +\Theta^*_{2,0}\frac{c_2}{\kappa}\frac{\tilde q_2}{\sqrt{s}}+  \frac{\Theta^*_{2,0} R_{0,1}^*}{ \kappa}\frac{1}{s^{\frac{3}{2}}} +  O \left(\frac{1}{s^{\frac{7}{4}}} \right),\\
 q_3' & =& -\frac{3}{2} q_3  + O\left( \frac{A^2}{s^\frac{3}{2}}\right),\\
 \tilde  q_3' & =& -\frac{1}{2} q_3  + O\left( \frac{A^2}{s^\frac{3}{2}}\right),
\end{eqnarray}
and
\begin{eqnarray*}
q_4'  &=&   -2 q_4  +  C_{4,2}  \frac{q_2}{\sqrt{s}} + D_{4,2} \frac{\tilde q_2}{\sqrt{s}} + \frac{R^*_{4,2}}{s^{\frac{3}{2}}} +  O\left( \frac{A^6}{s^{\frac{7}{4}}}\right),\\
 \tilde  q_4'  &=&   - \tilde q_4    +       \tilde C_{4,2} \frac{q_2}{\sqrt{s}} + \tilde D_{4,2} \frac{\tilde q_2}{\sqrt{s}} +   \frac{\tilde R_{4,2}^*}{s^{\frac{3}{2}}} + O\left( \frac{1}{s^{\frac{7}{4}}}\right).
 \label{odeq2q4}
\end{eqnarray*}
Let us recall the definition of $\tilde Q_0, Q_4$ and $\tilde Q_4$ defined as in Definition \ref{defthess} and using  the odes in \eqref{odeq2q4},  and \eqref{Node-principal-tilde-q-2} and the fact that $q \in V_A(s)$,    we can derive the following odes:
\begin{eqnarray}
\left| Q_2'(s)  + Q_2(s) \right| & \leq  & \frac{C}{s^\frac{7}{4}},\\
\left| Q_4'(s)  + 2 Q_4(s) \right| & \leq  & \frac{CA^6}{s^\frac{7}{4}},\\
\left| \tilde  Q_4'(s)  + Q_4(s) \right| & \leq  & \frac{C}{s^\frac{7}{4}},
\end{eqnarray}
and
\begin{eqnarray*}
\left| q_3'  + \frac{3}{2} q_3\right| \leq \frac{CA^2}{s^\frac{3}{2}},\\
\left| \tilde q_3'  + \frac{1}{2} \tilde q_3\right| \leq \frac{CA^2}{s^\frac{3}{2}}.
\end{eqnarray*}
 It is similar to the proof in Nouaili and Zaag, \cite{NZARMA18} for the case $\beta = 0$.  We use these estimates to conclude the result.
 
 Finally,  this also finish the proof of Proposition \ref{propode}.  $\blacksquare$
\subsubsection{ The infinite dimensional part: $q_-$}
We proceed as section 5.2 form \cite{MZ07}
We have 2 parts:
\begin{itemize}
\item In Part 1, we project equation \eqref{eqq} to get equations satisfied by $ q_-$. 
\item In Part 2, we prove the estimate on $q_-$.
\end{itemize}
\textbf{Part 1: Projection of equation \eqref{eqq} using the projector $P_-$}
In the following, we will project equation \eqref{eqq} term by term.

\medskip

\textbf{ First term: $\frac{\pa q}{\pa s}$}\\
From \eqref{decomp3}, its projection is
\beqtn
P_-(\frac{\pa q}{\pa s})=\frac{\pa q_-}{\pa s}
\eeqtn
\textbf{ Second term: $\tilde \Lg_{\beta,\delta} q$}\\
From \eqref{eqq}, we have the following, 
\[P_-( \Lg_{\beta,\delta} q)= \Lg_\beta q_-+P_-[(1+i\delta)\Re q_-].\]
\textbf{ Third term: $-i\left( \frac{\nu}{2 \sqrt s} +   \frac \mu s+\theta'(s)\right)q$}\\
Since $P_-$ commutes with the multiplication by $i$, we deduce that
\[P_-[-i\left(  \frac{\nu}{2 \sqrt s} + \frac \mu s+\theta'(s)\right)q]=-i\left(  \frac{\nu}{2 \sqrt s} + \frac \mu s+\theta'(s)\right)q_-.\]
\textbf{ Fourth term: $V_1 q$ and $V_2 \bar q$}\\
We have the following:
\begin{lemma}[Projection of $V_1 q$ and $V_2 \bar q$]
The projection of $V_1 q$ and $V_2\bar q$ satisfy for all $s\geq 1$,
\beqtn
\left\| \frac{P_-(V_1q)}{1+|y|^{M+1}}\right\|_{L^\infty}\leq \left(\|V_1\|_{L^\infty}+\frac{C}{s^{1/2}} \right)\left\|\frac{q_-}{1+|y|^{M+1}}\right\|_{L^\infty}
+\sum_{n=0}^{M}\frac{C}{s^{\frac{M+1-n}{4}}}(|\hat q_n|+|\tilde q_n|),
\eeqtn
and the same holds for $V_2 \bar q$.
\label{PVneg}
\end{lemma}
Using the fact that $q(s)\in \Vg_A(s)$, we get the following
\begin{corollary}
For all $A\geq 1$, there exists $s_{14}(A)$ such that for all $s\geq s_{14}$, if $q(s)\in \Vg_A(s)$ then,

\[\dsp   \left\|\frac{P_-(V_1 q)}{1+|y|^{M+1}}\right\|_{L^\infty}\leq \|V_1\|_{L^\infty}\left\| \frac{q_-}{1+|y|^{M+1}} \right \|_{L^\infty}+C\frac{A^M}{s^{\frac{M+2}{4}}},\]

and the same holds for $V_2 q$.
\end{corollary}

\textit{Proof of Lemma \ref{PVneg}:} We just give the proof for $V_1 q$ since the proof for $V_2\bar q$ is similar.

From Subsection \ref{sectidecompq}, we write $q=q_++q_-$ and
\[P_-(V_1 q)=V_1 q-P_+(V_1q_-)+P_-(V_1 q_+).\]
Moreover, we claim that the following estimates hold
\[
\begin{array}{lll}
\left\|\frac{V_1 q_-}{1+|y|^{M+1}}\right\|_{L^\infty} &\leq& \|V_1\|_{L^\infty}\left\|  \frac{q_-}{1+|y|^{M+1}}\right\|_{L^\infty}\\
\left\|\frac{P_+(V_1 q_-)}{1+|y|^{M+1}}\right\|_{L^\infty} &\leq& \frac{C}{s^{1/2}}\left\|  \frac{q_-}{1+|y|^{M+1}}\right\|_{L^\infty}.
\end{array}
\]
Indeed, the first one is obvious. To prove the second one, we use \eqref{bdVi1} to show that 

\[|P_{n,M}(V_1 q_-)|+|\tilde P_{n,M}(V_1 q_-)|\leq \frac{C}{s^{1/2}}\left\| \frac{q_-}{1+|y|^{M+1}}\right\|_{L^\infty}.\]
To control $P_-(V_1 q_+)=\sum_{n\leq M}P_-(V_1( q_n  h_n+\tilde q_n\tilde h_n))$, we argue as follows.\\
If $M-n$ is odd, we take $k=\frac{M-1-n}{2}$ in \eqref{dcVi}, hence
\[
\begin{array}{ll}
P_-(V_1( q_n  h_n+\tilde q_n \tilde h_n))=\sum_{j=1}^{k}\frac{1}{s^{\frac j2}} P_-\left (W_{1,j}( q_n h_n+\tilde q_n \tilde h_n)\right )+P_-\left (( q_n  h_n+\tilde q_n \tilde h_n)\tilde W_{1,k}\right)
\end{array}
\]
Since $2k+n\leq M$, we deduce that $P_-\left (W_{1,j}( q_n h_n+\tilde q_n \tilde h_n)\right )=0$ for all $0\leq j\leq k$. Moreover, using that
\[|\tilde W_{1,k}|\leq C\frac{(1+|y|^{2k+2})}{s^{\frac{k+1}{2}}}\]
 and applying Lemma  \ref{lemA3}, we deduce that

\beqtn
\left\|\frac{P_-(V_1( q_n  h_n+\tilde q_n \tilde h_n))}{1+|y|^{M+1}}\right\|_{L^\infty}\leq C \frac{| q_n|+|\tilde q_n|}{s^{\frac{M+1-n}{4}}}.
\label{estiodd}
\eeqtn
If $M-n$ is even, we take $k=\frac{M-n}{2}$ in \eqref{dcVi} and use that
\[|\tilde W_{1,k}|\leq C\frac{1+|y|^{2k+1}}{s^{\frac k2+\frac 14}},\]
to deduce that \eqref{estiodd} holds. This ends the proof of Lemma \ref{PVneg}.$\blacksquare$

\medskip

\textbf{Fifth term: $B(q,y,s)$.}\\
Using \eqref{estiquadinn}, we have the following estimate from Lemmas \ref{lemA3} and \ref{lemfifB}.
\begin{lemma} For all $K\geq 1$ and $A\geq 1$, there exists $s_{15}(K,A)$ such that for all $s\geq s_{15}$, if $q(s)\in \Vg_A(s)$, then
\beqtn
\left\|\frac{P_-(B(q,y,s))}{1+|y|^{M+1}} \right\|_{L^\infty}\leq C(M)    \left [\left(\frac{A^{M+2}}{s^{\frac 14}}\right)^{\bar p}+\frac{A^{5}}{s^{\frac 12}} \right]\frac{1}{s^{\frac{M+1}{4}}},
\label{Bmoins}
\eeqtn
where $\bar p=\min(p,2)$.
\end{lemma}
\medskip
\textit{Proof: }The proof is very similar to the proof of the previous lemma. From Lemma \ref{lemfifB}, we deduce that for all $s$ there exists a polynomial $B_M$ of degree $M$ in $y$ such that for all $y$ and $s$, we have
\beqtn
|B-B_M(y)|\leq C\left[\left(\frac{A^{M+2}}{s^{\frac 14}}\right)^{\bar p}+ \frac{A^{[5+(M+1)^2]}}{s^{\frac 12}}\right]\frac{(1+|y|^{M+1})}{s^{\frac{M+1}{4}}}.
\label{Bmoin}
\eeqtn
Indeed, we can take $B_M$ to be the polynomial
\[\dsp B_M=P_{+,M}\left [ \sum_{l=0}^{M}\sum_{\begin{array}{l}  0\leq j,k\leq M+1\\2\leq j+k\leq M+1   \end{array}  }\frac{1}{s^{l/2}} \left [ B_{j,k}^{l}(\frac{y}{s^{\frac 14}}) q_{+}^{j}\bar q_{+}^{k}\right ] \right ]  .\]
Then the fact that $B-B_M(y)$ is controlled by the right hand side of \eqref{Bmoin} is a consequence of the following estimates in the outer region and in the inner region.\\
First, in the region $|y|\geq s^{\frac 14}$, we have from Lemma \ref{lemfifB},
\[|B|\leq C|q|^{\bar p}\leq C\left(\frac{A^{M+2}}{s^{\frac 14}}\right)^{\bar p}\]
and from the proof of Lemma \ref{lembin}, we know that for $0\leq n\leq M$,
\[|\tilde P_{n}(B_M(q,y,s))|+| P_{n,M}(B_M(q,y,s))|\leq C\frac{A^n}{s^{\frac{n+2}{4}}}.\]
Beside, in the region $|y|\leq s^{\frac 14}$, we can use the same argument as in the proof of Lemma \ref{lemfifB} to deduce that the coefficients of degree $k\geq M+1$ of the polynomial
\[\sum_{l=0}^{M}\sum_{\begin{array}{l}  0\leq j,k\leq M+1\\2\leq j+k\leq M+1   \end{array}  }\frac{1}{s^{l/2}} \left [ B_{j,k}^{l}(\frac{y}{s^{\frac 14}}) q_{+}^{j}\bar q_{+}^{k}\right ] -B_M,\]
is controlled by $C\frac{A^k}{s^{\frac k 4+\frac 12}}$
 and hence
 \[\left | \sum_{l=0}^{M}\sum_{\begin{array}{l}  0\leq j,k\leq M+1\\2\leq j+k\leq M+1   \end{array}  }\frac{1}{s^{l/2}} \left [ B_{j,k}^{l}(\frac{y}{s^{\frac 14}}) q_{+}^{j}\bar q_{+}^{k}\right ] -B_M\right|\leq C\frac{A^{2M+2}}{s^{\frac{M+3}{4}}}(1+|y|^{M+1}),\]
in the region $|y|\leq s^{\frac 14}$.\\
Moreover, using that $ |q|\leq C\frac{A^{M+1}}{s^{\frac 14}}$ in the region $|y|\leq s^{\frac 14}$, we deduce that for all $s\geq 2 A^2$, we have 
\[\left|\sum_{l=0}^{M}\sum_{\begin{array}{l}  0\leq j,k\leq M+1\\2\leq j+k\leq M+1   \end{array}  }\frac{1}{s^{l/2}} \left [ B_{j,k}^{l}(\frac{y}{s^{\frac 14}}) q_{+}^{j}\bar q_{+}^{k}\right ] \right|\leq  C\frac{A^{2M+2}}{s^{\frac{M+3}{4}}}(1+|y|^{M+1}).\]
Finally, to control the term $|q|^{M+2}$, we use the fact that in the region $|y|\leq s^{\frac 14}$, we have the following two estimates $|q|\leq C\frac{A^{M+1}}{s^{\frac 14}}$
 and $|q|\leq A^5 \frac{1}{s^{\frac 34}}(1+|y|^{M+1})$ if $\sqrt {s}\geq 2 A^2$. Hence
 \[|q|^{M+2}\leq C\frac{A^5}{s^{\frac 34}}\left(\frac{A^{M+1}}{s^{\frac 14}}\right)^{M+1}(1+|y|^{M+1}).\]
This ends the proof of estimate \eqref{Bmoin} and conclude the proof of \eqref{Bmoins} by applying Lemma \ref{lemA3}.$\blacksquare$

\medskip

\textbf{Sixth term: $R^*(\theta',y,s)$}.\\
We claim the following:
\begin{lemma}
If $|\theta'(s)|\leq \frac{CA^5}{s^{3/2}}$, then the following holds
\[\dsp \left\|\frac{P_-(R^*(\theta',y,s))}{1+|y|^{M+1}}\right\|\leq C\dsp\frac{1}{s^{\frac{M+3}{4}}}\]
\end{lemma}
\textit{Proof: }Taking $n=\frac M2+1$ (remember $M$ is even), we write from Lemma \ref{decompR}
$R^*(\theta',y,s)=\Pi_n(\theta',y,s)+\tilde \Pi_n(\theta',y,s)$. Since $2n-2=M$, we see from subsection \ref{sectidecompq} that
\beqtn
\dsp |\tilde\Pi_n(\theta',y,s)|\leq\dsp C\frac{1+|y|^{2n-2}}{s^{\frac{n+1}{2}}}\leq\dsp C\frac{1+|y|^{M+1}}{s^{\frac{M+3}{4}}}
\label{RPIneq}
\eeqtn
in the region $|y|< s^{1/4}$. It is easy to see using \ref{Rexpect} and the definition of $\Pi_n$ that \eqref{RPIneq} holds for all $y\in\R$ and $s\geq 1$. Then applying Lemma \ref{lemA3}, we conclude easily.$\blacksquare$

\medskip

\textbf{Part 2: Proof of the last but one identity in (iv) of Proposition \ref{propode} (estimate on $q_-$)}\\ 
If we apply the projection $P_-$ to the equation \eqref{eqq} satisfied by $q$, we see that $q_-$ satisfies the following equation:
\[\dsp\frac{\pa q_-}{\pa s}=\Lg_{\beta,\delta} q_-+ P_-[(1+i\delta) \Re q_-]+P_-[-i(\frac\mu s+\theta'(s))q+ V_1 q+V_2\bar q+ B(q,y,s)+R^*(\theta',y,s)]. 
\]
Here, we have used the important fact that $P_-[(1+i\delta)\Re q_+]=0$.
As it was mentioned in \cite{MZ07}, in this part we use the operator $\Lg_{\beta}$, unlike the case for $q_n$ and $\tilde q_n$, where we used the properties of $\Lg_{\beta,\delta}$.
 The fact that $M$ is large as fixed in \eqref{boundM} is crucial in the proof.

Using the kernel of the semigroup generated by $\Lg_{\beta}$, we get for all $s\in [\tau,s_1]$,
\[\begin{array}{lll}
q_-(s)&=&\dsp e^{(s-\tau)\Lg_{\beta}}q_-(\tau)\\
&&+\dsp \int_{\tau}^{s}e^{(s-s')\Lg_{\beta}}P_-[(1+i\delta) \Re  q_-]ds'\\
&&+\dsp \int_{\tau}^{s}e^{(s-s')\Lg_{\beta}}P_-\left[-i(\frac\mu s+\theta'(s'))q+ V_1 q+V_2\bar q+ B(q,y,s')+R^*(\theta',y,s')\right]ds'.
\end{array}
\]
Using Lemma \ref{lemA2}, we get
\[\begin{array}{l}
\left\| \frac{q_-(s)}{1+|y|^{M+1}}\right\|_{L^\infty}\leq \dsp e^{-\frac{M+1}{2}(s-\tau)}\left\|\frac{q_-(\tau)}{1+|y|^{M+1}}\right\|_{L^\infty}\\
+\dsp \int_{\tau}^{s}e^{-\frac{M+1}{2}(s-s')} \sqrt{1+\delta^2} \left\|\frac{q_-}{1+|y|^{M+1}}\right\|_{L^\infty}ds'\\
+\dsp \int_{\tau}^{s}e^{-\frac{M+1}{2}(s-s')} \left\|\frac{P_-\left[-i(\frac\mu s+\theta'(s'))q+ V_1 q+V_2\bar q+ B(q,y,s')+R^*(\theta',y,s')\right]}{1+|y|^{M+1}}\right\|_{L^\infty}ds'
\end{array}
\]
Assuming that $q(s')\in \Vg_A(s')$, the results from Part 1 yields (use (i) of Proposition \ref{propode} to bound $\theta'(s)$)
\[\begin{array}{l}
\left\| \frac{q_-(s)}{1+|y|^{M+1}}\right\|_{L^\infty}\leq \dsp e^{-\frac{M+1}{2}(s-\tau)}\left\|\frac{q_-(\tau)}{1+|y|^{M+1}}\right\|_{L^\infty}\\
+\dsp \int_{\tau}^{s}e^{-\frac{M+1}{2}(s-s')} \left(\sqrt{1+\delta^2}+\||V_1 |+|V_2|\|_{L^\infty} \right)\left\|\frac{q_-}{1+|y|^{M+1}}\right\|_{L^\infty}ds'\\
+\dsp C(M)\int_{\tau}^{s}e^{-\frac{M+1}{2}(s-s')} \left [    \frac{A^{(M+1)^2+5}}{(s')^{\frac{M+3}{4}}}   +\frac{A^{(M+2)\bar p}}{(s')^{\frac{\bar p-1}{2}}}\frac{1}{(s')^{\frac{M+2}{2}}}    +\frac{A^M}{(s')^{\frac{M+2}{2}}}              \right]ds'.
\end{array}
\]
Since we have already fixed $M$ in \eqref{boundM} such that
\[\dsp M\geq 4\left( \sqrt{1+\delta^2}+1+2 \max_{i=1,2,y\in \R,s\geq 1}|V_i(y,s)|\right),\]
using Gronwall's lemma or Maximum principle and  \eqref{taus}, we deduce that
\[\begin{array}{ll}
e^{\frac{M+1}{2}s}\left\| \frac{q_-(s)}{1+|y|^{M+1}}\right\|_{L^\infty}&\leq\dsp  e^{\frac{M+1}{4}(s-\tau)}e^{\frac{M+1}{2}\tau}\left\| \frac{q_-(\tau)}{1+|y|^{M+1}}\right\|_{L^\infty}\\
&+\dsp e^{\frac{M+1}{2}s}  2^{\frac{M+3}{4}}   \left [    \frac{A^{(M+1)^2+5}}{s^{\frac{M+3}{4}}}   +\frac{A^{(M+2)\bar p}}{s^{\frac{\bar p-1}{2}}}\frac{1}{(s')^{\frac{M+2}{2}}}    +\frac{A^M}{s^{\frac{M+2}{2}}}              \right]
\end{array}
\]
which concludes the proof of the last but one identity in (iv) of Proposition \ref{propode}.

\subsubsection{The outer region: $q_e$}\label{outreg}
Here, we finish the proof of Proposition \ref{propode} by proving the last inequality in (iv). Since $q(s)\in \Vg_A(s)$ for all $s\in [\tau,s_1]$, it holds from Claim \ref{propshrinset} and Proposition \ref{propode} that
\beqtn
\|q(s)\|_{L^\infty(|y|<2K s^{1/4})}\leq C\frac{A^{M+1}}{s^{1/4}}\mbox{ and }|\theta'(s)|\leq \frac{CA^5}{s^{3/2}}.
\label{estheq}
\eeqtn
Then, we derive from \eqref{eqq1} an equation satisfied by $q_e$, where $q_e$ is defined by \eqref{defiqe}:
\beqtn
\begin{array}{ll}
\dsp \frac{\pa q_e}{\pa s}&= \Lg_{\beta} q_e-\frac{1}{p-1} q_e +(1-\chi)e^{\frac{i\delta}{p-1}s}\left\{L(q,\theta',y,s)+ R^*(\theta',y,s)\right\}\\
&\dsp -e^{\frac{i\delta}{p-1}s}q(s)\left( \pa_s \chi+(1+i\beta)\Delta \chi+\frac 12y\cdot\nabla \chi\right)+2 e^{\frac{i\delta}{p-1}s}(1+i\beta)div(q(s)\nabla \chi).
\end{array}
\eeqtn
Writing this equation in its integral form and using the maximum principle satisfied by $e^{\tau \Lg_{\beta}}$ (see Lemma \ref{lemA1}, see Appendix below), we write
\[
\begin{array}{ll}
\|q_e(s)\|_{L^\infty}& \leq \dsp e^{-\frac{s-\tau}{p-1}}\|q_e(\tau)\|_{L^\infty},\\
&+\dsp \int_{\tau}^{s}e^{-\frac{s-s'}{p-1}}\left(\|(1-\chi)L(q,\theta',y,s')\|_{L^\infty}+\|(1-\chi)R^*(\theta',y,s')\|_{L^\infty}\right)ds'\\
&+\dsp \int_{\tau}^{s} e^{-\frac{s-s'}{p-1}}\left\| q(s')\left(\pa_s \chi+(1+i\beta)\Delta \chi +\frac 12 y\cdot \nabla\chi\right)\right\|_{L^\infty}ds'\\
&+\dsp \int_{\tau}^{s} e^{-\frac{s-s'}{p-1}} \frac{1}{\sqrt{1-e^{-(s-s')}}}\| q(s')\nabla \chi\|_{L^\infty}ds'.
\end{array}
\]
Let us bounds the norms in the three last lines of this inequality.\\
First from \eqref{defchi14} and \eqref{estheq}
\beqtn
\begin{array}{ll}
\left\|q(s')\left(\pa_s \chi+(1+i\beta)\Delta \chi +\frac 12 y\cdot \nabla\chi\right)\right\|_{L^\infty}&\leq C(1+\frac{1}{K^2 s'})\|q(s')\|_{L^\infty(|y|<2 K {s'}^{1/4})}\\
&\dsp\leq C \frac{A^{M+1}}{{(s')}^{1/4}},
\label{estiqe1}
\end{array}
\eeqtn
\beqtn
\|q(s')\nabla \chi\|_{L^\infty}\leq \frac{C}{K {(s')}^{1/4}}\|q(s')\|_{L^\infty(|y|<2K {(s')}^{1/4})}\leq C\frac{A^{M+1}}{\sqrt{s'}},
\label{estiqe2}
\eeqtn
for $s'$ large enough.\\
Second note that the residual term $(1-\chi)R^*$ is small as well. Indeed, recalling the bound \eqref{estR*} on $R$, we write from the definition of $R^{*}$ \eqref{eqqd1} and \eqref{estheq}:
\beqtn
\|(1-\chi)R^*(\theta',y,s')\|_{L^\infty}\leq \frac{C}{{(s')}^{1/4}}+|\theta'(s')|\leq\frac{C}{{(s')}^{1/4}}
\label{estiqe3}
\eeqtn
for $s'$ large enough.\\
Third, the term $(1-\chi)L(q,\theta',y,s')$ given in \eqref{eqqd1} is less than $\epsilon |q_e| $ with $\epsilon =\frac{1}{2(p-1)}$. Indeed, it holds from \eqref{estheq} that: 
\beqtn
\begin{array}{l}
\dsp \|(1-\chi)L(q,\theta',y,s')\|_{L^\infty}\\
\leq \dsp C\|q_e(s')\|_{L^\infty}\left( \|\varphi(s')\|_{L^\infty(|y|\geq K{s'}^{1/4})}^{p-1} +\|q(s')\|_{L^\infty(|y|\geq K{s'}^{1/4})}^{p-1}+\frac{1}{s'}+|\theta'(s)|\right),\\
\leq \frac{1}{2(p-1)}\|q_e(s')\|_{L^\infty},
\label{estiqe4}
\end{array}
\eeqtn
whenever $K$ and $s'$ are large (in order to ensure that $\|\varphi(s')\|_{L^{\infty}(|y|\geq K{s'}^{1/4})}$ is small).\\
Notice that it is only here that we need the fact that $K$ is big enough. Using estimates \eqref{estheq}, \eqref{estiqe1}, \eqref{estiqe2}, \eqref{estiqe3} and \eqref{estiqe4}, we write
\[\begin{array}{ll}
\|q_e(s)\|_{L^\infty}&\leq \dsp e^{-\frac{s-\tau}{p-1}}\|q_e(\tau)\|_{L^\infty}\\
&\dsp + \int_{\tau}^{s}e^{-\frac{s-s'}{p-1}}\left(\frac{1}{2(p-1)}\|q_e(s')\|_{L^\infty}+C\frac{A^{M+1}}{(s')^{\frac 14}}+C\frac{A^{M+1}}{(s')^{\frac 12}}\frac{1}{\sqrt{1-e^{-(s-s')}}}\right)ds'.

\end{array}
\] 
Using Gronwall's inequality or Maximum principle, we end-up with
\[\|q_e(s)\|_{L^\infty}\leq e^{-\frac{(s-\tau)}{2(p-1)}}\|q_e(\tau)\|_{L^\infty}+\frac{CA^{M+1}}{ \tau^{\frac 14}}(s-\tau+\sqrt{s-\tau}),\]
which concludes the proof of Proposition \ref{propode}.


\section{Single point blow-up and final profile}
In this section, we prove Theorem \ref{thm1}\label{pth1}.
Here, we use the solution of problem \eqref{eqq}-\eqref{eqmod} constructed in the last section to exhibit a blow-up solution of equation \eqref{GL} and prove Theorem \ref{thm1}.\\
(i) Consider $(q(s), \theta(s))$ constructed in Section \ref{existence} such that  \eqref{limitev} holds. From \eqref{limitev} and the properties of the shrinking set given in Claim \ref{propshrinset}, we see that $\theta(s)\to \theta_0$ as $s\to \infty$ such that
\beqtn
\label{estitetas}
|\theta(s)-\theta_0|\leq CA^{10}\int_{s}^{\infty}\frac{1}{\tau^{\frac 54}}d\tau\leq \frac{CA^{10}}{ s^{\frac 14}} \mbox{ and }\|q(s)\|_{L^\infty(\R)}\leq \frac{C_0(K,A)}{\sqrt s}.
\eeqtn
Introducing $w(y,s)=e^{i(\nu \sqrt s+\mu\log s+\theta(s))}\left( \varphi(y,s)+q(y,s)\right)$, we see that $w$ is a solution of equation \eqref{equa-w} that satisfies for all $s\geq \log T$ and $y\in\R$,
\[|w(y,s)-e^{i\theta_0+i\nu\sqrt s+i\mu\log s}\varphi(y,s)|\leq C\|q(s)\|_{L^\infty}+C|\theta(s)-\theta_0|\leq \frac{C_0}{s^{\frac 14}}.\] 
Introducing
\[u(x,t)=e^{-i\theta_0}\kappa^{i\delta}(T-t)^{\frac{1+i\delta}{p-1}}w\left(\frac{y}{\sqrt{T-t}},-\log(T-t)\right),\]
we see from \eqref{chauto} and the definition of $\varphi$ \eqref{defi-varphi} that $u$ is a solution of equation \eqref{GL} defined for all $(x,T)\in\R\times [0,T)$ which satisfies \eqref{profile-u}.\\
If $x_0=0$. It remains to prove that when $x_0\not =0$, $x_0$ is not a blow-up point. The following result from Giga and Kohn \cite{GKCPAM89} allows us to conclude:
\begin{prop}[Giga and Kohn - No blow-up under the ODE threshold] For all $C_0>0$, there is $\eta_0>0$ such that if $v(\xi,\tau)$ solves
\[\left |  v_t -\Delta v\right |\leq C_0 (1+|v|^p)\]
and satisfies
\[|v(\xi,\tau)|\leq \eta_0(T-t)^{-1/(p-1)}\]
for all $(\xi,\tau)\in B(a,r)\times[T-r^2,T)$ for some $a\in \R$ and $r>0$, then $v$ does not blow up at ($a,T$).
\label{propositionGK}
\end{prop}
\textit{Proof: } See Theorem 2.1 page 850 in \cite{GKCPAM89}. $\blacksquare$\\

Indeed, we see from \eqref{profile-u} and \eqref{defi-varphi} that
\[\sup_{|x-x_0|\leq\frac{|x_0|}{2}}(T-t)^{\frac{1}{p-1}}|u(x,t)|\leq \left |  \varphi_0\left( \frac{|x_0|/2}{\sqrt{(T-t)}|\log(T-t)|}\right)\right |+\frac{C}{|\log (T-t)|^{\frac 14}}\to 0\]
as $t\to T$, $x_0$ is not a blow-up point of $u$ from Proposition \ref{propositionGK}. This concludes the proof of (i) of Theorem \ref{thm1}.\\
(ii)  Arguing as Merle did in \cite{FMCPAM92}, we derive the existence of a blow-up profile $u^*\in C^2(\R^*)$ such that $u(x,t)\to u^*(x)$ as $t\to T$, uniformly on compact sets of $\R^*$. The profile $u^*(x)$ is not defined at the origin. In the following, we would like to find its equivalent as $x\to 0$ and show that it is in fact singular at the origin. We argue as in Masmoudi and Zaag \cite{MZ07}. Consider $K_0>0$ to be fixed large enough later. If $x_0\neq 0$ is small enough, we introduce for all $(\xi,\tau)\in \R\times [-\frac{t_0(x_0)}{T-t_0(x_0)},1)$, 

\begin{align}
v(x_0,\xi,\tau)&=(T-t_0(x_0))^{\frac{1+i\delta}{p-1}} v(x,t),\\
\mbox{where,   }x&=x_0+\xi\sqrt{T-t_0(x_0)},\; t=t_0(x_0)+\tau(T-t_0(x_0)),
\label{defVN}
\end{align}

and $t_0(x_0)$ is uniquely determined by
\beqtn
|x_0|=K_0\sqrt{(T-t_0(x_0))|\log(T-t_0(x_0))|^{\frac 12}}.
\label{xt0N}
\eeqtn
From the invariance of problem (\ref{GL}) under dilation, $v(x_0,\xi,\tau)$ is also a solution of (\ref{GL}) on its domain. From (\ref{defVN}), \eqref{xt0N}, \eqref{defi-varphi}, we have
\[\sup_{|\xi|<2|\log(T-t_0(x_0))|^{1/8}}\left |v(x_0,\xi,0)-\varphi_0(K_0) \right|\leq \frac{C}{|\log(T-t_0(x_0))|^{\frac 18}}\to 0\mbox{ as }x_0\to 0. \]
Using the continuity with respect to initial data for problem (\ref{GL}) associated to a space-localization in the ball $B(0,|\xi|<|\log(T-t_0(x_0))|^{1/8})$, we show as in Section 4 of \cite{ZAIHPANL98} that
\[
\begin{array}{l}
\sup_{|\xi|\leq |\log(T-t_0(x_0))|^{1/8},\;0\leq\tau<1}|v(x_0,\xi,\tau)-U_{K_0}(\tau)|\leq \epsilon(x_0)\mbox{ as }x_0\to 0, \\
\end{array}
\]
where $U_{K_0}(\tau)=((p-1)(1-\tau)+b K_{0}^{2})^{-\frac{1+i\delta}{p-1}}$ is the solution of the PDE (\ref{GL}) with constant initial data $\varphi_0(K_0)$. Making $\tau \to 1$ and using (\ref{defVN}), we see that 
\[
\begin{array}{lll}
u^*(x_0)=\lim_{t\to T} v(x,t)&=&(T-t_0(x_0))^{-\frac{1+i\delta}{p-1}}e^{i\nu\sqrt {|\log (T-t_0(x_0))|}} |\log (T-t_0(x_0))|^{i\mu} \lim_{\tau \to 1} v(x_0,0,\tau)\\
&\sim&(T-t_0(x_0))^{-\frac{1+i\delta}{p-1}} e^{i\nu\sqrt {|\log (T-t_0(x_0))|}}|\log (T-t_0(x_0))|^{i\mu}U_{K_0}(1)
\end{array}
\]
as $x_0\to 0$. Since we have from (\ref{xt0N}) 
\[\log(T-t_0(x_0)) \sim 2 \log |x_0|\mbox{ and } T-t_0(x_0)\sim \frac{|x_0|^2}{\sqrt 2K_{0}^{2}\sqrt{|\log|x_0||}},\]
as $x_0\to 0$, this yields (ii) of Theorem \ref{thm1} and concludes the proof of Theorem \ref{thm1}. 

\medskip

Proof of Theorem \ref{Newthm}, we proceed as in the proof of (i) of Theorem \ref{thm1} and using results established in Section \ref{existence} we end the proof.

\appendix
\section{Spectral properties of $\Lg_{\beta}$}
In this Appendix, we recall from Appendix A of \cite{MZ07} some properties associated to the operator $\Lg_{\beta}$, defined in \eqref{eqqd1}. We recall that:
\[
  \Lg_{\beta} v=(1+i\beta)\Delta v-\frac 12 y\cdot \nabla v=\frac{1}{\rho_\beta}div(\rho_\beta\nabla w).
\]
where 
\[\rho_{\beta}(y)=\frac{e^{-\frac{|y|^2}{4(1+i\beta)}}}{(4\pi (1+i\beta))^{N/2}}.\]
Moreover, the operator $\Lg_\beta$ is self adjoint with respect to the weight  $\rho_{\beta}$ in the sense that
\beqtn
\int_{\R^N}u(y)\Lg_{\beta}w(y)\rho_{\beta}(y)dy=\int_{\R^N}w(y)\Lg_{\beta}u(y)\rho_{\beta}(y)dy.
\eeqtn
In one space dimension ($N=1$), the eigenfunction $f_n$ of $\Lg_{\beta}$ are dilations of standard Hermite polynomials $H_n(y)$:

\[f_n(y)=H_n(\frac{y}{2\sqrt {1+i\beta}})\mbox{,where } \Lg_{\beta} H_n=-\frac n 2 H_n.\]
If $N\geq 2$, its eigenfunction $f_{\alpha}(y_1,...,y_N)$ where
 $\alpha = (\alpha_1,...,\alpha_N)\in \N^N$ is a multi-indic are given by 
\[f_{\alpha}(y)=\Pi_{i=1}^{N}f_{\alpha_i}(y_i)=\Pi_{i=1}^{N}H_{\alpha_i}(\frac{y_i}{2\sqrt{1+i\beta}}).\]
The family $f_\alpha$ is orthogonal in the sense that for all $\alpha$ and $\xi\in \N^N$,
\[\int f_\alpha f_\xi \rho_\beta dy=\delta_{\alpha,\xi}\int f_{\alpha}^{2}\rho_\beta dy.\]

The semigroup generated by $\Lg_{\beta}$ is well defined and has the following kernel:
\beqtn
e^{s\Lg_\beta}(y,x)=\frac{1}{[4\pi(1+i\beta)(1-e^{-s})]^{N/2}}\exp\left [ -\frac{|x-ye^{-\frac s2}|^2}{4(1+i\beta)(1-e^{-s})}\right ].
\label{defsemgr}
\eeqtn
In the following, we give some properties associated to the kernel.
\begin{lemma}
a) The semigroup associated to $ \Lg_\beta$ satisfies the maximum principle:
\[
\|e^{s \Lg_\beta}\varphi\|_{L^\infty}\leq \|\varphi\|_{L^\infty}.
\]
b) Moreover, we have
\[\|e^{s \Lg_\beta} \div (\varphi)\|_{L^\infty}\leq \frac{C}{\sqrt{1-e^{-s}}}\|\varphi\|_{L^\infty},  \]
where $C$ only depends on $\beta$.
\label{lemA1}
\end{lemma}
\textit{Proof: } a) It follows directly by part, this also follows from the definition of the semigroup \eqref{defsemgr}.\\
b) Using an integration by part, this also follows from the definition of the semigroup \eqref{defsemgr}.
\begin{lemma}
There exists a constant $C$ such that if $\phi$ satisfies
\[\forall x\in\R\;\; |\phi(x)|\leq (1+|x|^{M+1})\]
then for all $y\in \R$, we have
\[|e^{s \Lg_\beta}P_-(\phi(y))|\leq Ce^{-\frac{M+1}{2}s}(1+|y|^{M+1})\]
\label{lemA2}
\end{lemma}
\textit{Proof: } This also follows directly from the semigroup's definition, through an integration by part, for a similar case see page 556-558 from \cite{BKN94}.

Moreover, we have the following useful lemma about $P_-$.
\begin{lemma}
For all $k\geq 0$, we have
\[\left\| \frac{P_-(\phi)}{1+|y|^{M+k}}    \right\|_{L^\infty}\leq C\left\|\frac{\phi}{1+|y|^{M+k}}     \right\|.\]
\label{lemA3}
\end{lemma}
\textit{Proof: } Using \eqref{eqQn}, we have
\[|\phi_n| \leq C   \left\|\frac{\phi}{1+|y|^{M+k}}\right\|_{L^\infty}. \]  
Since for all $m\leq M$, $|h_m(y)|\leq C(1+|y|^{m+k})$ and 
\[|\phi|\leq   C\left\|\frac{\phi}{1+|y|^{M+k}} \right\|_{L^\infty}(1+|y|^{m+k}),\]
the result follows from definition \eqref{decomp1} of $\phi$.

\section{Details of expansions of the potential terms: $V_1$ and $V_2$}\label{appendix-potential}

In this section, we aim at  giving some expansions of $V_1$ and $ V_2$  in  order to give the conclusion of item $(i)$ from Lemma \ref{lemma-projec-potential} and some related constants. Indeed, we recall the definition of $V_1$  and $ V_2$

 \begin{eqnarray*}
V_1 (y,s)  &=&  (1 + i\delta)  \frac{p+1}{2}  \left( |\varphi |^{p-1} - \frac{1}{p-1}  \right),\\
V_2 (y,s )  & =&  (1 +  i \delta) \frac{p-1}{2}  \left(  |\varphi|^{p-3} \varphi^2 - \frac{1}{p-1}\right),
  \end{eqnarray*}
where    
$$  \varphi (y,s)  = \varphi_0 ( y, s)   +  \frac{(1 + i \delta) a}{\sqrt{s}}=   \kappa \left(1  +  \frac{b}{p-1}  \frac{|y|^2}{ \sqrt{s}}    \right)^{-\frac{1}{p-1}}   +  \frac{(1 + i \delta) a}{\sqrt{s}},$$
and
$$  a =  \frac{2 \kappa b (1 - \delta \beta)}{(p-1)^2} . $$
Then, using Taylor expansion, we claim to the following asymptotic behaviors
\begin{equation}\label{expanding-of-V_1}
V_1 (y,s) =  \frac{1}{\sqrt{s}} W_{1,1} (y) +  \frac{1}{s} W_{1,2} + O \left( \frac{1 + |y|^6}{s^{\frac{3}{2}}} \right),
\end{equation}
and
\begin{equation}\label{expanding-of-V_2}
V_2 (y,s)  =  \frac{1}{\sqrt{s}} W_{2,1} (y)  +   \frac{1}{s} W_{2,2} (y) + O ( \frac{1 + |y|^6}{s^{\frac{3}{2}}}),
\end{equation}
where 

\begin{eqnarray*}
W_{1,1} (y)  &=&  (1 + i \delta) \frac{(p+1)}{2}  \frac{b}{(p-1)^2}  (- y^2  + 2 (1 - \delta \beta)  ),\\
W_{1,2} (y) & = &   (1 + i \delta)  \frac{(p+1)}{2} \frac{b^2}{(p-1)^3} \left\{ y^4 - \frac{(2 (1 - \delta \beta)(p-2 + \delta^2))}{p-1} y^2  \right.\\
& +& \left.    \frac{(p-1) (1 + \delta^2) (1 - \delta \beta)^2  + (p-3) (1 - \delta^2) (1 - \delta \beta)^2  }{ p-1}      \right\}\\
& =&   (1 + i \delta)  \frac{(p+1)}{2} \frac{b^2}{(p-1)^4}\left\{ (p-1) y^4 - [ 2 (1 - \delta \beta)(p-2 + \delta^2)] y^2  \right.\\
& +& \left.    2(p-2 + \delta^2)  (1 - \delta \beta)^2    \right\}
\end{eqnarray*}

 and
 
 \begin{eqnarray*}
& & W_{2,1} (y)   =    (1+ i\delta) \frac{b}{2  (p-1)^2}\{  (p-1 + 2 i \delta) ( - y ^2 + 2 (1 - \delta \beta))\},\\
& & W_{2,2} (y) =  (1 + i\delta) \frac{b^2}{ 2 (p-1)^4} \left\{  (p-1 + 2i\delta)(p-1 + i\delta) y^4      \right.\\
& - &    (2(p-1)(p-2 )  + (2p-10)\delta^2  + (8p -16) \delta i  )   (1 -\delta \beta) y^2\\
& +&  \left.   (1 - \delta \beta)^2 \left[  \frac{(p+1)(p-1)  (1 + i\delta)^2}{2} + (p+1) (p-3)(1 + \delta^2) + \frac{(p-3)(p-5) (1 - i \delta)^2}{2} \right] \right\}.
\end{eqnarray*}

For the proof  of \eqref{expanding-of-V_1} and \eqref{expanding-of-V_2}.  We give only  the second one because the first one can be obtained in the same way. Let us start  the proof \eqref{expanding-of-V_2}:  We first write 
$$  \varphi = \varphi_0 (z) + \frac{a (1 + i \delta)}{\sqrt s} \equiv  \varphi_0 + A.   $$

Using Taylor expansion of the function $\left| \varphi_0 + z  \right|^{p-3} ( \varphi_0 + z)^2 $ at $z = 0$, applying at $z = A$, we derive 
\begin{eqnarray*}
 | \varphi |^{p-3} \varphi^2  &=& |\varphi_0|^{p-3} \varphi_0^2   + \frac{p+1}{2} | \varphi_0|^{p-3} \varphi_0 A +  \frac{(p-3)}{2} | \varphi_0|^{p-5} \varphi_0^3  \bar A +   \frac{(p+1)(p-3) }{8}  |\varphi_0|^{p-3} A^2  \\
  &  +&   \frac{(p+1)(p-3) }{4} |\varphi_0 |^{p-5} \varphi_0^2  |A|^2 + \frac{(p-3)(p-5)}{8} |\varphi_0|^{p-7} \varphi_0^4 \bar{A}^2 + O \left( \frac{1}{s^\frac{3}{2}}\right).
\end{eqnarray*}
Besides that, we also have 
\begin{eqnarray*}
|\varphi_0|^{p-3} \varphi_0^2 & =&     \frac{1}{p-1} - \frac{(p-1 + 2 i \delta) b y^2}{ (p-1)^3 \sqrt{s}}  + \frac{(p-1 + 2i \delta) (p-1 + i \delta) b^2 y^4}{ (p-1)^5 s}  \\
&+&  O (\frac{1 + |y|^6}{s^{\frac{3}{2}}}),\\
\frac{ (p+1)}{2} |\varphi_0 |^{p-3} \varphi_0  A & = &   \frac{(p+1) (1 + i\delta) a}{2 \kappa  (p-1) \sqrt{s}}  -  \frac{(p+1) (1+ i \delta) (p-2 + i\delta) a b  y^2}{2 \kappa (p-1)^3 s} \\
&+& O \left( \frac{1 + |y|^6}{s^{\frac{3}{2}}}\right),\\
\frac{(p-3)}{2}  | \varphi_0|^{p-5} \varphi_0^3 \bar A & = &   \frac{(p-3) (1 - i\delta) a}{ 2 \kappa  (p-1) \sqrt{s}} -  \frac{(p-3) (1 - i\delta)  (p-2 + 3i \delta ) a b y^2}{2 \kappa  (p-1)^3 s }\\
& + & O \left( \frac{1  +|y|^6}{s^{\frac{3}{2}}}\right),\\
\frac{(p+1) (p-3)}{8} |\varphi_0|^{p-3} A^2 &=&  \frac{(p+1) (p-3)}{8} \frac{1}{\kappa^2 (p-1)} \frac{(1 + i\delta)^2 a^2}{s} +  O \left( \frac{1 + |y|^2}{s^{\frac{3}{2}}}\right), \\
\frac{ (p+1) (p-3)}{4} |\varphi_0|^{p-5} \varphi_0^2 |A|^2& =& \frac{ (p+1) (p-3)}{4} \frac{1}{ \kappa^2 (p-1)}  \frac{ (1 + \delta^2) a^2}{s} + O \left( \frac{1 + |y|^2}{s^{\frac{3}{2}}}\right),  \\
   \frac{(p-3)(p-5)}{8} |\varphi_0|^{p-7} \varphi_0^4 \bar{A}^2 &=&  \frac{(p-3)(p-5)}{8} \frac{1}{\kappa^2 (p-1)} \frac{(1 - i\delta)^2 a^2 }{s} + O \left( \frac{1 + |y|^2}{s^{\frac{3}{2}}}\right). 
\end{eqnarray*}
Thus, plugging these asymptotic in the formula of $V_2$,  we can derive \eqref{expanding-of-V_2}.  

\medskip
\noindent
In addition to that, we aim at determining the  constants given in item $(ii) $ of Lemma \ref{lemma-projec-potential}:

\begin{eqnarray*}
\tilde D_{4,2} & =&  \tilde P_{4,M} (W_{1,1} \tilde h_2  + W_{2,1} \bar{\tilde h}_2  )=\frac{b (\delta^2 - p)}{(p-1)^2},
\\
D_{2,2}  &=&  P_{2,M} (  W_{1,1} \tilde h_2 + W_{2,1} \bar{\tilde h}_2 ) \\
& =& -\frac{b}{2 (p-1)^2} \{  -24 p \delta  + 56 \delta^3 + 64 \delta^2 \beta +32 \delta +24 p\delta^2 \beta + 40 \delta^4 \beta  \},\\
\tilde D_{2,2} &=& \tilde P_{2,M}( W_{1,1} \tilde h_2 + W_{2,1} \bar{\tilde h}_2   ) = \frac{4b}{p-1} \delta \beta(1 + \delta^2),\\
 \tilde L_{2,4} & =&  \tilde P_{2,M} (i \tilde h_4) = 6\delta^2 \beta - 12 \delta -6 \beta,\\
D_{4,2}&   = &  P_{4,M} (W_{1,1} \tilde h_2 + W_{2,1} \bar{\tilde h}_2) = \frac{b}{ (p-1)^2}\left\{-2\delta(1+\delta^2) \right\}\\
\tilde D_{2,0} &=& \tilde P_{2,M} (W_{1,1} \tilde h_0 + W_{2,1} \bar{\tilde h}_0) = -\frac{b}{2(p-1)^2} ( 2p - 2 \delta^2)
\\
 \tilde L_{0,2} &=& \tilde P_{0,M}(i \tilde h_2)= -2\delta +\delta^2 \beta - \beta,
\\
 \tilde D_{0,2} & = & \tilde P_{0,M} (W_{1,1} \tilde h_2 + W_{2,1} \bar{\tilde h}_2) \\
 &=&-\frac{b}{2(p-1)^2} \{   -32 \delta \beta -12 p \beta^2  + 12 \delta^2 \beta^2 -16 \delta^2 +16 p - 4 \delta^4 \beta^2 +4 p\delta^2 \beta^2 - 32 p\delta \beta \},\\
 \tilde C_{2,2} & =  &  \tilde P_{2,M} (W_{1,1} h_2 + W_{2,1} \bar{h}_2)
 \\
 & = & -\frac{b}{2(p-1)^2} \{  -14 \delta^2 \beta +2p \beta -12 \beta   \},\\
  \tilde C_{2,4} & =& \tilde P_{2,M} (W_{1,1} h_4 + W_{2,1} \bar{h}_4) \\
 & =& - \frac{b}{2(p-1)^2}  \{ 96p \beta +224 \delta^3 \beta^2 -288\delta^2 \beta -128p\delta \beta^2-192\beta +96\delta \beta^2    \}\\
  \tilde D_{2,4} &=& \tilde P_{2,M} ( W_{1,1} \tilde h_4 + W_{2,1} \bar{\tilde h}_4  )\\
 & & \hspace{-1.5cm}= -\frac{b}{2(p-1)^2} \{  -96 p \delta^2 \beta^2 - 168  p \delta \beta + 96  p - 528 \delta \beta - 96 \delta^2 +216 \delta^2 \beta^2 -168 p \beta^2 +144 \delta^4 \beta^2 -360 \delta^3\beta  \}\\
\tilde F_{2,2} &  =   & \tilde P_{2,M} (  W_{1,2} \tilde h_2 + W_{2,2} \bar{\tilde h}_2) \\
 & =&   \frac{b^2}{2(p-1)^4} \left\{ -240 p +276 p^2 - 312p  \delta^2 - 204 \delta^4  + (-288p -552p^2 +696) \delta \beta   \right. \\
& +    & \left. (432-144p) \delta^3 \beta  + 144 \delta^5 \beta+(180p -180p^2) \beta^2 + (96p^2  + 288 p -96) \delta^2 \beta^2 + (108 + 36p) \delta^4 \beta^2 \right\},
\end{eqnarray*}
\begin{eqnarray*}
D_{0,2} &=& P_{0,M} (W_{1,1} \tilde h_2 + W_{1,2} \bar{\tilde h}_2)\\
& =& -\frac{b}{2(p-1)^2} \{  32 \delta +24 \delta^5 \beta^2 +64 \delta^2 \beta +48 \delta^3 \beta^2 +64 \delta^4 \beta +32 \delta^3 +24 \delta \beta^2  +96 p \delta^3 \beta^2 +96 p \delta \beta^2  \},\\[0.2cm]
\quad L_{0,2} & =& P_{0,M} ( i \tilde h_2 ) = 4 \delta \beta + 4 \delta^3 \beta.
\end{eqnarray*}

We would like to explain a little bit how we obtain these constants. For example, $\tilde D_{4,2} =  \tilde P_{4,M} (W_{1,1} \tilde h_2  + W_{2,1} \bar{\tilde h}_2  )$. Indeed,  we first write 
$$W_{1,1} \tilde h_2  + W_{2,1} \bar{\tilde h}_2 = \sum_{j=0}^k c_j y^j,$$
where $c_j \in \mathbb{C}$. Then, we  the definition of $h_j$ and $\tilde h_j$ to write

$$W_{1,1} \tilde h_2  + W_{2,1} \bar{\tilde h}_2 = \sum_{j=0}^k c_k y^j =   \sum_{j=0}^k ( a_j h_j + \tilde a_j \tilde h_j  ).   $$

Thus, we  obtain
\begin{eqnarray*}
P_{j,M} ( W_{1,1} \tilde h_2  + W_{2,1} \bar{\tilde h}_2) = a_j \text{ and } \tilde P_{j,M} ( W_{1,1} \tilde h_2  + W_{2,1} \bar{\tilde h}_2) = \tilde a_j.
\end{eqnarray*}

\section{Details of some expansions of  $B(q)$}\label{Appendix-B-q}

\medskip

In this section, we will give  the supplement of  the proof of     Lemma  \ref{lembin}  and some constants relating. Indeed, let us recall from \eqref{eqqd} that:
\[B(q,y,s)=(1+i\delta)\left ( |\varphi+q|^{p-1}(\varphi+q)-|\varphi|^{p-1}\varphi-|\varphi|^{p-1} q-\frac{p-1}{2}|\varphi|^{p-3}\varphi(\varphi \bar q+\bar\varphi q)\right).\]

\noindent
We first use the expansion in \eqref{expand-|a+a_0|p-1-a-a-0},  we first derive
 \begin{eqnarray*}
B(q) & = & \frac{(1 + i \delta)}{8 \kappa} \{ (p-3) \bar q^2 + 2 (p+1) \bar q q  + (p+1) q^2  \} \\
& = & \frac{(1+ i \delta)}{8 \kappa} \left\{ q_2^2  [  (p-3)  \bar h_2^2    + 2 (p+1)  \bar h_2  h_2  +   (p+1) h_2^2]      \right.\\
& +  &   \tilde q_2^2   [    (p-3)  \bar{\tilde h}_2^2 +  2 (p+1) \tilde h_2 \bar{\tilde h}_2 + (p+1) \tilde h_2^2   ]\\
& + & \tilde q_0 q_2 [   2  (p-3)  \bar{\tilde h}_0 \bar{h}_2 + 2 (p+1) ( \bar{\tilde h}_0 h_2  +  \tilde h_0 \bar{h}_2) + 2 (p+1) \tilde h_0 h_2]   \\
& +& \tilde q_0 \tilde  q_2 [   2  (p-3)  \bar{\tilde h}_0 \bar{\tilde h}_2 + 2 (p+1) ( \bar{\tilde h}_0 \tilde h_2  +  \tilde h_0 \bar{\tilde h}_2) + 2 (p+1) \tilde h_0 \tilde  h_2]  \\
&+  & \left.  q_2 \tilde q_2 [   2  (p-3)  \bar{\tilde h}_2 \bar h_2 + 2 (p+1) ( \bar{\tilde h}_2  h_2  +  \tilde h_2 \bar{ h}_2) + 2 (p+1)   h_2 \tilde h_2]        \right\}  + O(|q|^3).
\end{eqnarray*}

\noindent 
 Then, from then fact that $q \in \mathcal V_{A}(s)$, estimates \eqref{lembin2} and \eqref{lembin3} directly follow. In addition to that, we also derive the following  
\begin{eqnarray*}
\tilde P_{2,M}(B(q)) =  q_2^2 \tilde B_2 (q_2^2)  + \tilde q_2^2 \tilde B_2(\tilde q_2^2) + \tilde q_0 q_2 \tilde B_2(\tilde q_0q_2) + \tilde q_0 \tilde q_2 \tilde B_2(\tilde q_0 \tilde q_2) +  q_2 \tilde q_2 \tilde B_2 (q_2 \tilde q_2) + O \left(  \frac{1}{s^3}  \right).
\end{eqnarray*}
where

\begin{eqnarray*}
\tilde B_2(q_2^2)  &=& \frac{1}{8\kappa} \left(   32 - 64 \delta \beta      \right),\\
\tilde B_2 (  \tilde q_2^2 ) & =&    \frac{1}{\kappa} [   4(p-\delta^2) - \delta \beta ( 4p + 2 \delta ^2 + 6)  ], \\
\tilde B_2 (\tilde q_0 q_2)  &=& 0,\\
\tilde B_2 (\tilde q_0 \tilde q_2) & =& \frac{1}{8\kappa} [ 8 (p-\delta^2)  ]=\frac{1}{\kappa} [ p-\delta^2  ] , \\
\tilde B_2 (q_2 \tilde q_2)  & = & \frac{1}{8 \kappa} [   -56 \delta^2 \beta + 8 p \beta - 48 \beta  ] = \frac{1}{\kappa} [ -7 \delta^2 \beta  + p \beta  - 6 \beta ].
\end{eqnarray*}
 
 \noindent
Using the fact that $q \in \mathcal{V}_A(s)$,  we get
\beqtn
\tilde P_{2,M}(B)= \tilde B_2 \tilde q_{2}^{2}+B_1\frac{\tilde q_2}{s} +\frac{B_2}{s^2} +O \left(\frac{A^6}{s^{9/4}} \right),
\eeqtn
where
 \begin{eqnarray*}
\tilde B_2   & = &   \frac{1}{\kappa} \left\{    4(p - \delta^2) - \delta \beta (6+ 4 p+ 2\delta^2)    \right\}, \\
B_1 & =& \frac{1}{\kappa} \left\{    (-7\delta^2 \beta + p \beta -6 \beta) R_{2,1}^* - (  p- \delta^2)  \tilde R_{0,1}^*         \right\},\\
B_2 & = & \frac{(R_{2,1}^*)^2}{s^2} \frac{(32-64 \delta \beta)}{8 \kappa}.
\end{eqnarray*}

\section{Details of expansions of $R^*(y,s, \theta'(s))$ }\label{appendix-expansion-R}

Using the definition of $\varphi$, the fact that $\varphi_0$ satisfies \eqref{eqfi0} and \eqref{eqq}, we see that $R^*$ is in fact a function of $\theta'$, $z=\frac{y}{s^{1/4}}$ and $s$ that can be written as


\begin{eqnarray*}
R^* &=& \frac{(1 + i \beta)}{s^{\frac{1}{2}}} \Delta_z \varphi_0 (z) - \frac{1}{2} z \cdot \nabla \varphi_0 (z) - \frac{(1 + i \delta )}{p-1} \varphi_0 (z) - \frac{(1 + i \delta )^2}{p-1} \frac{a}{\sqrt s} \\
&+&  (1 + i\delta)F\left( \varphi_0 (z) + \frac{a}{\sqrt s} (1+i\delta) \right) + \frac{1}{4 s} z \cdot \nabla \varphi_0 + \frac{a}{2 s^{\frac{3}{2}}}(1+i\delta) \\
&-&  i \left( \frac{\nu}{2 \sqrt s} + \frac{\mu}{\sqrt s} + \theta'(s) \right) \varphi\\
&=& \frac{(1 + i \beta)}{s^{\frac{1}{2}}} \Delta_z \varphi_0 (z) + \frac{1}{4 s} z \cdot \nabla \varphi_0  + \frac{a }{2 s^{\frac{3}{2}}} (1+i\delta)-\frac{(1 + i \delta )^2 a}{(p-1) \sqrt s}   \\
&+& (1 + i \delta) \left(F\left( \varphi_0 (z) + \frac{a }{\sqrt s} (1+i\delta)\right) - F(\varphi_0(z)) \right) \\
&-&  i \left( \frac{\nu}{2 \sqrt s} + \frac{\mu}{ s} + \theta'(s) \right) \left( \varphi_0(z) + \frac{a}{\sqrt s} (1+i\delta)\right),\\
&=&  R_1^* (y,s)  + \theta'(s) \Theta(y,s), 
\end{eqnarray*}
where  $F(w) = |w|^{p-1} w$   and $\Theta(y,s) = - i \left( \varphi_0 (y,s) + \frac{a (1 + i \delta)}{\sqrt s} \right)$. 

\bigskip
\begin{center}
\textbf{ Expansion  of $R_1^* (y,s)$ in terms oh $h_j$ and $\tilde h_j$}
\end{center}

In fact, we   firstly  try to expand each term of $ R^*_1(y,s)$ as a power series in  $\frac {1}{s}$  and $y$   as $s\to \infty$, uniformly for $\dsp |y|\leq C\dsp s^{\frac 14}$, for some  $C$, a positive constant.  

As a matter of fact, we have the following

\begin{eqnarray}
& & |a_0 + u|^{p-1} (a_0 + u) = |a_0|^{p-1} a_0 + \frac{p+1}{2} |a_0|^{p-1} u + \frac{p-1}{2} |a_0|^{p-3} a_0^2 \bar u\nonumber\\
&+& \frac{(p-1)(p-3)}{8} |a_0|^{p-5} a_0^3 \bar u^2 + \frac{(p+1)(p-1)}{4} |a_0|^{p-3}a_0 u \bar u + \frac{(p+1)(p-1)}{8} |a_0|^{p-3} \bar a_0 u^2\nonumber\\
&+& \frac{(p-1)(p-3)(p-5)}{48} |a_0|^{p-7} a_0^4 \bar u^3 + \frac{(p+1)(p-1)(p-3)}{16} |a_0|^{p-5} a_0^2 u \bar u^2,\nonumber\\
&+& \frac{(p+1)(p-1)^2}{16}  |a_0|^{p-3} u^2 \bar u + \frac{(p+1)(p-1)(p-3)}{48} |a_0|^{p-5} \bar a_0^2 u^3 + O(u^4),\label{expand-|a+a_0|p-1-a-a-0}
\end{eqnarray}
as $u \to 0$.   

\medskip
  Let us  consider  $a_0 =  \varphi_0(z) = \kappa (1 + \frac{b}{p-1} z^2)^{- \frac{1 + i \delta}{p-1}} $, where $z=\frac{y}{\sqrt s}$ and  $ u = \frac{a }{\sqrt s}(1+i\delta)$, then  we have  in addition 
\begin{eqnarray*}
\frac{p+1}{2} |a_0|^{p-1} u &=& \frac{(p+1) a(1+i\delta)}{2(p-1) \sqrt s} \left( 1 + \frac{b}{p-1} z^2\right)^{-1} \\
&=& \frac{(p+1) a(1+i\delta)}{2(p-1) \sqrt s} \left( 1  - \frac{b}{p-1} z^2  + \frac{b^2 z^4}{ (p-1)^2} + O(z^6)\right),\\
\frac{p-1}{2} |a_0|^{p-3} a_0^2 \bar u &=& \frac{ a(1-i\delta)}{2  \sqrt s} \left( 1 + \frac{b}{p-1} z^2\right)^{- \frac{(p-1 + 2 i \delta)}{p-1}} \\
&=& \frac{ a(1-i\delta)}{2 \sqrt s}  \left( 1 - \frac{(p-1 + 2 i \delta) b}{(p-1)^2} z^2 + \frac{(p-1 + 2 i\delta)(p-1 + i\delta) b^2}{(p-1)^4} z^4  \right)\\
& +&  O(z^6)
\end{eqnarray*}
and 
\begin{eqnarray*}
\frac{(p-1)(p-3)}{8} |a_0|^{p-5} a_0^3  \bar u^2  &=& \frac{(p-3) (1-i\delta)^2 a^2 }{8 \kappa s} \left(1 + \frac{b}{p-1} z^2 \right)^{ - \frac{(p-2 + 3 i \delta)}{p-1}} \\
& =&  \frac{(p-3)(1-i\delta)^2  a^2 }{8 \kappa s}  \left( 1 - \frac{(p-2 + 3 i \delta)b z^2}{(p-1)^2} + O(z^4)\right),\\
\frac{(p+1)(p-1)}{4} |a_0|^{p-3}a_0 u \bar u &=& \frac{(p+1) (1+\delta^2) a^2}{4 \kappa s} \left(1 + \frac{b}{p-1} z^2 \right)^{- \frac{(p-2 + i \delta)}{p-1}} \\
&= & \frac{(p+1) (1+\delta^2) a^2}{4 \kappa s} \left(  1 - \frac{(p-2 + i \delta) b z^2}{(p-1)^2} + O(z^4) \right),\\
\frac{(p+1)(p-1)}{8} |a_0|^{p-3} \bar a_0 u^2 &=& \frac{(p+1)(1+i\delta)^2  a^2}{8 \kappa s} \left(1 + \frac{b}{p-1}z^2 \right)^{- \frac{(p-2 - i\delta)}{p-1}}\\
& = & \frac{(p+1) (1+i\delta)^2 a^2}{8 \kappa s}  \left( 1 - \frac{(p-2 - i \delta) b z^2}{(p-1)^2} + O(z^4)  \right),
\end{eqnarray*}

Besides that, we also get the following
\begin{eqnarray*}
& & \frac{1}{4 s} z \cdot \nabla \varphi_0 (z) \\
& = &  - \frac{ \kappa (1 + i \delta) b z^2}{2 (p-1)^2 s} + O\left( \frac{z^4}{s}\right),\\ 
& & \frac{(1 + i \beta )}{s^{\frac{1}{2}}} \Delta_z \varphi_0 (z) \\
& = &  -\frac{2 \kappa (1 + i \delta) (1 + i\beta)b}{(p-1)^2 \sqrt s} \left( 1 - \frac{(p+i \delta) b z^2}{(p-1)^2} + \frac{(p+i \delta)(2p -1+ i \delta) b^2 z^4}{2 (p-1)^4 } + O(z^6) \right) \\
&+& \frac{4 \kappa (1 + i \delta) (1 + i \beta)(p +i\delta) b^2 z^2}{(p-1)^4 \sqrt s} \left( 1- \frac{(2p -1 + i\delta) b z^2}{(p-1)^2} + O(z^4) \right)\\
&=& - \frac{2 \kappa (1 + i \delta) (1 + i \beta) b}{(p-1)^2 \sqrt s}  + i \frac{6  \kappa (p+1) \delta  (1 + \beta^2) b^2 y^2 }{(p-1)^4 s}\\
&+& \frac{5 \kappa (p+1) \delta (1 + \beta^2) ( \delta - i( 2p - 1)  ) b^3 y^4}{(p-1)^6 s^{\frac{3}{2}}} + O \left(\frac{1 + |y|^6}{s^2}\right),\\
& & \varphi_0 (z) +(1+i\delta) \frac{a }{s^{\frac{1}{2}}} \\
&=& \kappa  - \frac{\kappa (1 + i \delta) b z^2}{(p-1)^2} + \frac{\kappa (p+1) \delta (\beta + i) b^2 z^4}{2(p-1)^4} + (1+i\delta)\frac{a }{\sqrt s} + O(z^6).
\end{eqnarray*}

So,  in expansion of $R^*_1$ in series of $\frac{1}{s^j}$,   we have 

\medskip
\textbf{ + Order  $\frac{1}{\sqrt s}$: }
\begin{eqnarray*}
& &- \frac{2 \kappa (1 + i \delta) (1 + i \beta) b }{(p-1)^2} - \frac{(1 + i \delta)^2 a}{p-1} + (1 + i \delta) \left( \frac{(p+1)(1+i\delta)  a }{2 (p-1)} + \frac{ (1-i\delta) a}{2}\right) - i \frac{ \kappa \nu }{2} \\
&=&  - i \frac{ \kappa \nu }{2}     - \frac{2 \kappa (1 + i \delta) (1 + i \beta) b }{(p-1)^2} + (1+ i\delta) a,
\end{eqnarray*}

Then, we can write $R^*_1$ as follows

\begin{eqnarray*}
R_1^* (y,s) =  \frac{1}{\sqrt s}  \mathcal{R}_0 (y)  + \frac{1}{s} \mathcal{R}_1 (y) + \frac{1}{s^\frac{3}{2}} \mathcal{R}_2(y) + \frac{1}{s^2} \mathcal{R}_3(y) +  \tilde{\mathcal{R}}(y,s),
\end{eqnarray*}

where $ \tilde{\mathcal{R}} $ satisfies
$$ |  \tilde{\mathcal{R}} (y,s) |   \leq \frac{C(1+ |y|^8)}{s^{\frac{5}{2}}},$$

which implies that
$$   | P_{j,M} (\tilde{\mathcal{R}}) | +  | \tilde P_{j,M} (\tilde{\mathcal{R}}) | \leq \frac{C}{s^\frac{5}{2}}. $$

In addition to that, by using a explicit calculation, we can  show that $\mathcal{R}_j (y)$ is a polynomial of order  $2j$.    Besides that, these polynomials don't contain any odd order of $y$.

\medskip
 Moreover, we can   write  $\mathcal{R}_j$ in terms of $h_k, $ and $\tilde h_k$ as follows

\begin{eqnarray*}
\mathcal{R}_j(y) = \sum_{k=0}^{2j} ( R^*_{k,j} h_k +\tilde  R_{k,j}^*\tilde h_k  )
\end{eqnarray*}

More precisely, we have
\begin{eqnarray*}
\tilde R^{*}_{0,0} &=& a-2(1-\beta\delta ) \frac{b}{(p-1)^2}), \\
R^{*}_{0,0} &=& \kappa \frac{\nu}{2}-2\kappa\beta(1+\delta^2)\frac{b}{(p-1)^2},\\
\tilde R^{*}_{2,1} &=& \left( \frac{2(\delta^2 - p)a b - \kappa \delta \nu b}{2 (p-1)^2}\right)\\
R^{*}_{2,1} &=&\left( \frac{(1  + \delta^2)}{2 (p-1)^2} (\kappa \nu b - 4 \delta a b) +\frac{6 \kappa (p+1) \delta (1 + \beta^2) b^2}{(p-1)^4} \right)\\
\tilde R^{*}_{0,1} &=&  \frac{a^2(p -\delta^2) + \kappa \delta \nu a}{2 \kappa}  + \frac{(1 - \delta \beta)( 2 (\delta^2 - p)ab - \kappa \delta \nu b)}{(p-1)^2}\\
& -&  \frac{(1  + \delta^2)}{2 (p-1)^2} \beta(\kappa \nu b - 4 \delta a b) -\frac{6 \kappa (p+1) \delta \beta (1 + \beta^2) b^2}{(p-1)^4}\\
R^{*}_{0,1} &=& \frac{(1 + \delta \beta)(1 + \delta^2)(\kappa \nu b - 4 \delta a b)}{(p-1)^2} + \frac{12 \kappa (p+1) \delta(1 + \delta \beta)(1 + \beta^2)b^2}{(p-1)^4}\\
&+& \frac{(1 + \delta^2) (8 \delta a^2 - 4 \kappa \nu a)}{8 \kappa}  - \kappa \mu,\\
\tilde R_{2,2}^*  &=& \frac{5 \kappa (p+1) \delta (1 + \beta^2)b^3 }{(p-1)^6}   \left[ 12\delta - 6 \delta^2 \beta + 6 (2p -1) \beta \right]    \\
&  +&  \frac{\nu \kappa (p+1) \delta  b^2}{4(p-1)^4} [   12 - 6 \delta \beta + 6\beta^2] \\
& & \hspace{-1cm} + \frac{ab^2}{2 (p-1)^4} [ 24  p^2 -24 p+ (30 - 6p - 24p^2) \delta \beta -24p \delta^2 - 24 \delta^4 + (18 - 6p) \delta^3 \beta +12 \delta^5 \beta ] \\
&-& \frac{\kappa b}{ 2 (p-1)^2} -\frac{\mu \kappa \delta b }{(p-1)^2} 
-  \frac{a^2 b}{8 \kappa (p-1)^2} \left\{     4 p^2 - 8p - 12 \delta^4 \right\}.
\end{eqnarray*}

\label{page-constant-R^*-j-k}
Besides that, the other constants which will not give in this paper, because we need only  the existence of $\mu$. In addition to that, in order to make more clearly our computation, we also give below the calculations of  polynomial  $\mathcal{R}_j $. Then, the readers can check the accuracy of our constant.

\bigskip
+ Formula of $\mathcal{R}_0(y)$: As we mentioned above, $\mathcal{R}_0(y)$ is a polynomial of order $0$ in $y$,
\begin{eqnarray*}
& & \mathcal{R}_0(y) \\
& =&- \frac{2 \kappa (1 + i \delta) (1 + i \beta) b }{(p-1)^2} - \frac{(1 + i \delta)^2 a}{p-1} + (1 + i \delta) \left( \frac{(p+1)(1+i\delta)  a }{2 (p-1)} + \frac{ (1-i\delta) a}{2}\right) - i \frac{ \kappa \nu }{2} \\
&=&  - i \frac{ \kappa \nu }{2}     - \frac{2 \kappa (1 + i \delta) (1 + i \beta) b }{(p-1)^2} + (1+ i\delta) a,
\end{eqnarray*} 

+ Formula of $\mathcal{R}_1(y):$  As we mentioned, $\mathcal{R}_1(y)$ is a polynomial of order $2$

\begin{eqnarray*}
& & \mathcal{R}_1(y) =  \mathcal{R}_{1,0} +  \mathcal {R}_{1,2} y^2,
\end{eqnarray*}
where
\begin{eqnarray*}
& & \mathcal{R}_{1,0} \\
&=&(1 +  i \delta) \left( \frac{(p-3) (1  - i \delta)^2 a^2}{ 8 \kappa} + \frac{(p+1) (1 + \delta^2) a^2 }{4 \kappa} + \frac{(p+1) (1 + i \delta)^2 a^2}{8 \kappa}  - i \frac{  \nu a}{2}  \right) - i \mu \kappa \\
&=& (1 + i \delta) \left(  a^2 \frac{4p + 8 \delta i  + 4 \delta^2}{8\kappa} - i \frac{ 4 \kappa \nu a}{8 \kappa} \right) - i \mu \kappa \\
&= &  (1 + i \delta) \left(  \frac{a^2(p -\delta^2) + \kappa \delta \nu a}{2 \kappa} \right) \\
&+& i  \left(  \frac{(1 + \delta^2) (8 \delta a^2 - 4 \kappa \nu a)}{8 \kappa}  - \kappa \mu \right)\\
& &  \mathcal {R}_{1,2}  = \displaystyle i \frac{\kappa (1 + i \delta) \nu b }{2 (p-1)^2} + i \frac{6 \kappa (p+1) \delta (1 + \beta^2) b^2}{(p-1)^4}\\
 &+& (1 + i \delta) \left( \frac{(p+1) (1 + i \delta) a}{2 (p-1)} \left( - \frac{b}{p-1} \right)  \right) + \frac{(1 - i \delta) a }{2} \left( - \frac{(p-1 + 2 i \delta) b}{(p-1)^2} \right)
\end{eqnarray*}

+ Formula of  $\mathcal{R}_2(y)$: In fact, this polynomial is of order $4$ that we will give in the following

$$ \mathcal{R}_{2} (y)   =\mathcal{R}_{2,0} + \mathcal{R}_{2,1}  y^2 + \mathcal{R}_{2,2} y^4,$$

where 
\begin{eqnarray*}
\mathcal{R}_{2,1} & =& -\frac{\kappa (1 + i\delta) b}{2 (p-1)^2} \\
&+& (1 + i \delta) \left\{ -\frac{(p-3) (1-i\delta)^2 a^2}{8 \kappa}   \frac{(p-2 + 3i \delta ) b}{ (p-1)^2} \right.\\
& - & \left. \frac{(p+1) (1+ \delta^2)a^2}{4\kappa} \frac{(p-2 + i \delta) b}{(p-1)^2}  - \frac{(p+1) (1 + i \delta)^2 a^2 }{8 \kappa} \frac{(p-2 - i\delta) b}{(p-1)^2}  \right\}\\
& + &  \frac{i \mu \kappa (1 + i \delta) b}{(p-1)^2}\\
&=& -\frac{\kappa (1 + i\delta) b}{2 (p-1)^2} +  \frac{i \mu \kappa (1 + i \delta) b}{(p-1)^2}\\
& +& (1 + i\delta) \frac{a^2 b}{8 \kappa (p-1)^2}  \left\{  -(p-3) (1-i\delta)^2  (p-2 + 3i \delta )  \right.\\
& -& \left.    2 (p+1) (1+ \delta^2) (p-2 + i \delta)  - (p+1) (1 + i \delta)^2  (p-2 - i\delta)     \right\},
\end{eqnarray*}
 and
 \begin{eqnarray*}
 \mathcal{R}_{2,2}  &=& \frac{5 \kappa (p+1) \delta (1 + \beta^2)(\delta - i (2p -1)) b^3 }{(p-1)^6}\\[0.2cm]
 & +& \frac{\nu \kappa (p+1) \delta (1 - i\beta) b^2}{4(p-1)^4}\\
 & +& \frac{ab^2}{2 (p-1)^4} \left[  (p^2 -1)(1 + i \delta)^2 + (1 + \delta^2)(p-1 + 2i\delta )(p-1 + i \delta)   \right] \\
 &  =&  \frac{5 \kappa (p+1) \delta (1 + \beta^2)b^3 }{(p-1)^6}(\delta - i (2p -1)) \\[0.2cm]
 & +& \frac{\nu \kappa (p+1) \delta  b^2}{4(p-1)^4}(1 - i\beta) \\
 & +& \frac{ab^2}{2 (p-1)^4} \left\{    2 p^2 - 2p - 2p \delta^2 - 2 \delta^4 + i [(2p^2 + 3p - 5)\delta + (3p-3) \delta^3 ]      \right\},
 \end{eqnarray*}
and $\mathcal{R}_{2,0}$ is some constant in $\mathbb{C}$ which we don't need to calculate explicitly. Moreover, we will not need the exact formula of $\mathcal{R}_3(y)$, we just need to know that it  is a polynomial of order $6$.  Finally, in order to finish the calculation on $R^*(\theta'(s))$, it remains to expand the second term $\theta'(s) \Theta(y,s)$.

\begin{center}
\textbf{Expansion of $ \theta'(s) \Theta(y)$ }
\end{center}
We introduce
$$ \Theta (y,s) =  - i \left( \varphi_0 (y,s) + \frac{a (1+ i \delta)}{ \sqrt s} \right), $$
where  $\varphi_0 $ and $a$ defined as in  \eqref{defi-varphi} and \eqref{definitionq}, respectively. Using Taylor expansions, we write 
\begin{eqnarray*}
\Theta(y,s) & =&   -i\kappa - \kappa (\delta-i)\frac{y^2}{\sqrt s}\frac{b}{(p-1)^2} +a(\delta-i)\frac{1}{\sqrt s}\\
 &+& \kappa (1-i\beta)\delta (p+1)\frac{y^4}{ s}\frac{b^2}{2(p-1)^4}
 +\tilde \Theta (y,s),
\end{eqnarray*}

\noindent
where  $\tilde \Theta(y,s) $ satisfies the following
$$  \left|  \tilde \Theta (y,s) \right| \le \frac{C(1 + |y|^6)}{s^\frac{3}{2}},  $$
which yields 
$$  |P_{j,M} (\tilde \Theta)| +| \tilde P_{j,M} (\tilde \Theta)| \le \frac{C}{s^{\frac{3}{2}}} .$$
and
\begin{eqnarray}
& &  -i\kappa - \kappa (\delta-i)\frac{y^2}{\sqrt s}\frac{b}{(p-1)^2}  +a(\delta-i)\frac{1}{\sqrt s} 
 + \kappa (1-i\beta)\delta (p+1)\frac{y^4}{ s}\frac{b^2}{2(p-1)^4}\\
 & =&  \left( -\kappa + \frac{\Theta^*_{0,0}}{\sqrt{s}}\right) h_0 + \frac{\tilde \Theta^*_{0,0}}{\sqrt{s}} \tilde h_0 + \frac{\Theta^*_{2,0}}{\sqrt{s}} h_2 + \frac{\tilde \Theta^*_{2,0}}{\sqrt{s}}   \tilde h_2 + \frac{\tilde \Theta^*_{2,1}}{s}   \tilde h_2 .
\label{Def-Theta-ij}
\end{eqnarray}
In addition to that, we can calculate these constants and we obtain
\begin{eqnarray*}
   \Theta^*_{0,0} &=& 4(1+\delta^2) \delta \beta\frac{\kappa b}{(p-1)^2},\\
   \tilde \Theta^*_{0,0} &=& -\beta (1+\delta^2) \frac{\kappa b}{(p-1)^2},\\
   \tilde \Theta^*_{2,0} &=& - \delta\frac{\kappa b}{(p-1)^2},\\
    \Theta^*_{2,0} &=& (1+\delta^2)\frac{\kappa b}{(p-1)^2},\\
     \tilde \Theta^*_{2,1} &=& -3\delta(p+1)(-\beta^2+\beta\delta-2)\frac{\kappa b^2}{(p-1)^4}.
  \end{eqnarray*}
\section{Formal derivation of constant $b$ and $\mu$ }\label{expand-derivation-P-}
We begin by recalling that we need the following equation, for the determination of $b$ 

\begin{equation}\label{equa-R-2}
\begin{array}{l}
- \frac{1}{2} R_2' r - \frac{R_2}{p-1} + p |R_0|^{p-1}R_2 +  R_1'' - 2 R_0 \varphi_0' \varphi_1' - R_1 \varphi_0'^2 - 2 \beta R_0' \varphi_1' - 2\beta R_1' \varphi_0' \\
- \beta R_0 \varphi_1'' - \beta R_1 \varphi_0''
+ \frac{1}{4}  R_0' r + \frac{p(p-1)}{2}|R_0|^{p-3} R_0 R_1^2 = 0
\end{array}
\end{equation}

Let us introduce the following
\begin{eqnarray*}
R_0  &=& (p-1 + b r^2)^{- \frac{1}{p-1}},\\
R_1 &=& (A_1 + B_1 r^2 ) (p-1 + br^2)^{- \frac{p}{p-1}},\\
\varphi_0' &= & \tilde A_0 r (p-1 + b r^2)^{-1},\\
\varphi_0'' &=& \tilde A_0 (p-1 + br^2)^{-1} - 2 \tilde A_0 b r^2 (p-1 + br^2)^{-2}\\
\varphi_1' &= & \tilde A_1 r (p-1 + b r^2)^{-1}+ \tilde B_1 r (p-1 + b r^2)^{-2}. 
\end{eqnarray*}
This implies that
\begin{eqnarray*}
R_1' &=& 2 B_1 r (p-1 + b r^2)^{- \frac{p}{p-1}}  - \frac{2pA_1 b r}{p-1} (p-1 + b r^2)^{- \frac{2p -1}{p-1}}\\
&-  & \frac{2 p B_1 b r^3}{p-1} (p-1 + b r^2)^{- \frac{(2p -1)}{p-1}},\\
R_1'' &=& 2 B_1 (p-1 + br^2)^{- \frac{p}{p-1}} - \frac{8 p B_1b r^2}{p-1} (p-1 + br^2)^{- \frac{2p-1}{p-1}} \\
&+& (A_1 + B_1 r^2) \left( - \frac{2 p b }{p-1} (p-1 + br^2)^{- \frac{(2p-1)}{p-1}} + \frac{4p (2p -1)b^2 r^2}{(p-1)^2} (p-1 + br^2)^{- \frac{(3p -2)}{p-1}}\right) \\
\varphi_1'' &=& \tilde A_1 (p-1  + b r^2)^{-1}  + \tilde B_1 (p-1 + br^2)^{-2} \\
&-& 2 \tilde A_1 b r^2 (p-1 + b r^2)^{-2}   - 4 \tilde B_1 b r^2 (p-1 + br^2)^{-3}.
\end{eqnarray*}
Next, we have to compute
\begin{eqnarray*}
F_3 &=& R_1'' - 2 R_0 \varphi_0' \varphi_1' - R_1 \varphi_0'^2  - 2 \beta R_0' \varphi_1' - 2 \beta R_1' \varphi_0' - \beta R_0 \varphi_1'' - \beta R_1 \varphi_0'' \\
 &+& \frac{1}{4}R_0' r  + \frac{p(p-1)}{2}|R_0|^{p-2} R_1^2, 
\end{eqnarray*}
where
\begin{eqnarray*}
R_1'' & = &   2 B_1 (p-1 + br^2)^{- \frac{p}{p-1}} - \frac{8 p B_1b r^2}{p-1} (p-1 + br^2)^{- \frac{2p-1}{p-1}} \\
&+& (A_1 + B_1 r^2) \left( - \frac{2 p b }{p-1} (p-1 + br^2)^{- \frac{(2p-1)}{p-1}} + \frac{4p (2p -1)b^2 r^2}{(p-1)^2} (p-1 + br^2)^{- \frac{(3p -2)}{p-1}}\right) \\
 2 R_0 \varphi_0' \varphi_1' &=&   2 \tilde A_0 r (p-1 + b r^2)^{-\frac{p}{p-1}} \left(  \tilde A_1 r (p-1 + b r^2)^{-1} + \tilde B_1r (p-1 + br^2)^{-2} \right)\\
&=& 2 \tilde A_0 \tilde A_1 r^2 (p-1 + br^2)^{- \frac{(2p-1)}{p-1}} + 2 \tilde A_0 \tilde B_1 r^2 (p-1 + br^2)^{- \frac{(3p -2)}{p-1}}. \\
R_1 \varphi_0'^2 &=& (A_1 + B_1 r^2) (p-1 + br^2)^{- \frac{p}{p-1}} (\tilde A_0 r(p-1 + br^2)^{-1})^2 \\
 &= & (\tilde A_0^2 A_1 r^2 + \tilde A_0^2 B_1 r^4) (p-1 + br^2)^{- \frac{(3p -2)}{p-1}}.\\
 2 \beta R_0' \varphi_1' & = & 2 \beta \left( - \frac{2 b r }{p-1} (p-1 + br^2)^{- \frac{p}{p-1}}\right)  (\tilde A_1 r (p-1 + b r^2)^{-1}+ \tilde B_1 r (p-1 + b r^2)^{-2})\\
 &= & - \frac{4 \beta \tilde A_1 b r^2}{p-1} (p-1 + br^2)^{- \frac{(2p-1)}{p-1}}  - \frac{4 \beta \tilde B_1 b r^2}{p-1} (p-1 + b r^2)^{ - \frac{(3p-2)}{p-1}}.
\end{eqnarray*}

\begin{eqnarray*}
2 \beta R_1' \varphi_0' &=& 2 \beta \left(  2 B_1 r (p-1 + b r^2)^{- \frac{p}{p-1}} - \frac{2 p A_1 b r}{p-1} (p-1 + b r^2)^{- \frac{(2p-1)}{p-1 }}  - \frac{2 p B_1 b r^3}{p-1} (p-1 + b r^2)^{- \frac{(2p - 1)}{p-1}}  \right)\\
&\times & ( \tilde A_0 r (p-1 + br^2)^{-1})\\
&=& 4 \beta \tilde A_0 B_1 r^2 (p-1 + b r^2)^{- \frac{(2p - 1)}{p-1}} - \frac{4p \beta \tilde A_0 A_1 b r^2}{p-1} (p-1 + b r^2)^{- \frac{3p - 2}{p-1}} \\
&-& \frac{4 p \beta \tilde A_0 B_1 b r^4}{p-1} (p-1 + b r^2)^{- \frac{3p - 2}{p-1}}\\
\beta R_0 \varphi_1'' &=& \beta (p -1 + br^2)^{- \frac{1}{p-1}} \left(  \tilde A_1 (p-1 + b r^2)^{-1} + \tilde B_1 (p-1 + br^2)^{-2} \right)\\
&+ &\beta (p -1 + br^2)^{- \frac{1}{p-1}} \left( - 2 \tilde A_1 b r^2 (p-1 + b r^2)^{- \frac{(2p-1)}{p-1}}  - 4 \tilde B_1 b r^2 (p-1 + b r^2)^{- \frac{3p - 2}{p-1}}\right)\\
&=& \tilde A_1 \beta  (p-1 + b r^2)^{-\frac{p}{p-1}} + \tilde B_1  \beta (p-1 + br^2)^{-\frac{(2p-1)}{p-1}} \\
&-& 2 \beta \tilde A_1 b r^2 (p-1 + b r^2)^{- \frac{(2p-1)}{p-1}} -  4 \beta \tilde B_1 b r^2 (p-1 + b r^2)^{- \frac{(3p-2)}{p-1}}  \\
\beta R_1 \varphi_0''  & = &  \beta (A_1  + B_1 r^2) (p-1 + br^2)^{- \frac{p}{p-1}} \left( \tilde A_0 (p-1 + br^2)^{-1} - 2 \tilde A_0 b r^2 (p-1 + br^2)^{-2} \right)\\
&=&  ( \beta \tilde A_0 A_1  + \beta \tilde A_0 B_1 r^2 ) (p-1 + br^2)^{- \frac{(2p -1)}{p-1}} - (2 \beta \tilde A_0 A_1 b r^2 + 2 \beta \tilde A_0 B_1 b r^4) (p-1 + b r^2)^{- \frac{(3p-2)}{p-1}}\\
\frac{1}{4} R_0' r & = &- \frac{br^2}{2(p-1)} (p-1+ br^2)^{- \frac{p}{p-1}}. 
\end{eqnarray*}
\begin{eqnarray*}
\frac{p(p-1)}{2} |R_0|^{p-3}R_0 R_1^2 &= & \frac{p(p-1)}{2}(A_1 + B_1r^2)^2 (p-1 + br^2)^{-\frac{(3p -2)}{p-1}}.
\end{eqnarray*}
In the following, we decompse $F_3$  as: 
$$F_3 =  Q_1 (r)  + Q_2 (r) + Q_3(r),$$
Where
\begin{eqnarray*}
Q_1  &=& (p-1 + br^2)^{- \frac{p}{p-1}} \left(2B_1  - \tilde A_1 \beta - \frac{br^2}{2(p-1)} \right) \\
& = & Q_{1,1} (p-1 + br^2)^{- \frac{p}{p-1}}  + Q_{1,2} r^2 (p-1 + br^2)^{- \frac{p}{p-1}}.\\
Q_2 &= & (p-1 + br^2)^{- \frac{(2p -1)}{p-1}}\left( - \frac{2 p b A_1}{p-1}  - \beta \tilde B_1 - \beta \tilde A_0 A_1  \right) \\
 &+& (p-1 + b r^2)^{- \frac{(2p-1)}{p-1}} \left( -  \frac{10pB_1b }{p-1} - 2 \tilde A_0 \tilde A_1  + \frac{4 \beta \tilde A_1 b}{p-1}- 4 \beta \tilde A_0 B_1 + 2 \beta \tilde A_1b - \beta \tilde A_0 B_1 \right)r^2\\
 &=& Q_{2,1} (p-1 + b r^2)^{- \frac{(2p-1)}{p-1}} + Q_{2,2} r^2 (p-1 + b r^2)^{- \frac{(2p-1)}{p-1}}.\\
Q_3 &=& (p-1 + b r^2)^{- \frac{3p - 2}{p-1}} \left(  \frac{p (p-1)}{2} A_1^2 \right)\\
& + &  (p-1 + b r^2)^{- \frac{3p - 2}{p-1}}  r^2  \left(  \frac{4p (2p-1) A_1 b^2}{(p-1)^2}  - 2 \tilde A_0 \tilde B_1  - \tilde A_0^2 A_1 \right)\\
&+& (p-1 + b r^2)^{- \frac{3p - 2}{p-1}}  r^2\left( \frac{4 p \beta \tilde B_1 b}{p-1} + \frac{2 (3p-1) \beta \tilde A_0 A_1 b}{p-1}   + p (p-1) A_1 B_1 \right)\\
&+& (p-1 + b r^2)^{- \frac{3p - 2}{p-1}}  r^4  \\
& = & Q_{3,1} (p-1 + b r^2)^{- \frac{3p - 2}{p-1}} + Q_{3,2} r^2 (p-1 + b r^2)^{- \frac{3p - 2}{p-1}} + Q_{3,3} r^4 (p-1 + b r^2)^{- \frac{3p - 2}{p-1}}. 
\end{eqnarray*}
By using variation of constant, we deduce that
\begin{eqnarray*}
R_2 = H^{-1}(r) \left(  \displaystyle \int F_3 \frac{2 H}{r}  H dr \right),
\end{eqnarray*}
where
$$H(r)  = \frac{(p-1 + b r^2)^{ \frac{p}{p-1}}}{r^2}.$$
We see that  the 

\begin{eqnarray*}
\frac{2 H }{r} F_3= 2\left( Q_{1,2}   - \frac{ b Q_{2,1}}{(p -1)^2} + \frac{ Q_{2,2}}{p-1} - \frac{  2b Q_{3,1}}{(p-1)^3}  + \frac{Q_{3,2}}{ (p-1)^2} \right) \frac{1}{r} + \text{  integrable term }.
\end{eqnarray*}
 As a matter of fact, we always find  $R_2$ such that $\Delta R_2$ is continuous at $0$.  This implies that $R_2$ may not contain terms $\log |r|$, consequently, we obtain,
 \begin{eqnarray}
P = Q_{1,2}   - \frac{ b Q_{2,1}}{(p -1)^2} + \frac{ Q_{2,2}}{p-1} - \frac{ 2 b Q_{3,1}}{(p-1)^3}  + \frac{Q_{3,2}}{ (p-1)^2} = 0.\label{conditio,-b}
\end{eqnarray}
We recall that
\begin{eqnarray*}
A_1  &=& - \frac{2 b (\delta \beta - 1 )}{p-1} \text{ and } B_1 =  \mathcal{C} ,\\
\tilde A_0 &=& - \frac{2 \delta b}{p-1},\\
\tilde A_1 &= &  \frac{4 \beta ( 1 + \delta^2) b^2}{(p-1)^2} \text{ and } \tilde B_1 = \frac{4 b^2}{(p-1)^2} \left(  (p +  3) \delta + \beta(2 p + \delta^2(p-3)) \right ) + 2 (p-1) \delta \mathcal{C}.
\end{eqnarray*}
We also define  $\mathcal{C}_1 = 2 (p-1) \delta \mathcal{C}$.  Then, 
$$ \tilde B_1 = \frac{4 b^2}{(p-1)^2} \left(  (p +  3) \delta + \beta(2 p + \delta^2(p-3)) \right ) + \mathcal{C}_1 $$
\begin{eqnarray*}
Q_{1,2} &=& - \frac{b}{2 (p-1)},\\
Q_{2,1} &=&  - \frac{2 p b A_1}{p-1} - \beta \tilde B_1 - \beta \tilde A_0 A_1\\
&=&   - \frac{2 p }{p-1} b  \left( - \frac{2 b (\delta \beta -1)}{p-1}  \right) - \beta \left( \frac{4 b^2}{(p-1)^2} \left(  (p +  3) \delta + \beta(2 p + \delta^2(p-3)) \right ) + \mathcal{C}_1 \right)\\
& - & \beta \left(  - \frac{2 \delta b}{p-1}    \right) \left( - \frac{2 b (\delta \beta -1 )}{ p-1} \right)\\
&=&  \frac{4 p(\delta \beta -1) b^2}{ (p-1)^2} - \frac{4 \delta \beta (\delta \beta -1)b^2}{ (p-1)^2}- \frac{4 b^2}{ (p-1)^2} \left(  (p+3) \delta \beta  + \beta^2 ( 2p + (p-3) \delta^2)  \right) - \beta \mathcal{C}_1\\
& = &  \frac{4 b^2}{(p-1)^2}  \left(   p (\delta \beta - 1)  - \delta \beta (\delta \beta -1)  - (p+3) \delta \beta - \beta^2 ( 2 p + (p-3) \delta^2)   \right) - \beta \mathcal{C}_1 \\
& = &  \frac{4 b^2}{(p-1)^2} \left(  (2 - p) \delta^2 \beta^2 - 2 \delta \beta - p - 2 p \beta^2  \right) - 2 (p-1) \delta \beta \mathcal{C}
\end{eqnarray*}
and 
\begin{eqnarray*}
Q_{2,2} &=& - \frac{10p B_1 b}{p-1}  - 2 \tilde A_0 \tilde A_1  + \frac{4 \beta \tilde A_1 b}{p-1}- 4 \beta \tilde A_0 B_1 + 2 \beta \tilde A_1b - \beta \tilde A_0 B_1 \\
&=& - \frac{10 p b B_1 }{p-1}  + \frac{2 (p+1) \beta  \tilde A_1 b}{p-1} - 5 \beta \tilde A_0 B_1 - 2 \tilde A_0 \tilde A_1\\
&=&    - \frac{10 p b}{p-1} \mathcal{C} + \frac{2 (p+1)  }{p-1} \beta \left(  \frac{4 \beta (1 + \delta^2) b^2}{(p-1)^2}\right)b - 5 \beta \left( -\frac{2 \delta b}{p-1} \right) \mathcal{C} \\
&-& 2 \left( -\frac{2 \delta b}{p-1}\right) \left( \frac{4 \beta (1 + \delta^2) b^2}{(p-1)^2} \right)\\
& =& \frac{b^3}{(p-1)^3} \left(  8 (p+1) \beta^2 (1+\delta^2)  + 16 \delta \beta (1 +\delta^2)\right)  +  \left(   \frac{10 \delta \beta b}{p-1}-  \frac{10 p b}{p-1}  \right) \mathcal{C}\\
& = &  \frac{b^3}{ (p-1)^3}  \left(   16 \delta^3 \beta + 8(p+1) \delta^2 \beta^2  + 8 (p+1) \beta^2 + 16 \delta \beta \right)  +   \left(   \frac{10 \delta \beta b}{p-1}-  \frac{10 p b}{p-1}  \right) \mathcal{C}\\
Q_{3,1} &=& \frac{2 p (\delta \beta - 1 )^2 b^2}{p-1} = \frac{b^2}{p-1} ( 2 p (\delta \beta -1)^2), 
\end{eqnarray*}
and 
\begin{eqnarray*}
Q_{3,2} &=& \frac{4 p (2p -1) A_1 b^2}{(p-1)^2} - 2 \tilde A_0 \tilde B_1 - \tilde A_0^2 A_1 + \frac{4 p \beta \tilde B_1 b}{p-1} + \frac{2 (3p - 1) \beta  \tilde A_0 A_1 b}{p-1} + p (p-1) A_1 B_1 \\
&=& \frac{4 p (2p -1)}{(p-1)^2} \left( -\frac{2 b (\delta \beta -1) }{p-1} \right) b^2 - 2 \left( -\frac{2 \delta b }{p-1} \right) \left\{ \frac{4b^2}{(p-1)^2}( (p+3)\delta + \beta ( 2p + (p-3) \delta) ) + \mathcal{C}_1\right\}\\
& - &  \left(  - \frac{2 \delta b}{p-1}\right)^{2} \left(  - \frac{2 b (\delta \beta -1)}{p-1}  \right)+ \frac{4 p \beta b }{p-1} \left\{ \frac{4 b^2}{(p-1)^2} \left( (p+3) \delta + \beta (2p  + (p-3) \delta^2) \right) + \mathcal{C}_1\right\}   \\
&+ &\frac{2 (3p - 1) \beta }{p-1} \left( - \frac{2 \delta b}{p-1}\right) \left( - \frac{2 b (\delta \beta -1)}{p-1}\right) b + p (p-1) \left( - \frac{2 b (\delta \beta -1)}{ p-1} \right) \mathcal{C}\\
& =& - \frac{8 p (2p-1)(\delta \beta -1) b^3}{(p-1)^3} + \frac{16 \delta b^3}{(p-1)^3} \left(  (p+3) \delta + \beta (2p + (p-3) \delta^2) \right) + \frac{4  \delta b}{p-1} \mathcal{C}_1 + \frac{8 \delta^2 (\delta \beta -1) b^3}{(p-1)^3}\\
& +& \frac{16 p \beta b^3}{(p-1)^3} \left(  (p+3) \delta + \beta (2p + (p-3) \delta^2) \right) + \frac{4 p \beta b}{p-1} \mathcal{C}_1 + \frac{8 (3p-1) \delta \beta (\delta \beta -1) b^3}{(p-1)^3} - 2 p(\delta \beta -1) b \mathcal{C}\\
& =& \frac{b^3}{(p-1)^3} \left[ - 8 p (2p-1)(\delta \beta -1) + 16(p+3) \delta^2 + 16 \delta \beta (2p  + (p-3) \delta^2) + 8 \delta^2 (\delta \beta -1)   \right.\\
& + & \left. 16 p (p+3) \delta \beta  + 16 p \beta^2 (2p + (p-3) \delta^2) + 8 (3p-1) \delta \beta (\delta \beta -1)\right] \\
&+& \frac{4 \delta b}{p-1} \mathcal{C}_1 + \frac{4 p \beta b}{p-1} \mathcal{C}_1 - 2p (\delta \beta -1) b\mathcal{C}\\
&=& \frac{b^3}{(p-1)^3}  \left[   \delta^3 \beta (16 p - 40) + \delta^2 \beta^2 (16p^2 - 24 p - 8) + \delta^2 ( 16 p + 40)\right. \\
&+& \left.  32 p^2 \beta^2 + \delta \beta (64 p + 8) + 8p (2p-1) \right] \\
&  +& b (8 \delta^2 + 6 p \delta \beta + 2 p ) \mathcal{C}
\end{eqnarray*}
Finally, we obtain
\begin{eqnarray*}
P &=& Q_{1,2}  - \frac{b Q_{2,1}}{(p-1)^2} + \frac{Q_{2,2}}{p-1}  - \frac{2 b Q_{3,1}}{(p-1)^3} + \frac{Q_{3,2}}{(p-1)^2} \\
&=&  - \frac{b}{2 (p-1)}  - \frac{b}{(p-1)^2} \left\{   \frac{4 b^2}{ (p-1)^2} [  (2-p) \delta^2 \beta^2 - 2 p \beta^2 - 2 \delta \beta - p  ]- 2 (p-1) \delta \beta \mathcal{C}  \right\}\\
&+& \frac{1}{p-1} \left\{    \frac{b^3}{(p-1)^3} [  16 \delta^3 \beta + 8 (p+1) \delta^2 \beta^2 + 8 (p+1) \beta^2 + 16 \delta \beta ] + \left(\frac{  10 \delta \beta b - 10 pb}{ p-1}  \right)   \mathcal{C}   \right\}\\
& - & \frac{2 b}{(p-1)^3} \left\{    \frac{b^2}{ p-1}   ( 2p \delta^2 \beta^2 - 4p \delta \beta + 2p) \right\}\\
& +& \frac{1}{ (p-1)^2} \left\{  [ (16p -40) \delta^3 \beta + (16p^2 - 24 p -8) \delta^2 \beta^2 + (16p + 40) \delta^2    + 32 p^2 \beta^2 + (64 p + 8) \delta \beta  +16p^2 - 8p ]    \right.\\
&\times & \left.  \frac{b^3}{ (p-1)^3} +    b ( 8 \delta^2 + 6p \delta \beta + 2p) \mathcal{C} \right\} 
\end{eqnarray*}
We conclude that $P=0$, gives us the constant $b$ as follows;
\begin{eqnarray*}
  \frac{b}{2( p-1)}  &=&    \frac{b^3}{ (p-1)^5 } \left\{    - 4 (p-1)  [  (2-p) \delta^2 \beta^2 - 2 p \beta^2 - 2 \delta \beta - p  ] \right.  \\
  &  +&     (p-1) [  16 \delta^3 \beta + 8 (p+1) \delta^2 \beta^2 + 8 (p+1) \beta^2 + 16 \delta \beta ]\\
  & - & 2 (p-1) ( 2p \delta^2 \beta^2 - 4p \delta \beta + 2p)\\
  & +& \left.  [ (16p -40) \delta^3 \beta + (16p^2 - 24 p -8) \delta^2 \beta^2 + (16p + 40) \delta^2    + 32 p^2 \beta^2 + (64 p + 8) \delta \beta  +16p^2 - 8p ]\right\} \\
  & +& \frac{ \mathcal{C}}{(p-1)^2}  \left(   10 \delta \beta b- 10 p b + 2 (p-1) \delta \beta b + b(8 \delta^2 + 6p \delta \beta  + 2 p)   \right) \\
  & = & \frac{b^3}{ (p-1)^5} \left[  \delta^3 \beta (32p -56) +  \delta^2 \beta^2 (24 p^2   -32 p - 8)  + \beta^2 (48p^2  -8p -8) + \delta^2 (16 p + 40) \right.\\
  & +&  \left.  \delta \beta (8p^2 + 80 p - 16) + 16 p^2 - 8p   \right].
.\end{eqnarray*}
By the critical condition we can write $  \beta  = \frac{p- \delta^2}{ (p+1) \delta}$ and we obtain 
\begin{eqnarray*}
\frac{b}{2( p-1)}  &=&  \frac{b^5}{ (p-1)^5}  \left[   \delta^3 \frac{p- \delta^2}{ (p+1) \delta} (32p -56)   + \delta^2 \left(  \frac{p- \delta^2}{ (p+1) \delta}\right)^2  (24 p^2   -32 p - 8) \right.\\
&  +& \left. \left(  \frac{p- \delta^2}{ (p+1) \delta}\right)^2     (48p^2  -8p -8) \right.\\
& + & \left.  \delta^2 (16 p + 40)  + \delta   \frac{p- \delta^2}{ (p+1) \delta} (8p^2 + 80 p - 16 ) + 16p^2 -8p  \right]\\
& =  &  \frac{b^5}{ (p+1)^2 \delta^2 (p-1)^5} \left[ \delta^6 ( - 8p^2 - 8p +48 )  + \delta^4 ( -8p^3 + 72 p^2 -16 p + 48 ) \right.\\
&+& \left.  \delta^2 (  48 p^4 - 16p^3 + 72p^2 - 8p   )  + p^2 ( 48 p^2 - 8p - 8) \right].
\end{eqnarray*}

\medskip

\textbf{From the determination of $\mu$}

\medskip

The equation $F_4(0)=0$ gives us the value of $\mu$ as follows
which yields
\begin{eqnarray*}
\mu &=&   \varphi_1'' (0)  - 2 \beta \varphi_0' (0) \varphi_1'(0) + 2 R_0^{-1} (0) R_0' (0) \varphi_1' (0) + 2 R_0^{-1}(0) R_1'(0) \varphi_0' (0)\\
& - &  2 R_0^{-2} (0) R_0' (0) \varphi_0' (0) R_1 (0) +   \beta R_0^{-1} (0) R_1'' (0) - \beta R_0^{-2}R_0'' (0) R_1 (0)\\
& + & \delta (p-1) R_0^{p-2}  (0) R_2 (0) + \frac{\delta (p-1) (p-2)}{2} R_0^{p-3}(0) R_1^2 (0). 
\end{eqnarray*}

+  The calculation of   $\varphi_1'' (0):$
\begin{eqnarray*}
\varphi_1'' (0) &=&    \frac{4 \beta (1 + \delta^2) b^2}{ (p-1)^3} + \frac{4b^2}{(p-1)^4} \left\{ (p+3)\delta + \beta [ 2p + \delta^2 (p-3)] \right\} + \frac{2 \delta \mathcal{C}}{p-1}\\
& =& \frac{4b^2}{ (p-1)^4} \left\{  (p+3) \delta + \beta [ 3p-1 + \delta^2( 2p -4) ]  \right\} + \frac{2 \delta \mathcal{C}}{p-1}.
\end{eqnarray*}

+ The calculation of  $- 2 \beta \varphi_0' (0) \varphi_1'(0)$:
\begin{eqnarray*}
- 2 \beta \varphi_0' (0) \varphi_1'(0) & =  &  0, \text{ because of } \varphi_0'(0) =0. 
\end{eqnarray*}

+ The calculation of  $2 R_0^{-1} (0) R_0' (0) \varphi_1' (0) $:

\begin{eqnarray*}
2 R_0^{-1} (0) R_0' (0) \varphi_1' (0) & =  &  0, \text{ because of } R_0'(0) =0. 
\end{eqnarray*}

+ The calculation $  2 R_0^{-1}(0) R_1'(0) \varphi_0' (0)$:
\begin{eqnarray*}
 2 R_0^{-1}(0) R_1'(0) \varphi_0' (0) = 0  \text{ because of } \varphi_0'(0) =0.
\end{eqnarray*}

+ The calculation of $-2 R_0^{-2} (0) R_0' (0) \varphi_0' (0) R_1 (0)$
\begin{eqnarray*}
-2 R_0^{-2} (0) R_0' (0) \varphi_0' (0) R_1 (0) = 0 \text{ because of }   \varphi_0' (0) = 0.
\end{eqnarray*}

+ The calculation  of   $\beta R_0^{-1} (0) R_1'' (0)$

\begin{eqnarray*}
\beta R_0^{-1} (0) R_1'' (0) & = &  \beta (p-1)^{\frac{1}{p-1}} \left[  2 \mathcal{C} (p-1)^{-\frac{p}{p-1}}   +  \frac{4 p (\delta \beta -1)b^2}{(p-1)^4} (p-1)^{-\frac{1}{p-1}}\right] \\
& = & \frac{4 p \beta(\delta \beta -1)b^2}{(p-1)^4} + \frac{2 \beta \mathcal{C}}{p-1}. 
\end{eqnarray*}

+ The calculation  of $- \beta R_0^{-2}R_0'' (0) R_1 (0)$:
\begin{eqnarray*}
- \beta R_0^{-2} (0) R_0'' (0) R_1 (0) &=& - \beta (p-1)^{\frac{2}{p-1}} \left( -\frac{2b}{p-1} \right)(p-1)^{-\frac{p}{p-1}} \left( - \frac{(\delta \beta -1)2b}{ p-1} \right) (p-1)^{-\frac{p}{p-1}}\\
& = &  - \frac{4 \beta (\delta \beta -1) b^2}{ (p-1)^4}. 
\end{eqnarray*}

+ The calculation of   $\delta (p-1) R_0^{p-2}  (0) R_2 (0)$: The calculation  is more  difficult, we first determine  $ R_2 (0)$ .We derive  from \eqref{equa-R-2} that
\begin{equation}\label{equa-R-2-0}
R_2 (0) = - F_3 (0),
\end{equation}
where  $F_3 $  is defined in \eqref{def-of-F_3}.  We now find  $F_3 (0)$
\begin{eqnarray*}
F_3(0) & = &  R_1'' (0)  - 2 R_0 (0) \varphi_0' (0) - R_1 (0) [ \varphi_0' (0)]^2 -2 \beta R_0' (0) \varphi_1' (0) - 2 \beta R_1' (0) \varphi (0) - \beta R_0(0) \varphi_1'' (0)\\
& -& \beta R_1 (0) \varphi_0'' (0) + \frac{p (p-1)}{2} R^{p-2}_0 (0) R_1^2 (0)\\
& =& R_1'' (0)   - \beta R_0(0) \varphi_1'' (0) - \beta R_1 (0) \varphi_0'' (0) + \frac{p (p-1)}{2} R^{p-2}_0 (0) R_1^2 (0).
\end{eqnarray*}

We have 
\begin{center}
\begin{eqnarray*}
& & R_1'' (0)  =  \frac{4 p (\beta \delta -1)b^2}{ (p-1)^4}  (p-1)^{- \frac{1}{p-1}}  + \frac{2 \mathcal{C} }{p-1} (p-1)^{-\frac{1}{p-1}},\\
& &  - \beta R_0(0) \varphi_1'' (0)  =  - \left[  \tilde A_1 \beta (p-1)^{-\frac{p}{p-1}}  + \tilde B_1 \beta (p-1)^{-\frac{2p-1}{p-1}} \right]\\
 & = & - \left[ \frac{4 \beta^2 (1 + \delta^2) b^2}{(p-1)^3}  + \frac{4b^2}{(p-1)^4} ((p+3) \delta \beta + \beta^2 [ 2p + \delta^2 (p-3)])  + \frac{2 \delta \beta \mathcal{C}}{p-1} \right] (p-1)^{-\frac{1}{p-1}},\\
 & &- \beta R_1 (0) \varphi_0'' (0) =  - \left[  \beta \tilde A_0 A_1 (p-1)^{-\frac{2p-1}{p-1}}  \right]\\
 & =& -  \left[  \frac{4 \delta \beta  (\delta \beta -1)b^2}{(p-1)^4}\right] (p-1)^{-\frac{1}{p-1}},\\
 & & \frac{p (p-1)}{2} R_0^{p-2} (0) R_1^2 (0)  = \frac{p (p-1)}{2} A_1^2 (p-1)^{- \frac{3p - 2}{p-1}}\\
 & = & \frac{2 p (\delta \beta -1)^2}{(p-1)^4} (p-1)^{-\frac{1}{p-1}}  .
\end{eqnarray*}
\end{center}
Then,  we obtain the following 
\begin{eqnarray*}
& & F_3 (0)  =     R_1'' (0) - \beta R_0 (0) \varphi''_1 (0) - \beta R_1 (0) \varphi_0''(0) + \frac{p (p-1)}{2} R_0^{p-2} (0) R_1^2 (0) \\
&=& \frac{4p (\delta \beta -1)b^2}{(p-1)^4} (p-1)^{-\frac{1}{p-1}} + \frac{2 \mathcal{C}}{p-1} (p-1)^{-\frac{1}{p-1}}\\
& - & \frac{4 \beta^2 (1+\delta^2)b^2}{(p-1)^3} (p-1)^{-\frac{1}{p-1}} - \frac{4b^2}{(p-1)^4} ((p+3) \delta \beta  + \beta^2 [2P + \delta^2(p-3)]) (p-1)^{-\frac{1}{p-1}}\\
& - & \frac{2 \delta \beta \mathcal{C}}{p-1}(p-1)^{-\frac{1}{p-1}} - \frac{4 \delta \beta (\delta \beta-1) b^2}{(p-1)^4} (p-1)^{-\frac{1}{p-1}} + \frac{2 p(\delta \beta -1)^2 b^2 }{(p-1)^4} (p-1)^{-\frac{1}{p-1}}\\
& =&  \left[ \frac{2 \mathcal{C}}{p-1} - \frac{2 \delta \beta \mathcal{C}}{p-1}   \right] (p-1)^{-\frac{1}{p-1}}\\
&+& \left\{   4p (\delta \beta -1)  - 4 (p-1)\beta^2 (1+\delta^2) - 4 ((p+3) \delta \beta + \beta^2 [2p + \delta^2(p-3)]) -4 \delta \beta (\delta \beta -1) + 2p (\delta \beta-1)^2\right\} \\
& \times & \frac{b^2}{(p-1)^4} (p-1)^{-\frac{1}{p-1}}\\
& = & \left[ \frac{2 \mathcal{C}}{p-1} - \frac{2 \delta \beta \mathcal{C}}{p-1}   \right] (p-1)^{-\frac{1}{p-1}}\\
&+& \frac{b^2}{(p-1)^4} ( (\delta \beta -1) (2p + (2p -4) \delta \beta)  -4 (p+3) \delta \beta- \beta^2[ 12p -4 + \delta^2 (8p -16)] ) (p-1)^{-\frac{1}{p-1}}.
\end{eqnarray*}
We then get from equation $F_4(0)=0$:

\begin{eqnarray*}
R_2 (0) & =  & - F_3 (0) = \left[ \frac{2 \delta \beta \mathcal{C}}{p-1}   - \frac{2 \mathcal{C}}{p-1}  \right] (p-1)^{-\frac{1}{p-1}}  \\
& +& \frac{b^2}{(p-1)^4} \{  4 (p+3) \delta \beta +  \beta^2[ 12p -4 + \delta^2 (8p -16)]  -  (\delta \beta -1) (2p + (2p -4) \delta \beta)  \} (p-1)^{-\frac{1}{p-1}}
\end{eqnarray*}
Therefore, we get the following 
\begin{eqnarray*}
& & \delta (p-1) R_0^{p-2}  (0) R_2 (0) = \frac{2 \delta^2 \beta \mathcal{C}}{p-1}   - \frac{2\delta \mathcal{C}}{p-1} \\
& +& \frac{b^2}{(p-1)^4} \{  4 (p+3) \delta^2 \beta +  \delta \beta^2[ 12p -4 + \delta^2 (8p -16)]  -  (\delta \beta -1) (2p \delta + (2p -4) \delta^2 \beta)  \}.
\end{eqnarray*}

+ The calculation of $\frac{\delta (p-1) (p-2)}{2} R_0^{p-3}(0) R_1^2 (0)$. We have
\begin{eqnarray*}
\frac{\delta (p-1) (p-2)}{2} R_0^{p-3}(0) R_1^2 (0) &=& \frac{\delta (p-1) (p-2) }{2} (p-1)^{-\frac{p-3}{p-1}}  \left(  -\frac{2b (\delta \beta -1)}{p-1} (p-1)^{-\frac{p}{p-1}}    \right)^2\\
& =& \frac{2 (p-2)\delta (\delta \beta -1)^2 b^2 }{(p-1)^4}.
\end{eqnarray*}
Finally, we obtain 
\begin{eqnarray*}
& & \mu   =  \varphi_1''(0) + \beta R_0^{-1} R_1'' (0) - \beta R_0'' R_0^{-2} R_1(0) + \delta (p-1) R_0^{p-2}R_2(0) + \frac{\delta (p-1)(p-2)}{2}R_0^{p-3} R_1^2(0)\\
& = & \frac{b^2}{(p-1)^4} \left\{  4 (p+3) \delta + 4 (3p-1)\beta  + 4 (2p-4) \delta^2 \beta + 4p \beta^2 \delta - 4 p\beta - 4 \beta^2 \delta + 4 \beta  + 4 (p+3) \delta^2 \beta    \right.\\
& +& \left. (12 p-4) \delta \beta^2 + \delta^3 \beta^2 (8p-16)  - 2p \delta^2 \beta - (2p-4) \delta^3 \beta^2    +2 p\delta  - (2p-4) \delta^2 \beta + 2 (p-2) \delta^3 \beta^2   \right.\\
& -& 4(p-2)\delta^2 \beta + \left.  2 (p-2) \delta \right\} + \frac{2 \beta (1+\delta^2) \mathcal{C}}{p-1}\\
&  = & \frac{b^2}{(p-1)^4} \left\{   8 (p+1) \delta + 8p \beta  + (4p +8) \delta^2 \beta + (16p-8) \delta \beta^2 + (8p -16) \delta^3 \beta^2   \right\} + \frac{2 \beta (1 + \delta^2) \mathcal{C}}{p-1}.
\end{eqnarray*}

\section{Cancellation of some coefficient in the ODE  of $\tilde q_2$}\label{cancellation-ODE-tilde-q-2}

In this part, we aim at giving  some details of the computation some qualities   related to the construction $\tilde q_2$'s  ODE.

a) The cancellation of \eqref{equality-order-s-3-2-tilde-q-2}:  We now rewrite this one as follows
$$  \frac{\nu}{2} R_{2,1}^* + \tilde C_{2,2}R^*_{2,1} - \tilde D_{2,0} \tilde R_{0,1}^*  -  \frac{\delta b}{(p-1)^2} R^*_{0,1}  + \tilde R^*_{2,2}  = 0.  $$

Using  the definitions of the  constants  in  Appendices \ref{appendix-potential} and  \ref{appendix-expansion-R}  we can derive the right hand side  as follows

\begin{eqnarray*}
& &  \frac{\nu}{2} R_{2,1}^* + \tilde C_{2,2}R^*_{2,1} - \tilde D_{2,0} \tilde R_{0,1}^*  -  \frac{\delta b}{(p-1)^2} R^*_{0,1}  + \tilde R^*_{2,2} \\
 & = & - \frac{ \kappa b}{2 (p-1)^2}\\
 & -&  \frac{4 b^3 \kappa}{(p-1)^6} \left\{ -12\beta\delta p-9\beta^3\delta p+3\beta^3\delta^3 p\right.\\
 & & -3p\beta^2\delta^4-2p^2\beta\delta-12p^2 \beta^3\delta-5\delta^2\beta^2p^2-16\delta^2 p\beta^2+2\beta^2-5p^2 \\
 &+& \left. 4\delta^4-13\delta^2+5p+\delta \beta+8\beta\delta^3+3\beta^3\delta-5\delta^2\beta^2-3\delta^5\beta+3\beta^3\delta^3-p\beta^2-7p\delta^2   \right\}.
\end{eqnarray*} 

In addition to that, using the critical  condition 

$$  p - (p+1) \delta \beta - \delta^2 = 0, $$

then, we derive the following

\begin{eqnarray*}
0 & =&   \frac{\nu}{2} R_{2,1}^* + \tilde C_{2,2}R^*_{2,1} - \tilde D_{2,0} \tilde R_{0,1}^*  -  \frac{\delta b}{(p-1)^2} R^*_{0,1}  + \tilde R^*_{2,2} \\
  & = & - \frac{\kappa b }{ 2 (p-1)^2} \\
  & - &  \frac{ 8 b^3 \kappa}{(p-1)^6 (p+1)^2 \delta^2} (1 +\delta^2) (  (p^2 +p-6) \delta^4 +(p^3 -10 p^2 + p) \delta^2 +p^3-6p^4+p^2).
\end{eqnarray*}

This yields
\begin{eqnarray*}
b^2  &=&  \frac{(p-1)^4 (p+1)^2 \delta^2}{ - 16 (1 + \delta^2) L}, 
\end{eqnarray*}
where
\begin{eqnarray*}
L &=& (p^2 +p - 6) \delta^4  + p (p^2 - 10p + 1) \delta^2 +p^2 (- 6p^2 + p + 1).
\end{eqnarray*}

b/ Explain the decomposition      $\tilde\Ag_2(\tilde H_1+1)+\tilde H_2 $. Indeed, it is enough to prove that there no $\mu$ in $\tilde H_2$. Indeed, we recall formula of $\tilde H_2$ as follows: 
\begin{eqnarray*}
& &\tilde H_2 = \frac{\nu}{2}\left[  X_2 + c_4 [\tilde C_{4,2}R^*_{2,1} +\tilde R_{4,2}^*]  - D_{2,0} \tilde R_{0,1}\right]   - \frac{\nu}{2} \left[  \tilde K_{2,4} \left( \frac{C_{4,2} R^*_{2,1}}{2} + \frac{R^*_{4,2}}{2} \right)\right] \\
& -&  \frac{\nu}{2} \left[  \tilde L_{2,4} \left( \tilde C_{4,2} R^*_{2,1} + \tilde R_{4,2}^*  \right)\right] + \mu R^*_{2,1} + \frac{R^*_{0,1} R^*_{2,1}}{\kappa} \\ 
& + &   \tilde D_{2,0} \left[  -\tilde X_0 + \frac{ \nu \tilde K_{0,2} R_{2,1}^*}{2}       -  \tilde C_{0,2} R^*_{2,1}         \right]    + \tilde C_{2,2} \left[   X_2 + c_4 ( \tilde C_{4,2} R^*_{2,1} + \tilde R_{4,2})  - D_{2,0} \tilde R_{0,1}^*   \right]      \\
& +&  \tilde C_{2,4} \left[   \frac{C_{4,2} R^*_{2,1}}{2}  + \frac{R^*_{4,2}}{2}  \right]   + \tilde D_{2,4} (   \tilde C_{4,2} R^*_{2,1} + \tilde R^*_{4,2} ) +\tilde E_{2,2} R_{2,1}^* - \tilde F_{2,0} \tilde R_{0,1}^*  \\
& +& \frac{1}{8\kappa} (32 - 64 \delta \beta) (R_{2,1}^*)^2 \\
&-&   \frac{\delta b}{(p-1)^2} \left[ \frac{\nu}{2} (1 + \delta^2) \tilde R^*_{0,1} -\frac{\nu}{2} K_{0,2} R^*_{2,1} - D_{0,0} \tilde R^*_{0,1} + C_{0,2} R^*_{2,1}+R^*_{0,2} + \frac{\Theta^*_{0,0} R^*_{0,1}}{\kappa} \right] \\
& + &   \tilde R_{2,3}^*  +  \frac{    (p+1) \delta [ 12 - 6 \delta \beta + 6 \beta^2 ]   R^*_{0,1}    b^2}{2(p-1)^4}.
\end{eqnarray*}
Using the definitions of constants in $\tilde H_2$, we can write 
\begin{eqnarray*}
\tilde H_2  = & & \frac{\nu}{2}X_2   + \mu R^*_{2,1} + \frac{R^*_{0,1} R^{*}_{2,1} }{\kappa}  -  \tilde D_{2,0} \tilde X_0 +  \tilde C_{2,2} X_2 - \frac{\delta b}{(p-1)^2} R^*_{0,2} - \frac{\delta b}{(p-1)^2} \frac{\Theta^*_{0,0} R^*_{0,1}}{\kappa} \\
&+&  \tilde R^*_{2,3} + \frac{(p+1) \delta [ 12 - 6 \delta \beta + 6 \beta^2] R^*_{0,1} b^2}{2(p-1)^4} + B ( p, \delta, \beta),
 \end{eqnarray*}
where  $B$ doesn't contain $\mu$, and the first terms contain $\mu$.  We repeat the process of removing terms that do not contain, and we can get 
\begin{eqnarray*}
\tilde H_2  &= & \frac{\nu}{2} \left[  R_{2,2}^*  - \mu \Theta^*_{2,0} \right] - \tilde D_{2,0} \left[   \tilde R^*_{0,2} - \mu \tilde \Theta^*_{0,0} \right] + \tilde C_{2,2} \left[   R^*_{2,2} - \mu \Theta_{2,0}^* \right] - \frac{\delta b}{(p-1)^2} R^*_{0,2} \\
&+&  \frac{\delta b}{(p-1)^2} \mu \Theta^*_{0,0} + \tilde R^*_{2,3} - \frac{\mu \kappa b^2 (p+1) \delta [ 12 - 6 \delta \beta +6 \beta^2]}{2(p-1)^2} + B_1,
\end{eqnarray*}
where  $B_1$  doesn't contain $\mu$. In addition to that,  we have 
\begin{eqnarray*}
R^*_{2,2}  &=&  =  \mu \Theta_{2,0}^* + \text{`term no } \mu\text{'},  \\
 \tilde R_{0,2}^*  &=& \mu  \tilde \Theta_{0,0}^*  +  \text{`term no } \mu\text{'}, \\
R^*_{2,2}  &=&  =  \mu \Theta_{2,0}^* +\text{`term no } \mu\text{'},\\
 R_{0,2}^* & = &     \mu \Theta^*_{0,0} +  \text{`term no } \mu\text{'}, \\
 R_{2,3}^* & = &  \mu \tilde \Theta^*_{2,1} +  \text{`term no } \mu\text{'}.
\end{eqnarray*}

\noindent
This implies the fact that  $  \tilde H_2   $ doesn't contain $\mu$. Besides that, for the critical condition,  we have $\tilde H_1 \leq  -\frac{3}{2}$. Thus, we can write 
$$  \tilde\Ag_2(\tilde H_1+1)+\tilde H_2 =  a_0 (p,\delta) \mu + a_1(p, \delta, \beta), $$
with $a_0 \ne 0$.

\newpage
\bibliographystyle{alpha}

\end{document}